\documentclass[times,final,11pt]{elsarticle}

\usepackage{amssymb}
\usepackage{times}
\usepackage{soul}
\usepackage{url}
\usepackage{tikz}

\usepackage{fullpage}
\usepackage{graphicx} 
\graphicspath{{Figures/}} 
\usepackage[shortlabels]{enumitem} 

\usepackage{lipsum} 

\usepackage{subfig} 

\usepackage{amsmath,amssymb,amsthm,mathtools} 

\usepackage{varioref} 

\usepackage{multirow}

\usepackage[colorlinks=true, urlcolor=blue, linkcolor=red]{hyperref}

\theoremstyle{definition} 

\theoremstyle{plain} 
\newtheorem{theorem}{Theorem}

\theoremstyle{remark} 
\newtheorem{remark}{Remark}

\usepackage{accents}
\newcommand*{\dt}[1]{
  \accentset{\mbox{\large\bfseries .}}{#1}}

\newcommand{\be}{{\mathbf e}}
\newcommand{\bu}{{\mathbf u}}
\newcommand{\bv}{{\mathbf v}}

\newcommand{\by}{{\mathbf y}}
\newcommand{\bA}{{\mathbf A}}
\newcommand{\bB}{{\mathbf B}}
\newcommand{\bC}{{\mathbf C}}
\newcommand{\bD}{{\mathbf D}}
\newcommand{\bF}{{\mathbf F}}
\newcommand{\bI}{{\mathbf I}}

\newcommand{\bQ}{{\mathbf Q}}

\newcommand{\bU}{{\mathbf U}}
\newcommand{\bV}{{\mathbf V}}
\newcommand{\bX}{{\mathbf X}}
\newcommand{\bY}{{\mathbf Y}}

\newcommand{\bLambda}{{\boldsymbol \Lambda}}
\newcommand{\bPhi}{{\boldsymbol \Phi}}

\newcommand{\Ra}[1]{{{#1}}}
\newcommand{\Rb}[1]{{{#1}}}

\usepackage[textsize=tiny]{todonotes}

\begin{document}

\begin{frontmatter}



\title{Gradient Preserving Operator Inference: \\ Data-Driven Reduced-Order Models for Equations with Gradient Structure}



\author[1]{Yuwei Geng}
\ead{ygeng@email.sc.edu}
\author[1]{Jasdeep Singh}
\ead{jasdeep@email.sc.edu}
\author[1]{Lili Ju}
\ead{ju@math.sc.edu}
\author[2]{Boris Kramer}
\ead{bmkramer@ucsd.edu}
\author[1]{{Zhu Wang}\corref{cor1}}
\ead{wangzhu@math.sc.edu}
\cortext[cor1]{Corresponding author}
\address[1]{Department of Mathematics, University of South Carolina, Columbia, SC 29208, USA}
\address[2]{Department of Mechanical and Aerospace Engineering, University of California San Diego, La Jolla, CA 92093, USA}

\begin{abstract}
Hamiltonian Operator Inference has been introduced in [\textit{Sharma, H., Wang, Z., Kramer, B., Physica D: Nonlinear Phenomena, 431, p.133122, 2022}] to learn structure-preserving reduced-order models (ROMs) for Hamiltonian systems. This approach constructs a low-dimensional model using only  data and knowledge of the Hamiltonian function. Such ROMs can keep  the intrinsic structure of the system, allowing them to capture the physics described by the governing equations. In this work, we extend this approach to more general systems that are either conservative or dissipative in energy, and which possess a gradient structure. We derive the optimization problems for inferring structure-preserving ROMs that preserve the gradient structure.  We further derive an {\em a priori} error estimate for the reduced-order approximation. To test the algorithms, we consider semi-discretized partial differential equations with gradient structure, such as the parameterized wave and Korteweg-de-Vries equations, \Ra{and equations of  three-dimensional linear elasticity} in the conservative case and the one- and two-dimensional Allen-Cahn equations in the dissipative case. The numerical results illustrate the accuracy, structure-preservation properties, and predictive capabilities of the gradient-preserving Operator Inference ROMs. 

\end{abstract}



\begin{keyword}
dissipative and conservative systems\sep gradient-flow equation\sep Operator Inference\sep reduced-order modeling \sep data-driven modeling


\MSC[2020] 65P99\sep 65M15

\end{keyword}

\end{frontmatter}



\section{Introduction}

Consider evolutionary partial differential equations  (PDEs) of the following form 
\begin{equation}
\frac{\partial y(x, t; \mu)}{\partial t} = \mathcal{L} \ \frac{\delta \mathcal{H}[y; \mu]}{\delta y},
\label{eq:h-pde}
\end{equation}
where the state variable $y$ depends on the spatial variable $x\in \Omega \subset \mathbb{R}^d$ $(d=1, 2 \text{ or } 3)$, time variable $t\in I \subset [t_0, \infty)$ and parameter $\mu \in {P}\subset \mathbb{R}^p, p\geq 1$. Moreover, $\mathcal{L}$ represents a constant linear differential operator, and $\frac{\delta \mathcal{H}}{\delta y}$ is the variational derivative of the functional 
$$\mathcal{H}[y; \mu] := \int_{\Omega} H(y, y_x, y_{xx}, \ldots; \mu) \,dx$$
with $y_x$, $y_{xx}$, \ldots, denoting the first- and second-order (and higher order) partial derivatives of $y$ with respect to $x$. 
Depending on the properties of $\mathcal{L}$, the PDE can possess different features. 
If $\mathcal{L}$ is skew adjoint, $\mathcal{H}$, referred to as the Hamiltonian, is a constant and the PDE is conservative. 
If $\mathcal{L}$ is negative semi-definite (resp. definite), $\mathcal{H}$, referred to as a Lyapunov function or energy function, is nonincreasing (resp. monotonically decreasing) and the PDE is dissipative.

Either class of PDEs has broad applications ranging from mechanics \cite{marsden2013introduction}, plasma physics \cite{arakawa_1966_incompressible,pagliantini_2023_time_integration} to ocean modeling \cite{salmon1982shape,stewart2010multilayer,stewart2016energy,chen2021conservative}. 
Numerical methods have been developed to approximate and predict solutions to such PDEs.  It is well known that numerical schemes that do not preserve the continuous models' properties are prone to suffer from error accumulation and generate unphysical solutions over a long-term integration.  Therefore, structure preserving algorithms have been designed to keep the right features at a discrete level, see the review papers \cite{christiansen2011topics,sharma2020review} and the reference therein. In particular, structure-preserving spatial discretization schemes include finite difference methods such as \cite{fukao2017structure,cai2019structure} and finite element methods such as \cite{gross2005conservation,xu2008discontinuous,christiansen2011topics,sun2020structure,sanchez2021symplectic}. Structure preserving time discretization schemes include geometric integrators \cite{mclachlan1998six,mclachlan2006geometric,hairer2006geometric}, 
discrete gradient methods \cite{quispel1996discrete,mclachlan1999geometric}, 
average vector field (AVF) methods \cite{quispel2008new,celledoni2012preserving}, 
and discrete variational derivative methods \cite{furihata2010discrete}.

The cost (e.g., measured in CPU hours) to simulate a high-dimensional discretization of \eqref{eq:h-pde} can become a limiting factor when one wants to simulate either for long time horizons, and/or obtain many parametric solutions for different $\mu_i$ in the context of exploration and uncertainty quantification. 
To accelerate the simulation of high-dimensional systems, reduced-order models (ROMs) can provide orders-of-magnitude speedups while providing accurate approximate solutions, see the surveys/books \cite{antoulas2005approximation,quarteroni2015reduced,benner2015survey,antoulas2020interpolatory,brunton2022data}. 
However, ROMs constructed by projecting the equation onto a linear subspace spanned by a reduced basis, or onto a nonlinear manifold, could destroy the features of mechanical systems. In the context of Lagrangian structure, this has been shown in \cite{lall2003structure}, and motivated the development of projection-based model reduction techniques that preserve the Lagrangian structure \cite{carlberg2015preserving}. 
For systems with Hamiltonian and port-Hamiltonian systems, structure-preserving projection-based model reduction has been developed, including the interpolatory projection method proposed in \cite{beattie2009interpolatory} that preserves the symmetry of linear dynamical systems and is extended to port-Hamiltonian systems via tangential rational interpolation in \cite{gugercin2012structure}; the \Ra{proper orthogonal decomposition (POD)}-based Galerkin projection method in \cite{gong2017structure,karasozen2018energy,miyatake2019structure} for Hamiltonian systems and the POD/$\mathcal{H}_2$-based Petrov-Galerkin project method in \cite{beattie2011structure,chaturantabut2016structure} for port-Hamiltonian; the proper symplectic decomposition in \cite{peng2016symplectic} and a Petrov-Galerkin projection-based, variationally consistent approach in \cite{gruber2024variationally} 
 for canonical Hamiltonian systems; and the reduced basis method for Hamiltonian systems in \cite{afkham2017structure,hesthaven2021structure,pagliantini2021dynamical,hesthaven2022reduced}. Moreover, there has been increased interest lately in projection-based model reduction on nonlinear manifolds for Hamiltonian systems, see \cite{sharma2023symplectic,buchfink2023symplectic,yildiz2023data}.

ROMs can also be constructed non-intrusively, which is appealing in situations where only snapshot data of state variables is available, paired with potentially some prior knowledge of the PDE model. 
\Ra{In particular, the operator inference (OpInf) method provides such a nonintrusive framework for model order reduction, which was introduced in \cite{peherstorfer2016data} and expanded in many directions, c.f. the survey paper \cite{KPW_OPINF_survey2024}. The recent work~\cite{mcquarrie2023nonintrusive} extended OpInf to affine-parametric  systems of PDEs, yet the adaption of the framework to problems with unknown parametric structures is nontrivial. Within the OpInf framework,} learning canonical Hamiltonian ROMs on linear subspaces is considered in \cite{sharma2022hamiltonian}, where the linear part of the gradient of the Hamiltonian is inferred through constrained least-squares solutions. This has been extended to noncanonical Hamiltonian systems in \cite{Gruber2023canonical}, where the entire Hamiltonian function is assumed to be known, and the reduced operator associated to the linear differential operator $\mathcal{L}$ is inferred. Port-Hamiltonian \cite{morandin2023port} and Lagrangian ROMs can also be inferred directly from time-domain data, see \cite{SK_LagrangianOPINF,filanova2023operator,ShNaToKr-Lagrangian-OPINF-SPML-2024}.

Motivated by these developments, in this paper, we propose {\em gradient-preserving Operator Inference (GP-OPINF)}, a more general framework to nonintrusively construct ROMs that preserve the gradient structure of the PDE model \eqref{eq:h-pde}. In this nonintrusive setting, we assume that we are given snapshot data and some prior knowledge of the continuous model and the functional $\mathcal{H}$, yet that we do not have access to the spatially-discretized forms in high-fidelity solvers. We note that the proposed approach is not limited to Hamiltonian systems.
\Ra{The novelties of this work include the following:}   
\begin{enumerate}
\item We propose suitable optimization problems for the gradient-preserving Operator Inference method. In particular, for conservative equations, the solution to the optimization problem is post-processed for improving the gradient structure of ROMs without introducing any consistency errors. 
\item We derive an {\em a priori} error estimate for the inferred GP-OPINF ROM, which we demonstrate by numerical experiments.
\item Our framework applies to dissipative and conservative systems with a gradient structure. 
\end{enumerate}

The rest of this paper is organized as follows.
In Section~\ref{sec:rom}, we introduce the full-order model form that we consider and, for comparison, an intrusive structure-preserving ROM for equation \eqref{eq:h-pde}.
In Section~\ref{GP-OpInf}, we introduce the gradient-preserving Operator Inference and present the associated optimization problems.
In Section~\ref{sec:err} we derive an {\em a priori} error estimate for the learned GP-OPINF ROM.
We then demonstrate the effectiveness of the proposed ROM through several numerical examples in Section~\ref{sec:num}.
A few concluding remarks are drawn in the last section.

\section{Full-Order Model and Structure-Preserving Reduced-Order Model}\label{sec:rom}

We present the full-order model setting in Section~\ref{ss:FOM} and then discuss the traditional structure-preserving intrusive ROM approach in Section~\ref{ss:intrusiveROM}. This provides the background for the new method in the following Section~\ref{GP-OpInf}.

To ease the exposition, in situation when no confusion arises, we omit the explicit reference to the parameter and adopt the following notations: at any $\mu\in {P}$, for any vector-valued function $\mathbf{f}$ and scalar-valued function $g$, 
$$\mathbf{f}(t)\equiv \mathbf{f}(t;\mu),\,\,
g(\mathbf{f}) \equiv g(\mathbf{f}(t; \mu); \mu), \quad \forall t\in I.$$

\subsection{Full-Order Model} \label{ss:FOM}
For any parameter (vector) $\mu \in {P}\subset \mathbb{R}^p, p\geq 1$, we call the system of finite dimensional ODEs obtained from \eqref{eq:h-pde} after a structure-preserving spatial discretization the full-order model (FOM). Here, we consider the FOM
\begin{equation}
\dt{\by} = \bD\, \nabla_{\by} H(\by),
\label{eq:ham_fom}
\end{equation}
where $\by(t)\equiv \by(t;\mu) \in \mathbb{R}^n$ is the $n$-dimensional state vector and the initial condition is $\by(t^0)= \by^0$. 
The right-hand-side of \eqref{eq:ham_fom} represents the velocity field (or drift),  in which $\bD\in \mathbb{R}^{n\times n}$ is a real matrix, $H(\by) \equiv H(\by; \mu): \mathbb{R}^n \rightarrow \mathbb{R}$ is a continuously differentiable function of the state and $\nabla_{\by} H(\by)$ is the gradient of $H(\by)$. The gradient structure of equation~\eqref{eq:ham_fom} plays a key role in characterizing the dynamics:

\begin{enumerate}[(i)]
\item When $\bD$ is anti-symmetric (also referred to as skew-symmetric or skew-adjoint), that is, $\bD= -\bD^\intercal$, the system is Hamiltonian. 
A special case is when 
$\bD=\begin{bsmallmatrix} 0&\bI\\-\bI&0 \end{bsmallmatrix}$
then equation~\eqref{eq:ham_fom} represents a canonical Hamiltonian system and the solution flow is symplectic.  Due to such a structure, the internal energy of the system, $H(\by)$, is conserved. Mathematically, for any $t_1, t_2\in I$ with $t_1<t_2$, we have  
\begin{equation*}
H(\by(t_2)) - H(\by(t_1)) 
= \int_{t_1}^{t_2} \frac{{\rm d}}{{\rm d}t} H(\by(t))\, {\rm d}t 
= \int_{t_1}^{t_2} \left(\nabla_{\by} H(\by)\right)^\intercal  \bD\, \nabla_{\by} H(\by) \, {\rm d}t 
= 0.
\label{eq:h}
\end{equation*}
\item When $\bD$ is negative semi-definite, the system represents a gradient flow and is dissipative.  To see this, we can easily check  
\begin{equation*}
H(\by(t_2)) - H(\by(t_1)) 
= \int_{t_1}^{t_2} \frac{{\rm d}}{{\rm d}t} H(\by(t))\, {\rm d}t 
= \int_{t_1}^{t_2} \left(\nabla_{\by} H(\by)\right)^\intercal  \bD\, \nabla_{\by} H(\by) \, {\rm d}t 
\leq 0. 
\label{eq:gf}
\end{equation*}
Note that if $\bD$ is negative definite, $H(\by)$ is strictly decreasing.
Numerical schemes to solve these equations recognize the special gradient structure in time discretization, e.g., geometric integrators and average vector field methods. 
\end{enumerate}

If the size of the dynamical system is large, it is expensive to simulate \eqref{eq:ham_fom}, especially for many parameters or for long time intervals. We next review intrusive projection-based ROMs that can provide accurate approximations at much lower computational cost. To preserve physics of the system, they should possess the same gradient structure. 


\subsection{Intrusive Structure-Preserving ROMs} \label{ss:intrusiveROM}
We review an intrusive Galerkin-projection-based approach to build the structure-preserving ROM, termed {\em SP-G} ROM, which was introduced in  \cite{gong2017structure}.
The first step in the ROM construction is to extract the low-dimensional basis in which we derive the ROM. Here, we use the widely used snapshot-based POD method \cite{berkooz1993proper}. 
For a non-parametric problem, we consider $\mu$ as fixed, and the snapshots (state vectors at select time instances) are stored in the matrix
\begin{equation}
\bY \equiv \bY(\mu) \coloneqq \left[ \by(t^0; \mu), \by(t^1; \mu), \ldots, \by(t^s; \mu) \right].
\label{eq:data_y}
\end{equation}
For simplicity, we assume the snapshot set contains the FOM solutions at all the uniformly distributed discrete time discretization points, that is, $t^{i+1} = t^i+\Delta t$ for $i=0, \ldots, s-1$. However, in practice one may use only a portion of them as the snapshots. 
For a parametric problem, one is interested in a set of parameters $\{\mu^0, \mu^1, \cdots, \mu^p\}$ from ${P}$, for which full-order solutions are generated. The associated snapshot data matrix is defined as  
\begin{equation}
\bY \coloneqq \left[\bY(\mu^0), \bY(\mu^1), \ldots, \bY(\mu^p)\right]. 
\label{eq:data_y2}
\end{equation}
The POD basis matrix $\bPhi= [\mathbf{\phi}_1, \mathbf{\phi}_2, \ldots, \mathbf{\phi}_r]$ contains the leading left singular vectors of $\bY$ as its columns, which can be computed by the method of snapshots \cite{sirovich1987turbulence} or singular value decomposition (SVD) and related randomized algorithms.

The reduced-order approximation of the state variable is defined by a linear ansatz:  
\begin{equation}
\Rb{\widehat{\by}(t; \mu)} \coloneqq \bPhi \by_r(t;\mu),
\label{eq:lin}
\end{equation}
where $\by_r(t; \mu) \in \mathbb{R}^r$ are the unknowns, which are determined by integrating the reduced-order model in time. 
\Ra{A centering trajectory (e.g., the average of snapshots) $\overline{\by}$ is usually subtracted from each snapshot before calculating the POD modes, which helps in enforcing fixed inhomogeneous Dirichlet boundary conditions. In that case, $\widehat{\by}(t; \mu) = \overline{\by}+\bPhi \by_r(t;\mu)$. 
We assume $\overline{\by}=\mathbf{0}$ in this work for simplicity. 
This linear ansatz has been adapted in \cite{gunzburger2007reduced} for reduced order modeling of problems with parametric, time-varying, inhomogeneous Dirichlet-type boundary conditions.
}

The {\em SP-G} ROM applies the ansatz \eqref{eq:lin}. To distinguish it from the non-intrusive reduced-order ROMs developed later in this paper, we let $\widetilde{\by}_r(t; \mu)$ denote the {\em SP-G} ROM's reduced-order state. The ROM is then constructed by first projecting the FOM onto the low-dimensional basis given by $\bPhi$, and subsequently adjusting it to ensure the correct gradient structure appears in the resulting ROM. This yields the {\em SP-G} ROM  
\begin{flalign}
\dt{\widetilde{\by}}_r = \widetilde{\bD}_r\, \nabla_{\widetilde{\by}_r} {H}_r(\widetilde{\by}_r) 
\label{eq:ham_rom}
\end{flalign}
for $\widetilde{\by}_r(t) \equiv \widetilde{\by}_r(t; \mu) \in \mathbb{R}^r$  with $\widetilde{\by}_r(t^0) = \bPhi^\intercal \by^0$,  
where $\widetilde{\bD}_r = \bPhi^\intercal \bD \bPhi$, $H_r(\widetilde{\by}_r) \coloneqq H(\bPhi \widetilde{\by}_r)$, and $\nabla_{\widetilde{\by}_r} H_r(\widetilde{\by}_r) = \bPhi^\intercal \nabla_{\by} H(\bPhi \widetilde{\by}_r)$. 
With this construction, the  right-hand side of the ROM \eqref{eq:ham_rom} has the same gradient structure as the FOM \eqref{eq:ham_fom}. Thus, the ROM has the same properties as its full-order counterpart. 
If the system is Hamiltonian, $H_r$ is also guaranteed to be constant. 

We remark that an error exists between $H_r(\widetilde{\by}_r)$ and $H(\by)$, which is due to the projected initial condition $\widetilde{\by}_r(t^0)$ in the ROM. However, this error shrinks as $r$ increases, as a larger ROM basis produces a more accurate ROM initial condition. 
Specifically, when the initial condition is independent of any parameters, it can be used as the centering trajectory, eliminating the error caused by projecting the initial condition.  
%
For details, the reader is referred to \cite{gong2017structure}. 

The above construction of {\em SP-G} ROM requires access to the full-order coefficient matrix $\bD$, knowledge of $H$ or its gradient $\nabla_{\by} H$, and snapshot data. 
In the next section, we introduce an Operator Inference approach to build a data-driven ROM for \eqref{eq:ham_fom} with the appropriate gradient structure, but without accessing the full-order coefficient matrix.

\section{Gradient-Preserving Operator Inference} \label{GP-OpInf}
A structure-preserving Operator Inference for canonical Hamiltonian systems, named {\em H-OpInf}, has been developed in \cite{sharma2022hamiltonian}, based on the unconstrained operator inference from~\cite{peherstorfer2016data,KPW_OPINF_survey2024}. The symplecticity-preserving method learns the reduced-order operator associated with the quadratic part of the Hamiltonian function and assumes knowledge of the nonquadratic part of the Hamiltonian function. 
\Ra{The motivational example in \cite[Sec. 3.1]{sharma2022hamiltonian} shows that the OpInf ROM without preserving the symplectic structure produces an increasing Hamiltonian energy for a canonical linear wave equation, resulting in inaccurate and unphysical reduced-order simulations. This occurs already within the training data regime.}
Non-canonical Hamiltonian systems are considered in \cite{Gruber2023canonical}, where the authors construct a ROM that preserves the (Hamiltonian) energy under the assumption of prior knowledge of the entire Hamiltonian function. 
In this work, we propose a more general gradient-preserving Operator Inference (GP-OpInf) approach for both conservative and dissipative systems governed by the general FOM~\eqref{eq:ham_fom}. 
In Section~\ref{ss:GP-OPINF-overview} we introduce the new framework. Section~\ref{ss:GPOPINF_conservative} presents solution approaches for the conservative case and Section~\ref{ss:GPOPINF_dissipative} for the dissipative case.

\subsection{Proposed GP-OPINF Framework} \label{ss:GP-OPINF-overview}
 We assume that snapshot data $\bY$ given in equations \eqref{eq:data_y} or \eqref{eq:data_y2} as well as the symbolic expression of $H(\by)$ or the associated gradient is given. 
The goal is to develop a non-intrusive gradient-preserving ROM. We do this by employing the linear ansatz $\eqref{eq:lin}$, characterized by the expression $\Rb{\widehat{\by}(t; \mu)} = \bPhi \by_r(t;\mu)$, to find $\by_r(t)\equiv \by_r(t; \mu) \in \mathbb{R}^r$ satisfying
\begin{flalign}
\dt{\by}_r = \bD_r\, \nabla_{\by_r} {H}_r(\by_r) 
\label{eq:ham_rom_opinf}
\end{flalign}
with $\by_r(t^0) = \bPhi^\intercal \by^0$. Here $\bD_r$ is unknown and will be inferred from data. Specifically, $\bD_r$  is anti-symmetric if the system is conservative; and is negative semi-definite if the system is dissipative.  
We refer to equation~\eqref{eq:ham_rom_opinf} as the {\em GP-OpInf} ROM since it has the same gradient structure as the FOM~\eqref{eq:ham_fom}.

To learn $\bD_r$, the gradient data \Ra{are} computed based on knowledge of the gradient of the Hamiltonian function and available snapshot data. In the non-parametric case,   
\begin{flalign}
\bF \equiv \bF(\mu) &\coloneqq \left[ \nabla_{\by} H(\by(t^0; \mu); \mu), \nabla_{\by} H(\by(t^1; \mu); \mu), \ldots, \nabla_{\by} H(\by(t^s; \mu); \mu) \right].
\label{eq:F}
\end{flalign}
The first-order time derivative data can be generated by applying a finite difference operator, denoted as $\mathcal{D}_t[\cdot]$, on the state vectors: 
\begin{flalign*}
\dt{\bY} \equiv \dt{\bY}(\mu) &\coloneqq \left[ \mathcal{D}_t[{\by}(t^0; \mu)], \mathcal{D}_t[{\by}(t^1; \mu)], \ldots, \mathcal{D}_t[{\by}(t^s; \mu)] \right].
\end{flalign*}
For instance, 
$\mathcal{D}_t[\cdot]$ can be chosen to be the following second-order finite difference operator, satisfying  
\begin{flalign}
\mathcal{D}_t[{\by}(t^j; \mu)] = 
\left\{
\begin{array}{ll}
(-3\by(t^0; \mu)+4\by(t^1; \mu)-\by(t^2; \mu))/(2\Delta t) & j=0, \\
(\by(t^{j+1}; \mu)-\by(t^{j-1}; \mu))/(2\Delta t) & j=1, \ldots, s-1,\\
(\by(t^{s-2}; \mu)-4\by(t^{s-1}; \mu)+3\by(t^s; \mu))/(2\Delta t) & j=s. 
\end{array}
\right.
\label{eq:dt}
\end{flalign}
In the parametric case, 
$\bF \coloneqq \left[\bF(\mu^0), \bF(\mu^1), \ldots, \bF(\mu^p) \right]$ 
and 
$\dt{\bY} = \left[\dt{\bY}(\mu^0), \dt{\bY}(\mu^1), \ldots, \dt{\bY}(\mu^p)\right]$. 
In either case, using the POD basis $\bPhi$, we project each of the above data matrices onto the POD basis and obtain the projected data
\begin{flalign*}
\bY_r &= \bPhi^\intercal \bY, \quad 
\bF_r = \bPhi^\intercal \bF, \quad \text{and} \quad
\dt{\bY}_r= \bPhi^\intercal \dt{\bY}.
\end{flalign*}
The {\em GP-OpInf} learns the reduced operator $\bD_r$ from the following optimization problems: 

\begin{enumerate}[(i)]
\item  When the system is conservative: 
\begin{flalign}
\min_{\bD_r\in \mathbb{R}^{r\times r}} \frac{1}{2} \|\dt{\bY}_r - \bD_r \bF_r\|_F^2 \quad \text{ s.t. }\quad \bD_r^\intercal = - \bD_r.
\label{eq:optc}
\end{flalign}
\item When the system is dissipative: 
\begin{flalign}
\min_{\bD_r\in \mathbb{R}^{r\times r}} \frac{1}{2} \|\dt{\bY}_r - \bD_r \bF_r\|_F^2 \quad \text{ s.t. }\quad \bD_r \preccurlyeq 0 ,
\label{eq:opt3}
\end{flalign}
where $\bD_r \preccurlyeq 0$ indicates that $\bD_r$ is semi-negative definite. 
\end{enumerate}
The next two sections present suitable solutions to these respective optimization problems.

\subsection{Solution of the GP-OpInf Optimization for the Conservative Case} \label{ss:GPOPINF_conservative}
The solution to the optimization problem \eqref{eq:optc} can be derived using the method of Lagrange multipliers. Consider the Lagrangian function  
\begin{flalign*}
L(\bD_r) = \frac{1}{2} \left<\dt{\bY}_r - \bD_r \bF_r, \dt{\bY}_r - \bD_r \bF_r\right>_F + \left<\bLambda, \bD_r+\bD_r^\intercal \right>_F,
\end{flalign*}
where the Frobenius inner product is defined as $\left<\bA, \bB\right>_F:=\sum_{i,j}A_{ij}B_{ij}$ for two real matrices $\bA$ and $\bB$. The stationary point of the Lagrangian then satisfies 
\begin{flalign}
(\dt{\bY}_r - \bD_r \bF_r) \bF_r^\intercal & = \bLambda+\bLambda^\intercal, \label{eq:case1-op1}\\
\bD_r^\intercal & = - \bD_r. \label{eq:case1-op2}
\end{flalign}
Using the symmetry of the right hand side in equation~\eqref{eq:case1-op1} together with equation \eqref{eq:case1-op2}, we get the Lyapunov equation 
\begin{flalign}
\bF_r \bF_r^\intercal \bD_r + \bD_r \bF_r \bF_r^\intercal
=  \dt{\bY}_r \bF_r^\intercal - \bF_r \dt{\bY}_r^\intercal.
\label{eq:case1-op3}
\end{flalign}
We next outline three solution approaches to obtain $\bD_r$, each leading to a different ROM. We then provide numerical comparisons in Section~\ref{sec:num}. \\

{\em GP-OpInf-V}: This straightforward (``vanilla") approach directly solves the Lyapunov equation~\eqref{eq:case1-op3}.  Given square matrices $\bA$, $\bQ \in \mathbb{R}^{n\times n}$, for a general continuous-time Lyapunov equation 
$$
\bA \bX + \bX \bA^\intercal = \bQ,
$$ 
it is known that (1) if $\bA$ and $-\bA^\intercal$ have no eigenvalues in common, then the solution $\bX$ is unique; 
(2) if $\bQ$ is (anti-)symmetric, then $\bX$ must be (anti-)symmetric, see, e.g., \cite{antsaklis1997linear}.
Since $\bF_r\bF_r^\intercal$ in \eqref{eq:case1-op3} is only positive semi-definite, there might be zero eigenvalues, violating (1), and the system would become singular. 
\Rb{In this case, standard software libraries for numerical linear algebra, such as LAPACK~\cite{kaagstrom1996lapack}, provide a unique solution $\bD_r$, yet to a slightly perturbed system. This solution may not be symmetric to machine precision anymore, and}
due to numerical errors in the solver, $\bF_r\bF_r^\intercal$ can have small negative eigenvalues. 
This could introduce additional errors in the numerical simulation of the GP-OpInf ROM. \Rb{To resolve this issue, we next consider a regularized optimization.}\\

{\em GP-OpInf-P}: 
\Rb{We introduce a regularization term to the optimization: for some $\epsilon >0$ 
\begin{flalign}
\min_{\bD_r\in \mathbb{R}^{r\times r}} \frac{1}{2} \|\dt{\bY}_r - \bD_r \bF_r\|_F^2 + \frac{\epsilon}{2}\|\bD_r\|_F^2 \quad \text{ s.t. }\quad \bD_r^\intercal = - \bD_r.
\label{eq:opt}
\end{flalign}
We use similar arguments as above and define the Lagrangian
\begin{flalign*}
L(\bD_r) = \frac{1}{2} \left<\dt{\bY}_r - \bD_r \bF_r, \dt{\bY}_r - \bD_r \bF_r\right>_F + \frac{\epsilon}{2} \left<\bD_r, \bD_r\right>_F + \left<\bLambda, \bD_r+\bD_r^\intercal \right>_F, 
\end{flalign*}
whose stationary point satisfies 
\begin{flalign*}
(\dt{\bY}_r - \bD_r \bF_r) \bF_r^\intercal - \epsilon \bD_r & = \bLambda+\bLambda^\intercal, \\ 
\bD_r^\intercal & = - \bD_r. 
\end{flalign*}
This leads to the ``perturbed" Lyapunov equation
\begin{flalign}
(\bF_r \bF_r^\intercal +\epsilon \bI) \bD_r + \bD_r (\bF_r \bF_r^\intercal+\epsilon \bI)
=  \dt{\bY}_r \bF_r^\intercal - \bF_r \dt{\bY}_r^\intercal,
\label{eq:case1-op4}
\end{flalign}
where $\bI$ is the $r\times r$ identity matrix.}

\Rb{In order to prevent the eigenvalues of the coefficient matrix $\bF_r \bF_r^\intercal +\epsilon \bI$ from being non-positive, we define
\begin{equation}
\epsilon = \left\{ \begin{array}{ll}
\max(c_0,\, \alpha |\min(\lambda(\bF_r\bF_r^\intercal))|) & \text{ if } \min(\lambda(\bF_r\bF_r^\intercal))\leq c_0,\\
0 & \text{ otherwise},
\end{array}
\right.
\label{eq:epsilon}
\end{equation}
with a small positive constant $c_0$ and a positive constant $\alpha$. This approach  improves the anti-symmetry of the inferred $\bD_r$. However, comparing with \eqref{eq:case1-op3}, the perturbed equation could introduce a consistency error. 
} 
\\

{\em GP-OpInf}:
We propose an alternative approach to eliminate this common consistency error.
Suppose that $\widehat{\bD}_r$ satisfies \eqref{eq:case1-op3}, that is, 
\begin{flalign}
\bF_r \bF_r^\intercal \widehat{\bD}_r + \widehat{\bD}_r \bF_r \bF_r^\intercal
=  \dt{\bY}_r \bF_r^\intercal - \bF_r \dt{\bY}_r^\intercal.
\label{eq:case1-op-r1}
\end{flalign}
\Rb{We assume that the equation can be uniquely solved. In the case where $\bF_r\bF_r^\intercal$ has zero eigenvalues, the approach used in LAPACK should be applied to provide a unique solution \cite{kaagstrom1996lapack}. }
As discussed above, $\widehat{\bD}_r$ may not be exactly anti-symmetric. 
After transposing both sides of \eqref{eq:case1-op-r1} and taking the negative, we obtain 
\begin{flalign}
- \bF_r \bF_r^\intercal \widehat{\bD}_r -\widehat{\bD}_r \bF_r \bF_r^\intercal 
=  \dt{\bY}_r \bF_r^\intercal - \bF_r \dt{\bY}_r^\intercal.
\label{eq:case1-op-r2}
\end{flalign}
By averaging \eqref{eq:case1-op-r1} and \eqref{eq:case1-op-r2}, we find that \eqref{eq:case1-op3} admits the anti-symmetric solution 
\begin{flalign}
\bD_r =  (\widehat{\bD}_r -\widehat{\bD}_r^\intercal)/2.
\label{eq:case1-op5}
\end{flalign}
Thus, the inferred $\bD_r$, defined in \eqref{eq:case1-op5}, is superior to the one obtained from \eqref{eq:case1-op4} because it possesses the right property without causing any consistency errors to the optimization problem.

\subsection{Solution of the GP-OpInf Optimization for the Dissipative Case} \label{ss:GPOPINF_dissipative}
The solution to the optimization problem \eqref{eq:opt3} for the dissipative case can be obtained through an interior-point method \cite{renegar2001mathematical,helmberg1996interior}. 
By introducing  
    $
        f(\bX) =  \lVert \dt\bY_r - \bX\bF_r \lVert_{F}^2
    $
    and the logarithmic barrier 
    \begin{equation*}
    g(\bX) =\left \{
    \begin{aligned}
    &\sum_{i=1}^r \log(-\lambda_i(\bX)) & \text{if  }  \lambda_i(\bX) < 0 \\
    &-\infty & \text{otherwise},
    \end{aligned}
    \right.
    \end{equation*}
    where $\bX$ is an $r\times r$ matrix and $\lambda_i(\bX)$ is the $i$th eigenvalue of the matrix $\bX$, we can solve, instead of \eqref{eq:opt3}, an unconstrained optimization with an augmented objective function: 
    \begin{equation}\label{eq:sdp}
        \bD_r = \arg\min_{\widetilde{\bD}_r\in \mathbf{R}^{r\times r}} \left \{  f_{\beta}(\widetilde{\bD}_r) \coloneqq f(\widetilde{\bD}_r) - \beta g(\widetilde{\bD}_r) \right \}, 
    \end{equation}
    where $\beta>0$ is the barrier parameter. The corresponding (disspative) ROM \eqref{eq:ham_rom_opinf} using $\bD_r$ determined from \eqref{eq:sdp} is referred to as the {\em GP-OpInf} ROM.

 %

\Ra{
\begin{remark}
The online simulation time of the nonintrusive {\em GP-OpInf} ROM \eqref{eq:ham_rom_opinf} is similar to that of the intrusive {\em SP-G} ROM for the same $r$, as both ROMs have the same structure yet may have different model parameters. 
If the Hamiltonian is nonlinear, one can use hyper-reduction techniques to accelerate the evaluation of the nonlinear term. For instance, the DEIM approach presented in \cite{pagliantini2023gradient} can be used to reduce the cost for evaluating the nonlinear Hamiltonian while preserving the structure. 

Compared to {\em SP-G} ROM, the {\em GP-OpInf} ROM requires additional offline computation to infer $\mathbf{D}_r$. In the conservative case, we need to solve a $r\times r$ Lyapunov equation. Since the  matrix $\bF_r\bF_r^\intercal$ is usually dense, the Bartels-Stewart algorithm can be applied \cite{bartels1972algorithm} with an associated computational complexity of  $\mathcal{O}(r^3)$ flops. In the dissipative case, one can solve \eqref{eq:sdp} by a gradient descent method as outlined in \ref{app1}. The associated computational complexity is $\mathcal{O}(r^2sp)$ per iteration, where $sp$ represents the number of snapshots. 
\end{remark}
}

\section{Error Estimation}\label{sec:err}
Previous work in \cite{sharma2023symplectic,Gruber2023canonical} established an asymptotic (as $r\rightarrow n$) convergence of learned reduced operators to the associated projection-based reduced operators. However, this analysis can sometimes not explain the error in the solution of the ROM, as small errors in the operators may cause large solution errors, specifically for structured (e.g., gradient) systems. 
Here, we take a step further and analyze the approximate \textit{solution error} of the inferred ROMs in the \textit{preasymptotic} regime. Specifically, we estimate the {\em a priori} error of the {\em GP-OpInf} ROM approximation relative to the FOM \eqref{eq:ham_rom_opinf}.  The derivation mainly follows the error estimation developed in \cite{gong2017structure} for analyzing the projection-based {\em SP-G} ROM. However, in addition to considering the error from POD truncation, there is a need to estimate the error caused by Operator Inference.  

Given a mapping $f: \mathbb{R}^n\rightarrow \mathbb{R}^n$, the Lipschitz constant and the logarithmic Lipschitz constant of the mapping are defined as
$$
\mathcal{C}_{\mathrm{Lip}}[f]:= \sup_{\bu\neq \bv}\frac{\|f(\bu)-f(\bv)\|}{\|\bu-\bv\|},
\quad \text{and} \quad
\mathcal{C}_{\mathrm{log-Lip}}[f]:= \sup_{\bu\neq \bv}\frac{\left<\bu-\bv, f(\bu)-f(\bv)\right>}{\|\bu-\bv\|^2},
$$
where  $\left<\cdot, \cdot\right>: \mathbb{R}^n\times \mathbb{R}^n \rightarrow \mathbb{R}$ for any positive integer $n$ denotes the Euclidian inner product.  The logarithmic Lipschitz constant could be negative, which was used in \cite{chaturantabut2012state} to show that the error of reduced-order solution is uniformly bounded over the time interval of simulations.

\begin{theorem}
Let $\by(t)$ be the solution of the FOM \eqref{eq:ham_fom} on the time interval $[0, T]$ and ${\by_r}(t)$ be the solution of the ROM  \eqref{eq:ham_rom_opinf} on the same interval. Suppose $\nabla_{\by} H(\by)$ is Lipschitz continuous, then the ROM approximation error satisfies
\begin{equation}
\int_0^T \hspace{-0.1cm} \|\by-\bPhi{\by_r}\|^2\, { \rm d} t 
\leq  
 C(T) \ \bigg( 
\underbrace{\int_0^T \hspace{-0.1cm} \|\by-\bPhi\bPhi^\intercal \by\|^2 { \rm d}t}_{\text{ \em projection error}} 
+ 
\underbrace{\int_0^T \hspace{-0.1cm} \|\dt{\by}-\mathcal{D}_t[\by]\|^2\, { \rm d}t}_{\text{\em data error}} 
+ 
\underbrace{\int_0^T \hspace{-0.1cm} \|\bPhi^\intercal \mathcal{D}_t[\by]  - \bD_r \bPhi^\intercal \nabla_{\by} H(\by)\|^2\, { \rm d}t}_{\text{\em optimization error}} 
\bigg), \label{eq:error_bound}
\end{equation}
where $C(T)= \max\{1+C_2^2, 2\} T \alpha(T)$, 
$\alpha(T)= 2\int_0^T e^{2C_1 (T-\tau)}\,{ \rm d} \tau$, and the constants $C_1= \mathcal{C}_{\mathrm{log-Lip}}[\bPhi \bD_r \bPhi^\intercal\nabla_{\by}H]$ and 
$C_2= \|\bPhi \bD_r \bPhi^\intercal\|\, \mathcal{C}_{\mathrm{Lip}}[\nabla_{\by}H]$.
\label{thm:1}
\end{theorem}

There are three terms in the error bound: the first term $\int_0^T \|\by-\bPhi\bPhi^\intercal \by\|^2 { \rm d} t$ is the {\bf projection error} caused by projecting $\by(t; \mu)$ onto the subspace spanned by the POD basis $\bPhi$, that interprets the POD approximability and has been analyzed thoroughly in \cite{kunisch2001galerkin,singler2014new}; the second one $\int_0^T \|\dt{\by}-\mathcal{D}_t[\by]\|^2\, { \rm d} t$ is the {\bf data error} caused by the approximation of time derivative snapshots using the finite difference method; the third one $\int_0^T \|\bPhi^\intercal \mathcal{D}_t[\by]  - \bD_r \bPhi^\intercal \nabla_{\by} H(\by)\|^2\, { \rm d} t$ represents the {\bf optimization error} generated by learning $\bD_r$ from data.

\begin{proof}
Consider the ROM of a fixed dimension $r$, and define its approximate solution error by 
\begin{equation*}
\be\coloneqq \by-\bPhi {{\by_r}}.
\label{eq:e0}
\end{equation*}
It can be decomposed as $\be =\boldsymbol{\rho} +\boldsymbol{\theta}$ with $\boldsymbol{\rho}\coloneqq \by-\bPhi\bPhi^\intercal \by$ and $\boldsymbol{\theta} \coloneqq \bPhi\bPhi^\intercal \by-\bPhi \by_r$. 
Define 
$\boldsymbol{\zeta} \coloneqq \dt{\by} - \mathcal{D}_t[{\by}]$ and 
$\boldsymbol{\eta} \coloneqq \bPhi^\intercal \mathcal{D}_t[{\by}]  - \bD_r \bPhi^\intercal \nabla_{\by} H(\by)$. 
Consider the time derivative of $\boldsymbol{\theta}$
\begin{eqnarray}
\dt{\boldsymbol{\theta}} &=& \bPhi \bPhi^\intercal \dt{\by} - \bPhi \dt{\by}_r \nonumber \\
        &=& \bPhi \bPhi^\intercal\dt{\by}  -
          \bPhi \bD_r\, \nabla_{\by_r} {H}_r(\by_r) \nonumber \\
		  &=& \bPhi \bPhi^\intercal \dt{\by}  - \bPhi \bD_r\, \nabla_{{\by_r}} H(\bPhi {\by_r}) \nonumber \\
		  &=& \bPhi \bPhi^\intercal\dt{\by}  - \bPhi \bPhi^\intercal \mathcal{D}_t[\by] \nonumber \\
		  &&+ \bPhi \bPhi^\intercal \mathcal{D}_t[\by]  - \bPhi \bD_r \bPhi^\intercal \nabla_{\by} H(\by) \nonumber \\
		  &&+ \bPhi \bD_r \bPhi^\intercal \nabla_{\by} H(\by) - \bPhi \bD_r \bPhi^\intercal \nabla_{\by} H(\bPhi \bPhi^\intercal \by) \nonumber \\
		  &&+ \bPhi \bD_r \bPhi^\intercal \nabla_{\by} H(\bPhi \bPhi^\intercal \by) - \bPhi \bD_r \bPhi^\intercal \nabla_{\by} H(\bPhi {\by_r}) .
\label{eq:d_theta00}		
\end{eqnarray}
Note that
\begin{equation}
\frac{{ \rm d}}{{ \rm d} t}\|\boldsymbol{\theta}\| 
= \frac{1}{2\|\boldsymbol{\theta}\|} \frac{{ \rm d}}{{ \rm d}t}\|\boldsymbol{\theta}\|^2
= \frac{1}{\|\boldsymbol{\theta}\|} \left<\boldsymbol{\theta},\,\dt{\boldsymbol{\theta}}\right>.
\label{eq:d_theta}
\end{equation}
Taking the inner product of equation~\eqref{eq:d_theta00} with $\theta$, the terms on the right-hand side of \eqref{eq:d_theta00} are
\begin{equation}
\left<\boldsymbol{\theta}, \bPhi \bPhi^\intercal\dt{\by}  - \bPhi \bPhi^\intercal \mathcal{D}_t[{\by}] \right> 
= \left<\boldsymbol{\theta}, \bPhi \bPhi^\intercal (\dt{\by} - \mathcal{D}_t[{\by}])\right>
 \leq \|\boldsymbol{\theta}\|\, \|\boldsymbol{\zeta}\|,
 \label{eq:err_projH0}
\end{equation}
\begin{equation}
\left<\boldsymbol{\theta}, \bPhi \bPhi^\intercal \mathcal{D}_t[{\by}]  - \bPhi \bD_r \bPhi^\intercal \nabla_{\by} H(\by) \right> 
= \left<\boldsymbol{\theta}, \bPhi \boldsymbol{\eta}\right>
 \leq \|\boldsymbol{\theta}\|\, \|\boldsymbol{\eta}\|,
 \label{eq:err_projH}
\end{equation}
\begin{equation}
\left<\boldsymbol{\theta}, \bPhi \bD_r \bPhi^\intercal\, \left(\nabla_{\by} H(\by) - \nabla_{\by} H(\bPhi \bPhi^\intercal\by)\right) \right>
\leq \|\boldsymbol{\theta}\|\, \|\bPhi \bD_r \bPhi^\intercal\|\, \mathcal{C}_{\mathrm{Lip}}[\nabla_{\by}H]\, \|\boldsymbol{\rho}\|,
 \label{eq:err_proju}
\end{equation}
\begin{equation}
\left<\boldsymbol{\theta}, \bPhi \bD_r \bPhi^\intercal\, \left(\nabla_{\by} H(\bPhi \bPhi^\intercal\by) - \nabla_{\by} H(\bPhi {\by_r})\right) \right>
\leq \mathcal{C}_{\mathrm{log-Lip}}[ \bPhi \bD_r \bPhi^\intercal\nabla_{\by}H]\, \|\boldsymbol{\theta}\|^2.
 \label{eq:err_MF}
\end{equation}
Let $C_1= \mathcal{C}_{\mathrm{log-Lip}}[\bPhi \bD_r \bPhi^\intercal\nabla_{\by}H]$,
$C_2= \|\bPhi \bD_r \bPhi^\intercal\|\, \mathcal{C}_{\mathrm{Lip}}[\nabla_{\by}H]$, and combine \eqref{eq:d_theta} with \eqref{eq:d_theta00} and \eqref{eq:err_projH0}-\eqref{eq:err_MF}, we have
\begin{equation*}
\frac{{ \rm d}}{{ \rm d} t}\|\boldsymbol{\theta}\| \leq C_1 \|\boldsymbol{\theta}\| + C_2 \|\boldsymbol{\rho}\| + \|\boldsymbol{\zeta}\| + \|\boldsymbol{\eta}\|.
\label{eq:d_theta0}
\end{equation*}
%
%
\Rb{
Applying the classical differential version of Gronwall lemma (see, for example, Appendix~B in \cite{evans2022partial}) over the interval $[0, t]$, for $t\leq T$, and 
using the fact that $\boldsymbol{\theta}(0) = 0$, we get}
\begin{equation*}
\|\boldsymbol{\theta}(t)\| \leq \int_0^t e^{C_1 (t-\tau)} (C_2 \|\boldsymbol{\rho}\| + \|\boldsymbol{\zeta}\| + \|\boldsymbol{\eta}\|)\, { \rm d} \tau.
\label{eq:d_theta1}
\end{equation*}
\Rb{Applying the Cauchy-Schwarz inequality on the RHS and squaring both sides, we have
\begin{flalign*}
\|\boldsymbol{\theta}(t)\|^2 
&\leq \int_0^t e^{2C_1 (t-\tau)}\,  { \rm d} \tau \int_0^t [C_2 \|\boldsymbol{\rho}\| +  (\|\boldsymbol{\zeta}\| + \|\boldsymbol{\eta}\|)]^2\, { \rm d} \tau \\
& \leq \alpha(t) \left[ \int_0^T C_2^2 \|\boldsymbol{\rho}\|^2 { \rm d} t+ \int_0^T (\|\boldsymbol{\zeta}\| + \|\boldsymbol{\eta}\|)^2\, { \rm d} t\right] \\
& \leq \alpha(T)\left[C_2^2 \int_0^T \|\boldsymbol{\rho}\|^2 { \rm d} t+ \int_0^T (\|\boldsymbol{\zeta}\| + \|\boldsymbol{\eta}\|)^2\, { \rm d} t\right],
\label{eq:d_theta2}
\end{flalign*}
where $\alpha(t)= 2\int_0^t e^{2C_1 (t-\tau)}\,  { \rm d} \tau$.
}
Hence,
\begin{equation*}
\int_0^T\|\boldsymbol{\theta}(t)\|^2\, { \rm d} t \leq T \alpha(T)\left[C_2^2 \int_0^T \|\boldsymbol{\rho}\|^2 { \rm d} t+ \int_0^T (\|\boldsymbol{\zeta}\| + \|\boldsymbol{\eta}\|)^2\, { \rm d} t\right].
\label{eq:d_theta3}
\end{equation*}
This, together with the triangular inequality and orthogonality of $\rho$ and $\theta$, yields
$$
\int_0^T \|\be(t)\|^2\, { \rm d} t \leq \left(1+ C_2^2 T \alpha(T)\right) \int_0^T \|\boldsymbol{\rho}\|^2 { \rm d} t+ 2 T \alpha(T) \int_0^T (\|\boldsymbol{\zeta}\|^2 + \|\boldsymbol{\eta}\|^2)\, { \rm d} t,
$$
which proves the theorem.
\end{proof}

\Ra{
\begin{remark}\label{re:2}
In comparison to the asymptotic error estimation of the {\em SP-G} ROM in \cite{gong2017structure}, the {\em GP-OpInf} ROM is different in that it only involves the projection error of the state, whereas the error of projecting the gradient of the Hamiltonian onto the reduced basis space is included in the {\em SP-G}. Additionally, it incorporates the data error and optimization error. 
We will numerically investigate these errors in the next section. 
\end{remark}
}

\begin{remark}
The FOM \eqref{eq:ham_fom} and ROM \eqref{eq:ham_rom_opinf} are continuous in time. In numerical simulations, time-marching schemes with identical or different time step sizes are applied. Consequently, the time discretization errors appear in the upper bound of the ROM approximation error. This time discretization error, along with the data error, can be controlled by choosing a sufficiently small time step size. 
\end{remark} 

\Rb{
\begin{remark}
It is worth mentioning that, in addition to \eqref{eq:F}, one could evaluate the gradient data on the projected snapshots, e.g., 
\begin{flalign}
\bF \equiv \bF(\mu) &\coloneqq \left[ \nabla_{\by} H(\bPhi \bPhi^\intercal \by(t^0; \mu); \mu), \nabla_{\by} H(\bPhi \bPhi^\intercal \by(t^1; \mu); \mu), \ldots, \nabla_{\by} H(\bPhi \bPhi^\intercal \by(t^s; \mu); \mu) \right]
\label{eq:F2}
\end{flalign}
when inferring the reduced operator. 
Using it does not alter the proposed {\em GP-OpInf} framework, yet for the error estimation, equations \eqref{eq:err_projH} and \eqref{eq:err_proju} can be grouped together. This results in a similar error bound to \eqref{eq:error_bound}, with the only difference being that the optimization error term becomes $\int_0^T \|\bPhi^\intercal \mathcal{D}_t[\by]  - \bD_r \bPhi^\intercal \nabla_{\by} H(\bPhi \bPhi^\intercal\by)\|^2$. 
\end{remark} 
}

\section{Numerical Experiments\label{sec:num}}
In this section, we provide many numerical experiments to demonstrate the effectiveness of the proposed method. In Section~\ref{sec:num_errors} we first define the error measures used throughout to assess quality and accuracy of a learned ROM. Section~\ref{sec:we} contains conservative dynamical systems and Section~\ref{sec:num_dissipative} evaluates dissipative dynamical systems.

\subsection{Error Measures}
\label{sec:num_errors}
Given a set of snapshots, the performance of a ROM depends on its dimension and the numerical schemes used for simulations.  
To test the performance of the various {\em GP-OpInf} ROMs and to illustrate the error analysis, for a given time discretization, we compute the (squared) \textit{ROM approximation error} as a function of $r$:  
\Rb{
\begin{equation}
   \mathcal{E}(r) := \frac{T\Delta A}{N} \sum_{j=1}^N \|\by(t^j; \mu)-\bPhi\by_r(t^j; \mu)\|^2, 
   \label{eq:e_approx}
\end{equation}
}
where $\|\cdot\|$ represents the vector 2-norm, $\by(t^j; \mu)$ and $\by_r(t^j; \mu)$ are the full-order and reduced-order solutions at $t^j$, respectively, for $j=1,\ldots, N$, $\bPhi= [\mathbf{\phi}_1, \mathbf{\phi}_2, \ldots, \mathbf{\phi}_r]$ is the reduced basis matrix, \Rb{and $\Delta A$ represents the associated element size in a uniform grid}. Following the error bound ~\eqref{eq:error_bound}, we additionally compute the (squared) \textit{projection error} and (squared) \textit{optimization error}: 
\Rb{
\begin{align}
\mathcal{E}_{\text{proj}}(r) & := \frac{T\Delta A}{N} \sum_{j=1}^N \|\by(t^j; \mu)-\bPhi\bPhi^\intercal \by(t^j; \mu)\|^2, \label{eq:e_proj} \\
\mathcal{E}_{\text{opt}}(r) & := \frac{T\Delta A}{N} \sum_{j=1}^N \|\bPhi^\intercal \mathcal{D}_t[\by(t^j; \mu)]- \bD_r \bPhi^\intercal \nabla_{\by} H(\by(t^j; \mu); \mu)\|^2.
\label{eq:e_opt}
\end{align}
}
The norms are evaluated by the same quadrature rule such as the composite rectangle method. 
In the following numerical experiments, we frequently use two matrices in discrete schemes:
$$
\mathbf{L}_n= 
\left( \begin{array}{cccccc}
-2 & 1 & 0 & \cdots & 1 \\
1 & -2 & 1 & \cdots & 0 \\
   &  \ddots & \ddots & \ddots &  \\
 0   &    \cdots          &  1 & -2  & 1 \\
 1 &    \cdots     &  0  &  1 & -2
\end{array} \right)_{n\times n}
\,\,
\text{ and } 
\,\,
\mathbf{S}_n = 
\left( \begin{array}{ccccc}
 0    & 1 & 0 & \cdots & -1 \\
 -1   & 0 & 1 & \cdots & 0 \\
 \vdots & \ddots & \ddots & \ddots & \vdots \\
    0 & \ldots &  -1 & 0  & 1 \\
    1 & 0 &  \ldots  &  -1 & 0
\end{array} \right)_{n\times n};
$$ 
and we denote the spatial domain of an equation by $\Omega$ and the time domain by $I$. Specifically, in the FOM simulations the time interval $I=[0, T_{\text{FOM}}]$, while in the ROM simulations $I=[0, T_{\text{ROM}}]$.

\Ra{For our numerical tests, the Lyapunov equation is solved by using the $\mathsf{lyap}$ function when using $\mathsf{Matlab}$, or the $\mathsf{linalg.solve\_continuous\_lyapunov}$ function when using $\mathsf{python/scipy}$.}
\Ra{Moreover, we use the licensed MOSEK \cite{mosek} optimization software developed by CVXPY \cite{diamond2016cvxpy,agrawal2018rewriting} to solve the semi-definite programming problem \eqref{eq:sdp}. If this software is unavailable, one can obtain the optimal solution through gradient descent, as discussed in \ref{app1}.  
}


\subsection{Conservative PDEs\label{sec:we}}
We consider the wave equation, the Korteweg-de Vries (KdV) equation and \Ra{the linear elastic equations} as a test-bed for investigating the {\em GP-OpInf} ROM when approximating conservative systems. After a spatial discretization of the PDEs, the former yields a canonical Hamiltonian system while the latter results in a non-canonical Hamiltonian system. We note that both equations have also been considered in \cite{sharma2022hamiltonian,Gruber2023canonical} with similar configurations.

\subsubsection{Parameterized wave equation}
Consider the one-dimensional linear wave equation with a constant wave speed $c$:
\begin{equation*}
u_{tt}= c^2 u_{xx}, \quad (x, t)\in \Omega\times I, 
\end{equation*}
where $u=u(x, t; \mu)$ is the parameter-dependent solution and $\mu\in P$ a parameter appearing in the initial condition, defined below alongside the boundary and initial conditions.
The wave equation can be recast to the {\em canonical} Hamiltonian formulation (see, e.g. \cite{celledoni2012preserving}) with the Hamiltonian $\mathcal{H}(u, v) = \int_\Omega (\frac{1}{2}v^2+\frac{c^2}{2} u_x^2) \, {\rm d}x$, so that
\begin{equation*}
\left[
\begin{array}{c}
u_t\\
v_t
\end{array}
\right] =
\left[
\begin{array}{cc}
0 & 1 \\
-1 & 0
\end{array}
\right]
\left[
\begin{array}{c}
\frac{\delta \mathcal{H}}{\delta u}\\
\frac{\delta \mathcal{H}}{\delta v}
\end{array}
\right],
\label{eq:lin_wave}
\end{equation*}
 which has the symplectic (gradient) structure. 

Taking a uniform partition in the spatial domain with the mesh size $\Delta x$ and defining a consistent discrete Hamiltonian $\overline{H}(t) \coloneqq \Delta x H(\bu, \bv)$ with $H(\bu, \bv) = \sum_j [\frac{1}{2} \bv_j^2 + \frac{c^2}{2(\Delta x)^2} (\bu_{j+1}-\bu_j)^2 ]$ leads to the 
Hamiltonian system of ODEs   
\begin{equation}
\left[
\begin{array}{c}
\dt{\bu}\\
\dt{\bv}
\end{array}
\right] =
\left[
\begin{array}{cc}
0 & \bI_n \\
-\bI_n & 0
\end{array}
\right]
\left[
\begin{array}{c}
-\bA \bu\\
\bv
\end{array}
\right]. 
\end{equation}
Here, $n$ is the number of degrees of freedom, $\bI_n$ the $n\times n$ identity matrix, and $\bA$ is a discretization of the scaled, one-dimensional second order differential operator. 
The above system is of the form \eqref{eq:ham_fom} with $\by = [\bu^\intercal, \bv^\intercal]^\intercal$ and 
$\bD$ is the $2n\times 2n$ skew-symmetric block matrix from above. This implies the system is conservative and the Hamiltonian is a constant function. 

\begin{figure}[htb!]
\centering
\begin{minipage}[ht]{0.32\linewidth}
\includegraphics[width=1\textwidth]{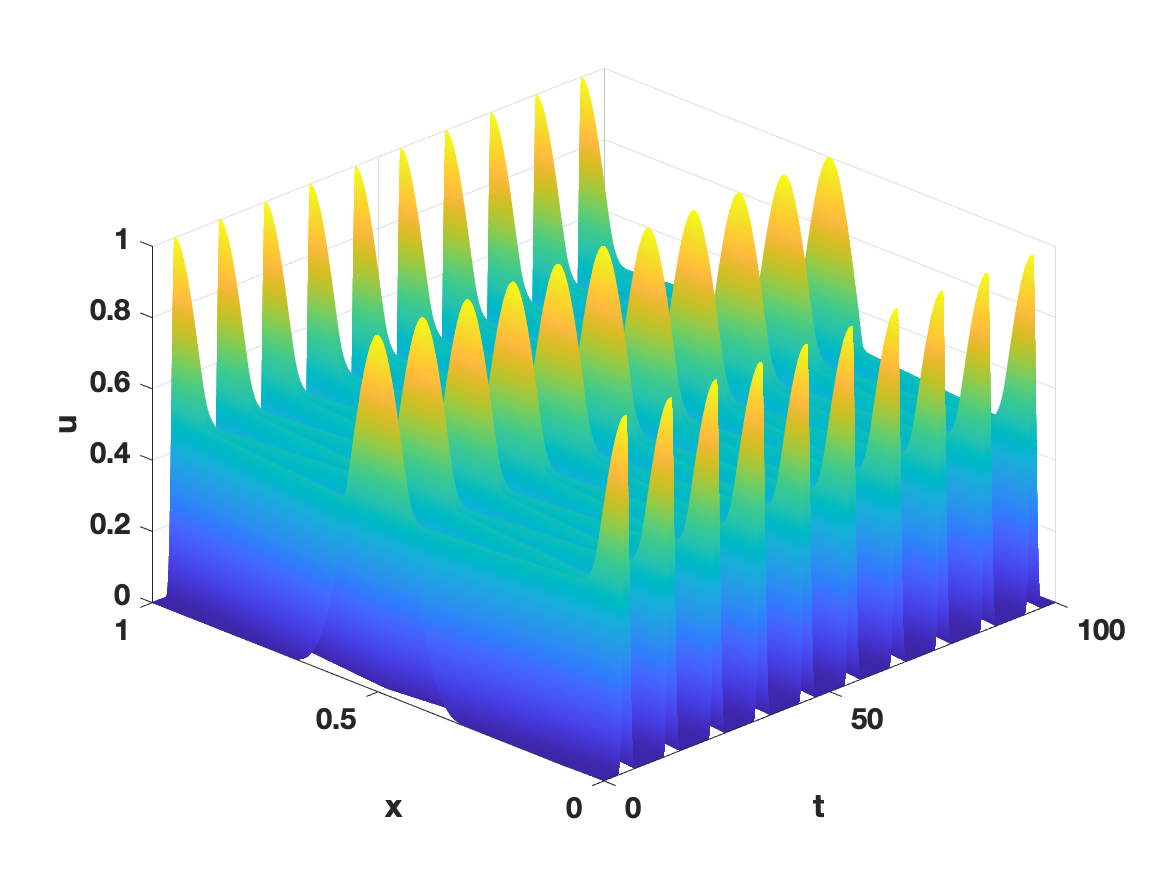}
\end{minipage}
\begin{minipage}[ht]{0.32\linewidth}
\includegraphics[width=1\textwidth]{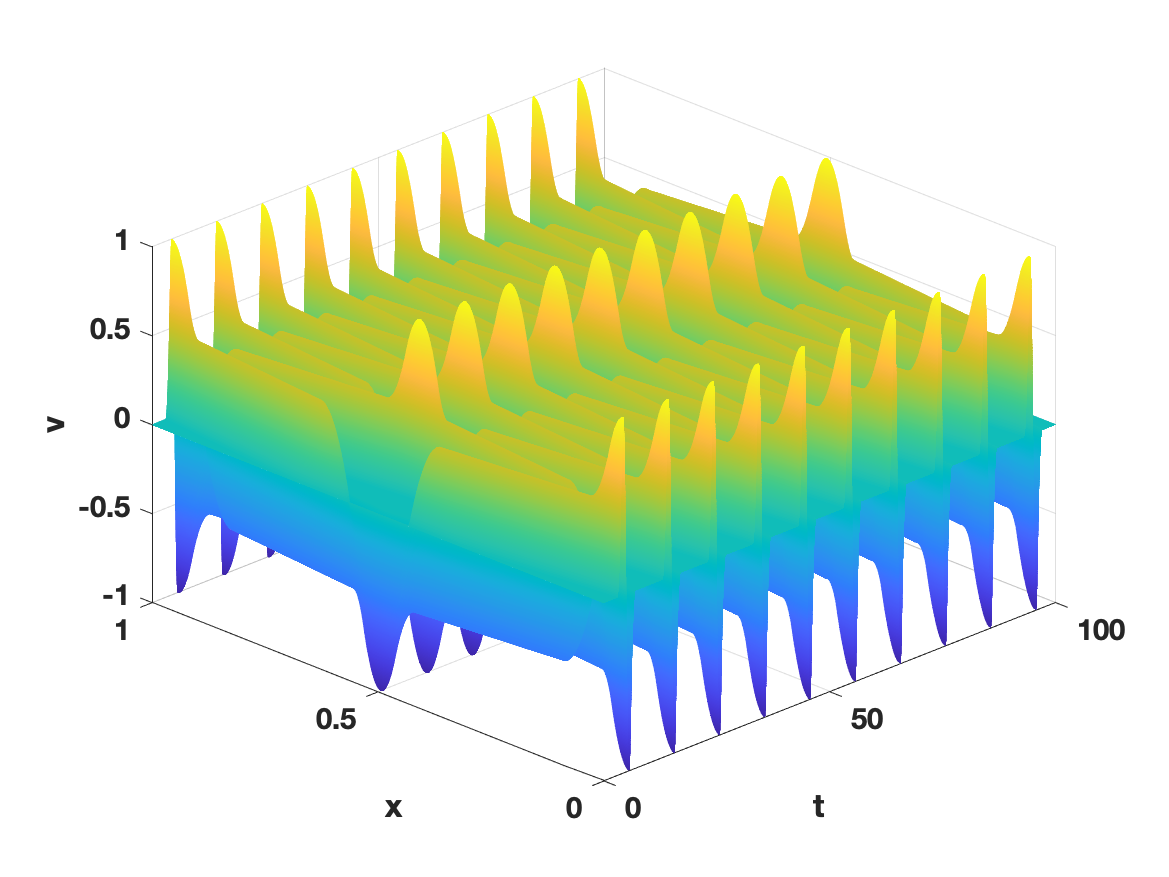}
\end{minipage}
\begin{minipage}[ht]{0.32\linewidth}
\includegraphics[width=1.2\textwidth]{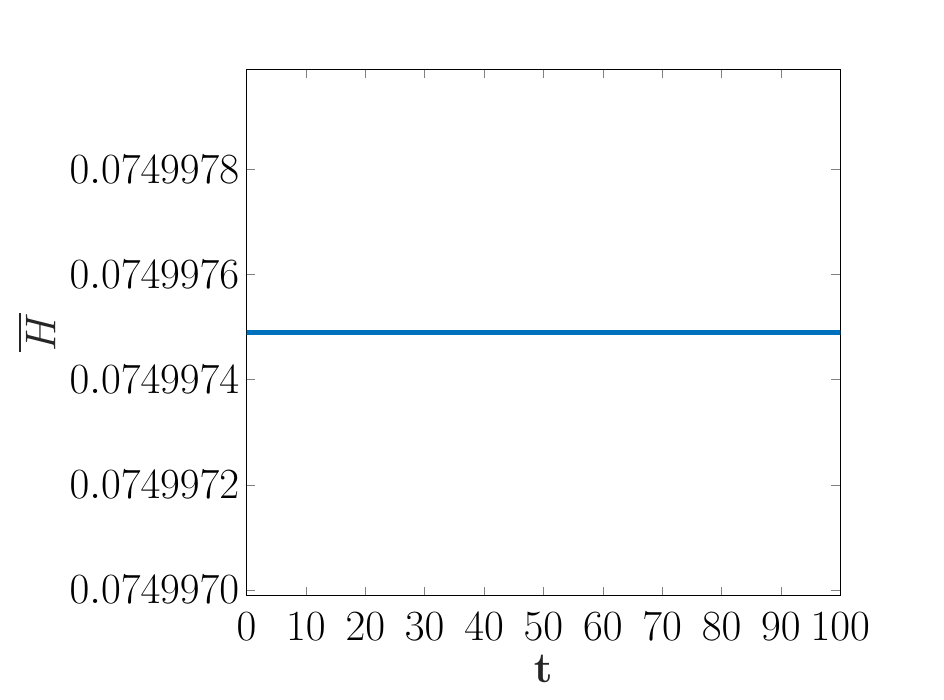}
\end{minipage}
\caption{
(Wave Equation) Full-order model simulation when $T_{\text{FOM}} = 100$: time evolution of $\bu$ (left), $\bv$ (middle), and energy $\overline{H}(t)$ (right).
}\label{fig:lin_wave_full}
\end{figure}

\paragraph{Computational setting} 
Consider the case in which $c= 0.1$ and $\Omega = [0, 1]$. 
The boundary condition is set to be periodic and the initial condition satisfies $u_0(x; \mu)=h(s(x; \mu))$ and $v_0(x; \mu) =0$, in which $h(s)$ is a cubic spline function defined by
\begin{equation*}
h(s) =
\left\{
\begin{array}{ll}
1-\frac{3}{2}s^2+\frac{3}{4}s^3 & \text{\quad if  \,\,} 0\leq s\leq 1, \\
\frac{1}{4}(2-s)^3 		       & \text{\quad if  \,\,} 1< s\leq 2, \\
0					       & \text{\quad if  \,\,} s> 2,
\end{array}
\right.
\quad \text{ and } \quad s(x; \mu)= \mu\left|x-\frac{1}{2}\right|
\end{equation*}
with $\mu$ a parameter that can take values from ${P}=[5, 15]$. 
In the FOM, the finite difference method is used for spatial discretization and 
$$
\bA= \frac{c^2}{\Delta x^2}\mathbf{L}_n.
$$ 
For time integration, the midpoint rule is applied with the time step $\Delta t$. 
The resulting linear system is solved by the built-in direct solver in MATLAB. 
For $T_{\text{FOM}} = 100$, $\Delta t= 1\times 10^{-3}$ and $\mu=10$, the time evolution of the full-order states and the energy $\overline{H}(t)$ are plotted in Figure~\ref{fig:lin_wave_full}.
Note that the Hamiltonian energy is preserved as $\overline{H}(t) \approx 7.5\times 10^{-2}$.

\paragraph{ROM construction} We assemble the snapshot matrices $\bU$ and $\bV$ from full-order simulations and generate the basis matrices $\bPhi_u\in \mathbb{R}^{n\times r_1}$ and $\bPhi_v\in \mathbb{R}^{n\times r_2}$ that contain, respectively, the $r_1$ and $r_2$ leading left singular vectors of $\bU$ and $\bV$. We define the POD basis $\bPhi =\begin{bsmallmatrix}
    \bPhi_u & \mathbf{0} \\ \mathbf{0} & \bPhi_v
\end{bsmallmatrix}$ 
and obtain the projected data $\bY_r = \bPhi^\intercal \bY = \bPhi^\intercal [\bU^\intercal, \bV^\intercal]^\intercal$ as well as the associated gradient of the Hamiltonian and the time derivative data matrices (using $\mathcal{D}_t$ defined in \eqref{eq:dt}). 
For simplicity, we choose $r_1=r_2=r$ in the following tests. However, different basis sizes can be chosen to better represent data. 
\Ra{We note that other types of reduced basis can be used, such as the cotangent lift \cite{peng2016symplectic}. We refer the interested readers to \cite{Gruber2023canonical} for a numerical investigation on different types of reduced bases.}
Next, we provide several comparisons of the ROMs generated by {\em GP-OpInf}, {\em GP-OpInf-V} and {\em GP-OpInf-P} and investigate their performance over $[0, T_{\text{ROM}}]$.

\paragraph{Test 1. Effects of regularization} 
We first consider the non-parametric case in which $\mu=10$ is fixed. The mesh size for the FOM is $\Delta x = 10^{-3}$ (correspondingly, $n = 1,000$) and the same time step $\Delta t = 1\times 10^{-3}$ is used in both the FOM and the ROMs. 
We use {\em GP-OpInf}, {\em GP-OpInf-V} and {\em GP-OpInf-P} to create ROMs of dimension $2r$ ($r$ bases for $u$ and $r$ bases for $v$) and let $r$ vary from $r=5$ to $r=200$.  

\begin{figure}[htb!]
\centering
\begin{minipage}[ht]{0.48\linewidth}
\includegraphics[width=1\textwidth]{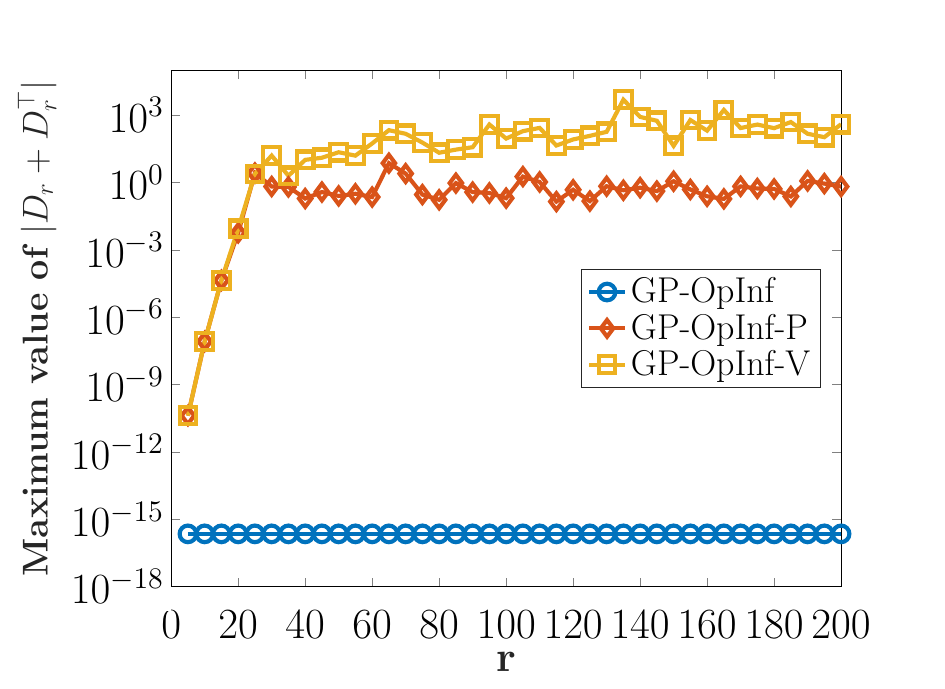}
\end{minipage}
\begin{minipage}[ht]{0.48\linewidth}
\includegraphics[width=1\textwidth]{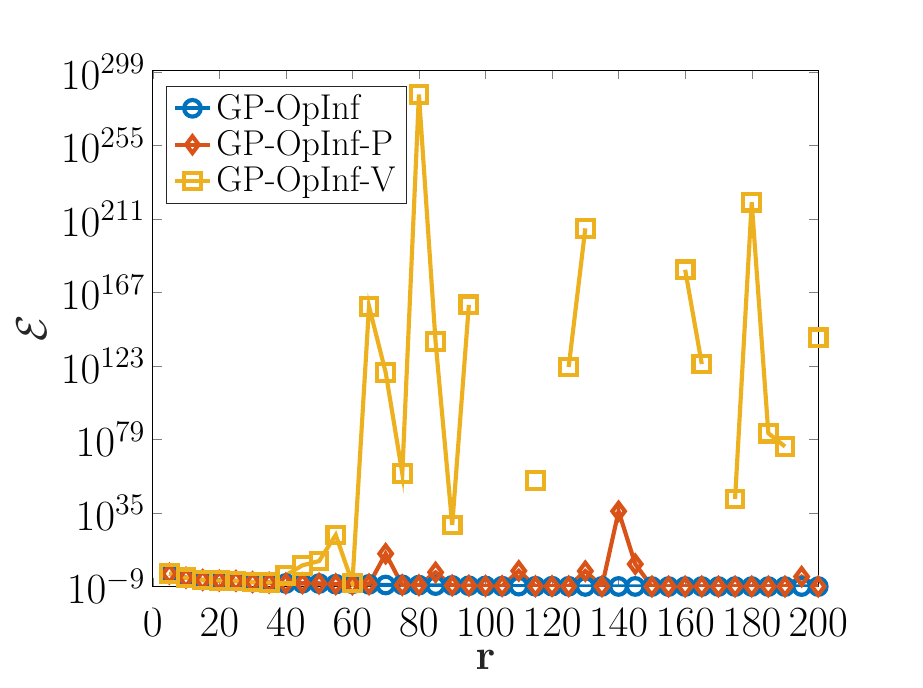}
\end{minipage}
\caption{(Wave Equation) ROMs of $2r$-dimensions when $T_{\text{FOM}} = 5$ and $T_{\text{ROM}} = 5$: 
(left) maximum magnitude of $|\bD_r+\bD_r^\intercal|$ for the inferred $\bD_r$ when $r$ varies;   
(right) ROM approximation error $\mathcal{E}$ from \eqref{eq:e_approx} when $r$ varies from 5 to 200 in increments of five. The inserted plots cover the portion of $r$ between 75 and 175. 
}\label{fig:lin_wave_reg}
\end{figure}

We first select $T_{\text{FOM}} = 5$, $T_{\text{ROM}} = 5$ and check whether the inferred $\bD_r$ possesses the desired skew-symmetric structure by plotting the maximum magnitude of $|\bD_r+\bD_r^\intercal|$ in Figure~\ref{fig:lin_wave_reg} (left). The inferred $\bD_r$ from {\em GP-OpInf-P} (with $\alpha=2$ and $c_0 = 10^{-13}$ in the definition of $\epsilon$ \eqref{eq:epsilon}) achieves better anti-symmetric structure than {\em GP-OpInf-V} when $r>30$, but both methods are not producing satisfactory results and are far off from being skew-symmetric.  In contrast, the proposed {\em GP-OpInf} produces exact anti-symmetric matrices up to machine precision.
The associated ROM approximation errors $\mathcal{E}$ are plotted in Figure~\ref{fig:lin_wave_reg} (right). For the {\em GP-OpInf-V} method, the error grows unbounded (no entry indicates an NaN) when $r> 40$. The {\em GP-OpInf-P} method yields better performance as it provides stable numerical solutions at almost all of the considered $r$-values, yet it becomes inaccurate at $r= 105$ and $r=155$. In contrast, the {\em GP-OpInf} method is able to obtain stable numerical approximations for all the $r$-values.   

\begin{figure}[htb!]
\centering
\begin{minipage}[ht]{0.48\linewidth}
\includegraphics[width=1\textwidth]{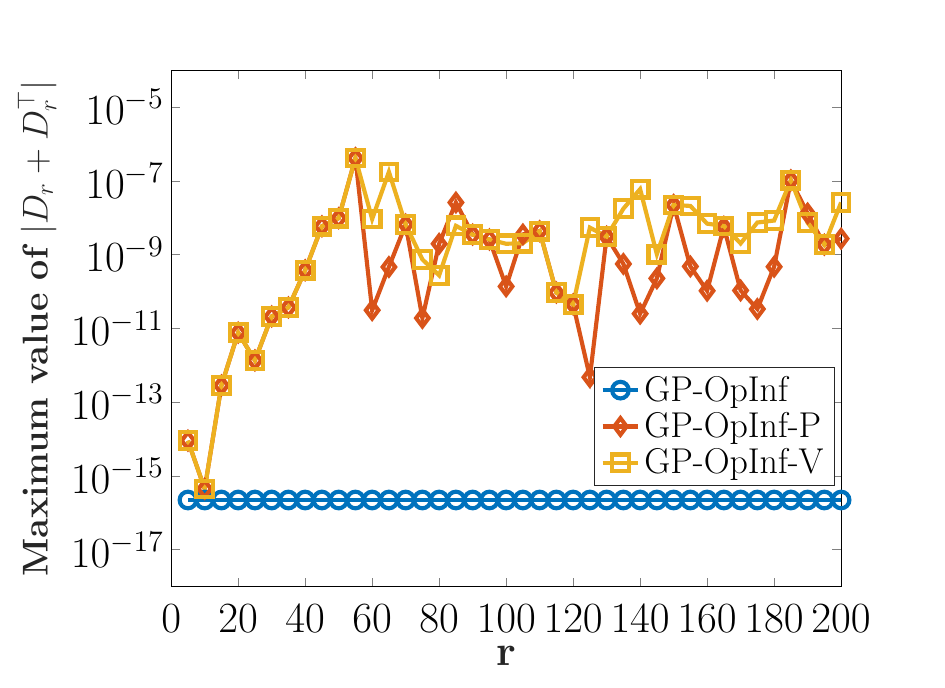}
\end{minipage}
\begin{minipage}[ht]{0.48\linewidth}
\includegraphics[width=1\textwidth]{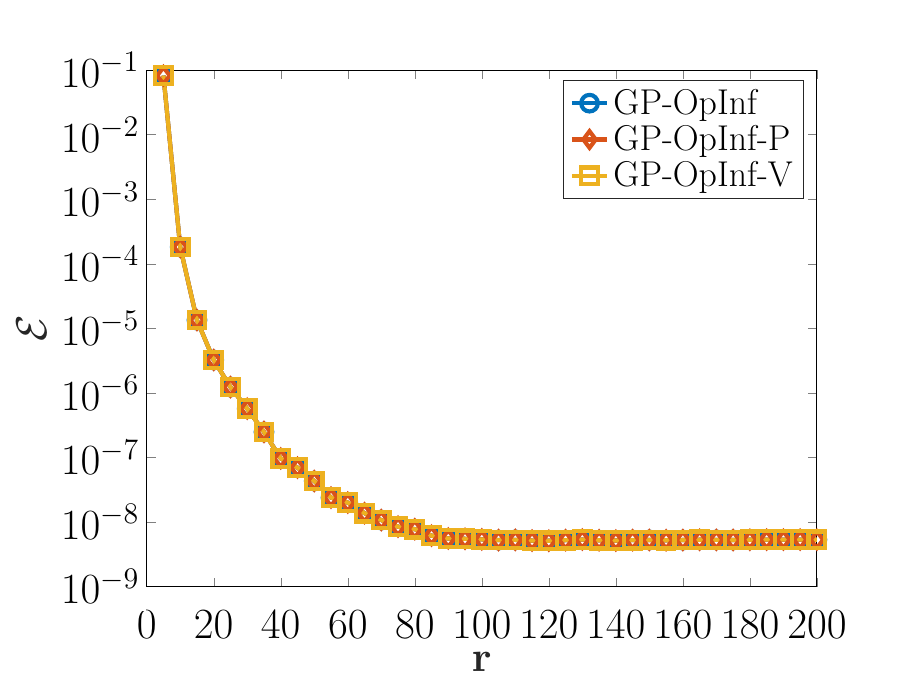}
\end{minipage}
\caption{(Wave Equation) ROMs of $2r$-dimensions when $T_{\text{FOM}} = 10$ and $T_{\text{ROM}} = 10$: 
(left) maximum magnitude of $|\bD_r+\bD_r^\intercal|$ for the inferred $\bD_r$ when $r$ varies;   
(right) numerical error $\mathcal{E}$ from \eqref{eq:e_approx}  when $r$ varies from 5 to 200 in increments of five. 
}\label{fig:lin_wave_reg2}
\end{figure}

We now consider the case of $T_{\text{FOM}} = 10$ and $T_{\text{ROM}} = 10$.  The associated maximum magnitude of $|\bD_r+\bD_r^\intercal|$ and the approximation errors of the ROMs based on three inferred $\bD_r$ are plotted in Figure~\ref{fig:lin_wave_reg2}. From the results, we find that all three methods learn $\bD_r$ with skew-symmetric structure, however, with noticeable differences. The maximum magnitude of  $|\bD_r+\bD_r^\intercal|$ is about $10^{-6}$ in {\em GP-OpInf-V} and {\em GP-OpInf-P}. Once again {\em GP-OpInf} is the most accurate yielding machine precision, $10^{-16}$, for the skew-symmetry test. From Figure~\ref{fig:lin_wave_reg2} (right), we observe that the three ROMs yield similar numerical accuracy in this case. 

\Rb{The different behavior of GP-OpInf-V in the second case is mainly a result of the inferred $\bD_r$ having a better structure. In Table~\ref{tab:mini_eig}, the minimum eigenvalues of $\bF_r\bF_r^\intercal$ are listed at several small $r$ values. We observe that small negative eigenvalues appear when $r=25$ and $r=30$ which result in the direct solution to the Lyapunov equation, $\bD_r$,  not being skew-symmetric. }
\begin{table}[htb!]
\centering
  \caption{Smallest eigenvalues of $\bF_r\bF_r^\intercal$ at several small $r$ values.} \label{tab:mini_eig}
		\begin{tabular}{lcccccc}
		\hline
		{Case} & $r= 5$ & $r= 10$ & $r= 15$ & $r= 20$ & $r= 25$ & $r= 30$  \\
		\hline
		$T_{\text{FOM}} = T_{\text{ROM}} = 5$ & 1.36e-07 & 8.16e-10 & 3.55e-12 & 9.57e-15 & -6.97e-13 & -6.41e-12 \\
		\hline
		$T_{\text{FOM}} = T_{\text{ROM}} = 10$ & 1.16e-04 & 6.74e-09 & 2.48e-10 & 1.98e-11 & 5.55e-12 & 1.26e-12  \\
		\hline
		\end{tabular}
\end{table}

In sum, we find that the {\em GP-OpInf} approach computes $\bD_r$ with skew-symmetric structure up to machine precision in both tests, which is markedly better than {\em GP-OpInf-V} and {\em GP-OpInf-P}. While more training data improved the {\em GP-OpInf-V} and {\em GP-OpInf-P} structural properties, noticeable (10 orders of magnitude) differences remained. The resulting {\em GP-OpInf} ROMs are more accurate, which is especially evident in the first test case where the other two approaches failed. In the following tests, we thus only focus on the {\em GP-OpInf} approach and investigate it in more detail.

\paragraph{Test 2. Illustration of error estimation}
Here, we investigate the error estimation for the {\em GP-OpInf} ROM for the same non-parametric case with $\mu= 10$, and compare the approximation error with the intrusive structure-preserving ROM ({\em SP-G} ROM) of the same dimension. 

\begin{figure}[htb!]
\centering
\begin{minipage}[ht]{0.48\linewidth}
\includegraphics[width=1\textwidth]{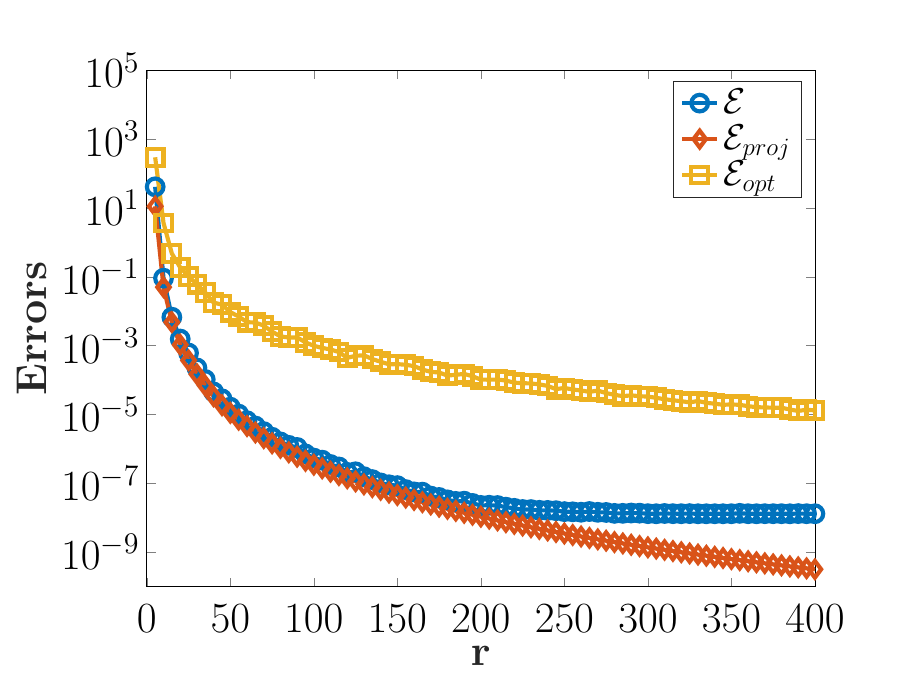}
\end{minipage}
\begin{minipage}[ht]{0.48\linewidth}
\includegraphics[width=1\textwidth]{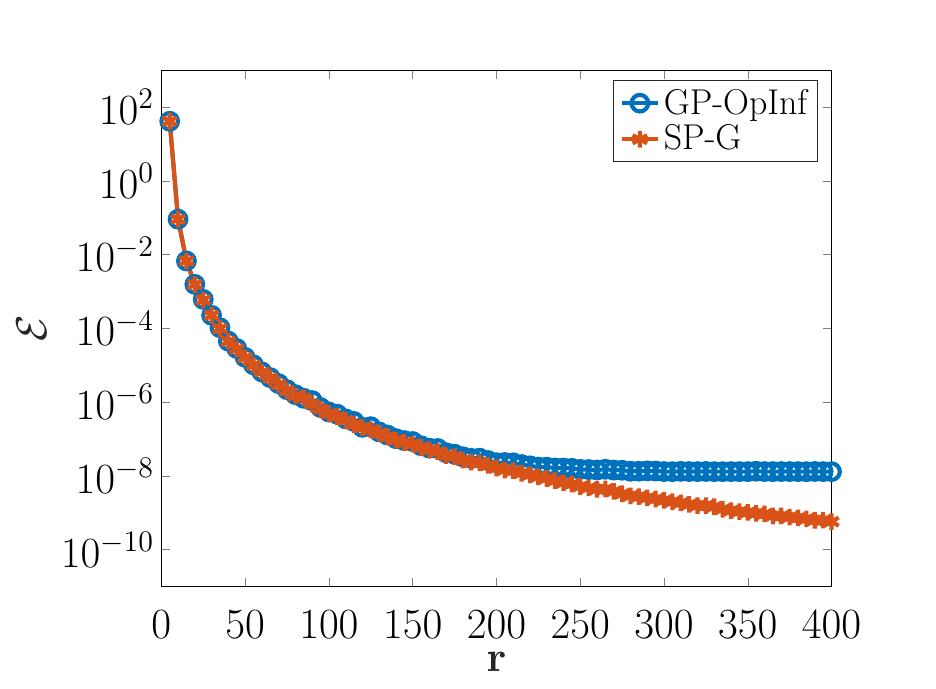}
\end{minipage}
\caption{(Wave Equation) Numerical errors of the $2r$-dimensional {\em GP-OpInf} ROM when $T_{\text{FOM}}=T_{\text{ROM}}=10$ and $\Delta x = \Delta t = 2\times 10^{-4}$: (left) the ROM approximation error \eqref{eq:e_approx} together with POD projection error \eqref{eq:e_proj} and optimization error \eqref{eq:e_opt}; (right) comparison of the ROM approximation error for the {\em GP-OpInf}  and {\em SP-G} ROM. 
}\label{fig:lin_wave_T10Tr10err}
\end{figure}

We choose $\Delta x = \Delta t=2\times 10^{-4}$ so that the discretization errors are negligible relative to the POD projection error and the optimization error.  We set $T_{\text{FOM}}=T_{\text{ROM}}=10$ and vary the dimension $r$ of the ROM. The three ROM errors \eqref{eq:e_approx}--\eqref{eq:e_opt} are plotted in Figure~\ref{fig:lin_wave_T10Tr10err} (left), which shows that the error decays monotonically as $r$ increases. When $r>200$, the ROM approximation error saturates, although the POD projection error is still decreasing, which indicates the error is dominated by the saturated optimization error and the fixed time discretization error.  
We compare the {\em GP-OpInf} ROM approximation error with the {\em SP-G} ROM in Figure~\ref{fig:lin_wave_T10Tr10err} (right). The {\em SP-G} approximation error decreases in a tendency similar to the POD projection error in Figure~\ref{fig:lin_wave_T10Tr10err} (left), and {\em GP-OpInf} achieves similar numerical accuracy to {\em SP-G} when $r$ is less than 200. However, as $r$ increases, the {\em GP-OpInf} ROM approximation accuracy stays around $10^{-8}$. Nevertheless, such large $r$ values lead to impractically expensive ROMs, and are of little practical use. 

In this numerical test, $T_{\text{FOM}}=T_{\text{ROM}}$, as assumed when deriving the error estimation. Since Operator Inference is a data-driven technique, the optimization error and POD projection error depend on the given snapshot data. Consequently, there is no general accuracy guarantee for the ROM simulations beyond the time duration on which the snapshots are gathered. Nonetheless, if the snapshot data represent the dynamical system behavior well, it is possible to attain accurate ROM predictions. 
Nevertheless, the predictive capabilities of the ROMs outside the training intervals are enhanced by enforcing the correct structure
in the inferred ROMs. We investigate predictions for longer time intervals in the next test.

\begin{figure}[htb!]
\centering
\begin{minipage}[ht]{0.48\linewidth}
\includegraphics[width=1\textwidth]{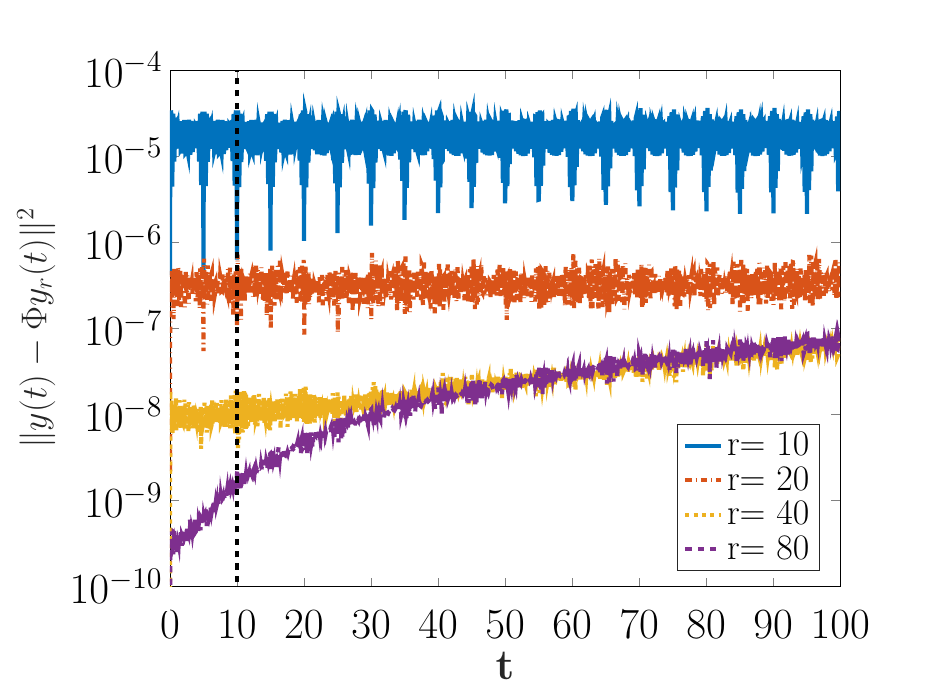}
\end{minipage}
\begin{minipage}[ht]{0.48\linewidth}
\includegraphics[width=1\textwidth]{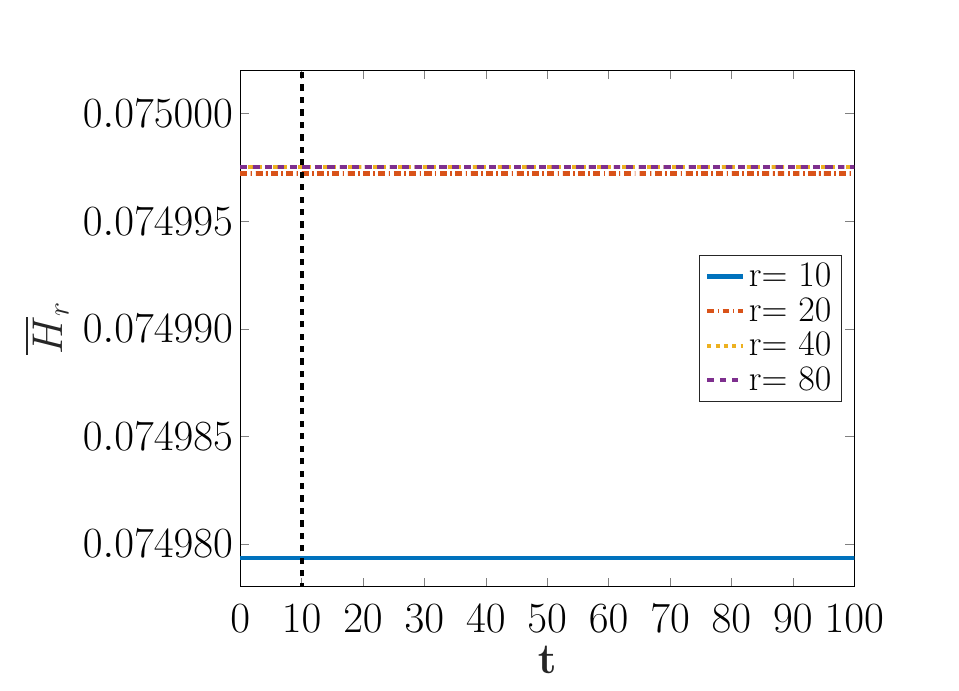}
\end{minipage}
\caption{(Wave Equation) Numerical accuracy of the $2r$-dimensional {\em GP-OpInf}: (left) time evolution of the summand of the ROM approximation error~\eqref{eq:e_approx}; (right) time evolution of the reduced-order approximate Hamiltonian energy. The vertical dashed line indicates the end of the training interval.  
}\label{fig:lin_wave_T10Tr100comp}
\end{figure}

\paragraph{Test 3. Long-term predictive capabilities of the GP-OpInf ROM}
We set $\Delta x = \Delta t = 10^{-3}$ and fix $\mu=10$ and generate snapshots from the FOM with the final time $T_{\text{FOM}}= 10$. However, we simulate the {\em GP-OpInf} ROM on a much longer interval with $T_{\text{ROM}}= 100$ to demonstrate the long-term predictive capabilities of structure-preserving ROMs. The time evolution of the ROM approximation errors, using the FOM solution as the benchmark, at $r= 10, 20, 40$ and $80$ are plotted in Figure~\ref{fig:lin_wave_T10Tr100comp} (left). The associated approximate Hamiltonian values are plotted in Figure~\ref{fig:lin_wave_T10Tr100comp} (right), where the dashed line indicates the end of the training time interval. We observe that for $r= 10$ and $20$, the errors stay at the same level over time, while for $r= 40$ and $80$, the errors increase gradually as $t$ increases. Meanwhile, since the appropriate (gradient) structure is captured in the inferred {\em GP-OpInf} ROM, the approximate Hamiltonian is conservative, which approaches the benchmark value (from the FOM) as $r$ increases.  We remark that these Hamiltonian approximations can be improved by choosing $\overline{\by}= \by^0$ and using a POD basis generated from the shifted snapshots as introduced in~\cite{gong2017structure}. 

\paragraph{Test 4. Parametric predictions away from training data}
We parameterize the initial condition $u_0(x; \mu) = h(s(x; \mu))$ and let $\mu\in [5, 15]$. To generate snapshots, 11 uniformly distributed training samples are collected from the parameter interval and the FOM is simulated at all samples with $\Delta x = \Delta t = 10^{-3}$ and $T_{\text{FOM}}=10$. Then, we compute the POD basis from the collected snapshots, learn the {\em GP-OpInf} ROM, from which we predict solutions at different test parameters over the same time interval ($T_{\text{ROM}}=T_{\text{FOM}}$). We consider three randomly selected test parameters, $\mu = 6.7, 9.5$ and $14.1$.

\begin{figure}[htb!]
\centering
\begin{minipage}[ht]{0.32\linewidth}
    \includegraphics[width=1.1\textwidth]{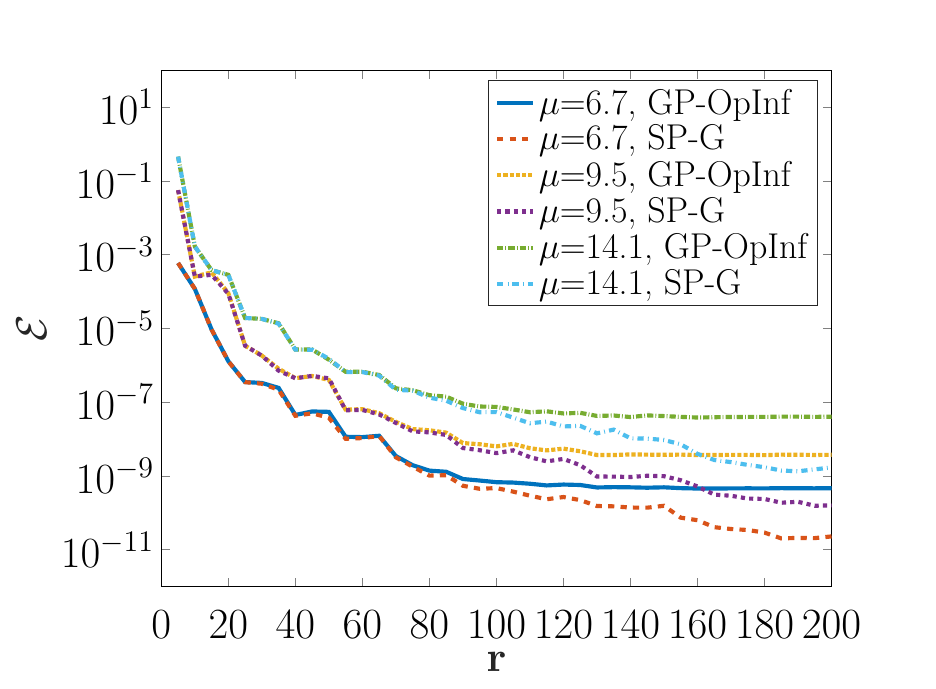}
\end{minipage}
\begin{minipage}[h]{0.32\linewidth}
	\includegraphics[width=1.1\textwidth]{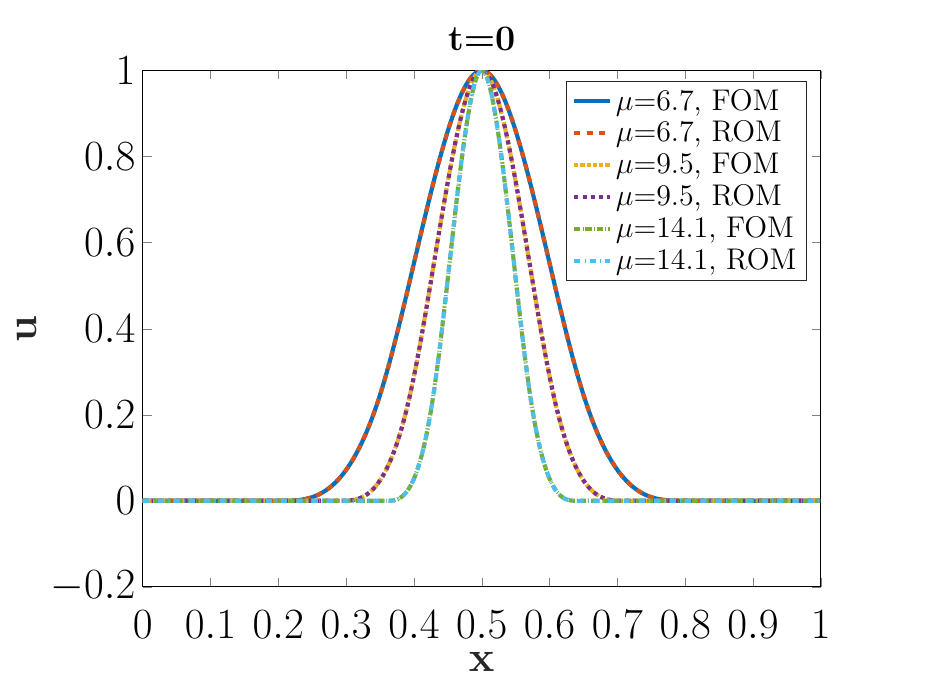}
\end{minipage}
\\
\begin{minipage}[h]{0.32\linewidth}
	\includegraphics[width=1.1\textwidth]{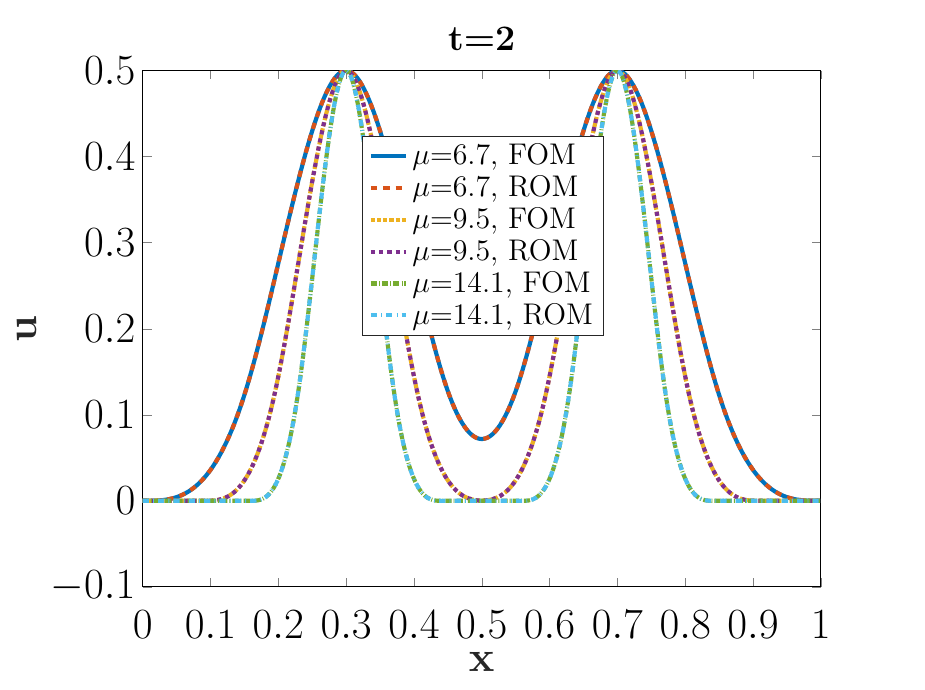}
\end{minipage}
\begin{minipage}[h]{0.32\linewidth}
	\includegraphics[width=1.1\textwidth]{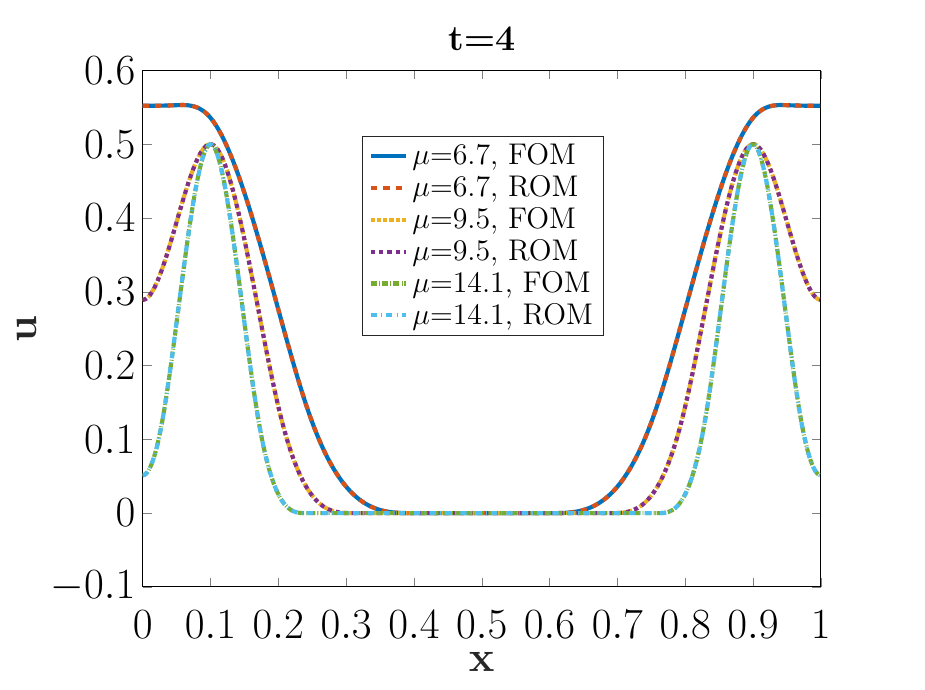}
	\end{minipage}
\begin{minipage}[h]{0.32\linewidth}
	\includegraphics[width=1.1\textwidth]{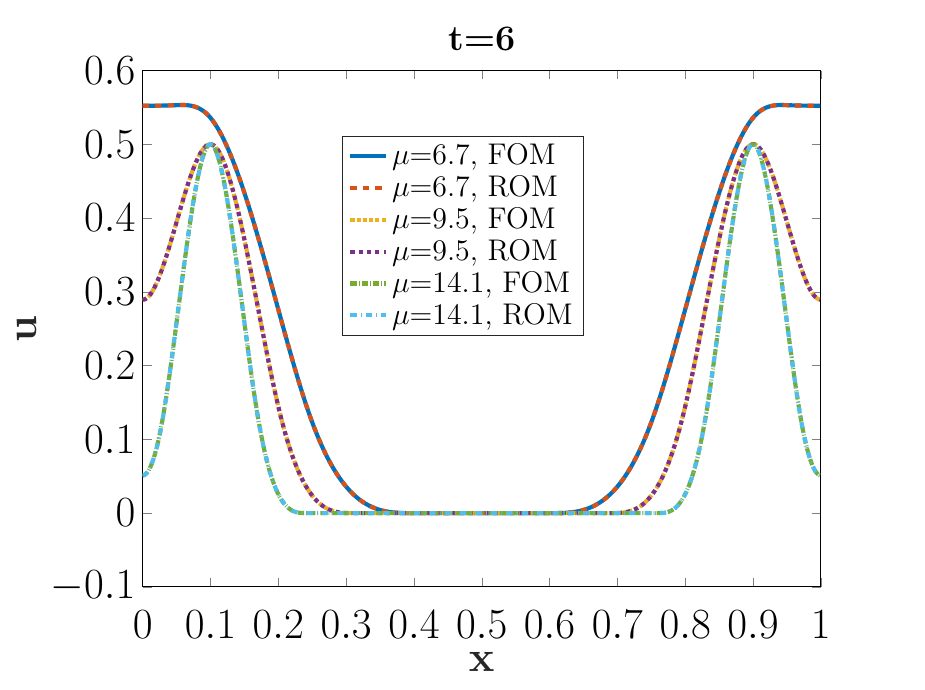}
\end{minipage}
\caption{(Wave Equation) Performance at three test parameters: (tope left) comparison of the ROM approximation errors \eqref{eq:e_approx} for the {\em GP-OpInf} with the {\em SP-G} ROM of the same dimension $2r$; (top right and bottom) the {\em GP-OpInf} ROM solution and the FOM solution at $t=0, 2, 4, 6$. 
}\label{fig:lin_wave_mu1}
\end{figure}

Figure~\ref{fig:lin_wave_mu1} (top left) compares the approximation errors of the {\em SP-G} ROM and the {\em GP-OpInf}. The figure shows that at all test parameters $\mu$, the error of the {\em GP-OpInf} ROM is close to that of {\em SP-G} ROM when $r\leq 80$. The error, however, saturates when $r$ becomes bigger and, thus, turns larger than that of the {\em SP-G} ROM. Nevertheless, such large basis dimensions produce inefficient ROMs and are of little practical value. 
We choose the more realistic $r=40$ case and plot the benchmark FOM solution and the ROM solution in Figure~\ref{fig:lin_wave_mu1} (top right and bottom); only a few snapshots are shown for illustration. In each case, the {\em GP-OpInf} ROM achieves accurate approximations close to the FOM solution. 

\begin{figure}[htb!]
\centering
\begin{minipage}[ht]{0.32\linewidth}
\includegraphics[width=1.1\textwidth]{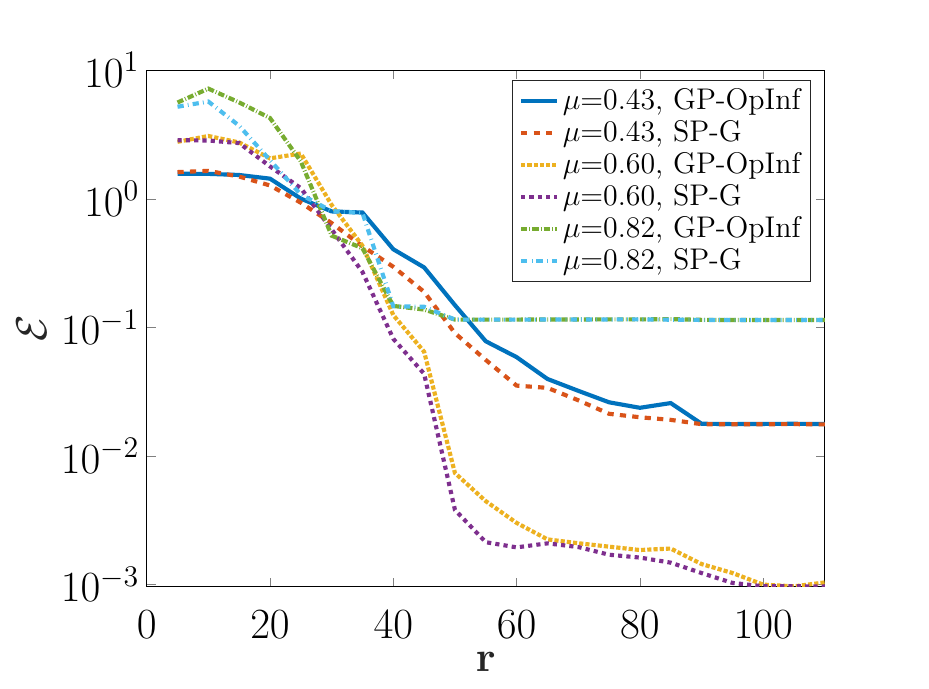}
\end{minipage}
\begin{minipage}[h]{0.32\linewidth}
	\includegraphics[width=1.1\textwidth]{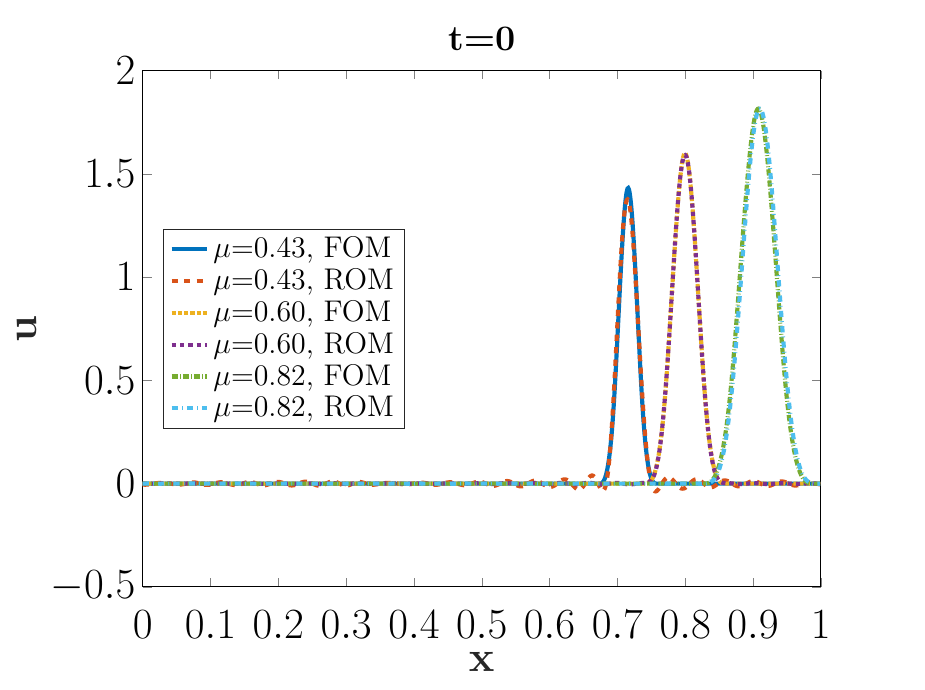}
\end{minipage}
\\
\begin{minipage}[h]{0.32\linewidth}
	\includegraphics[width=1.1\textwidth]{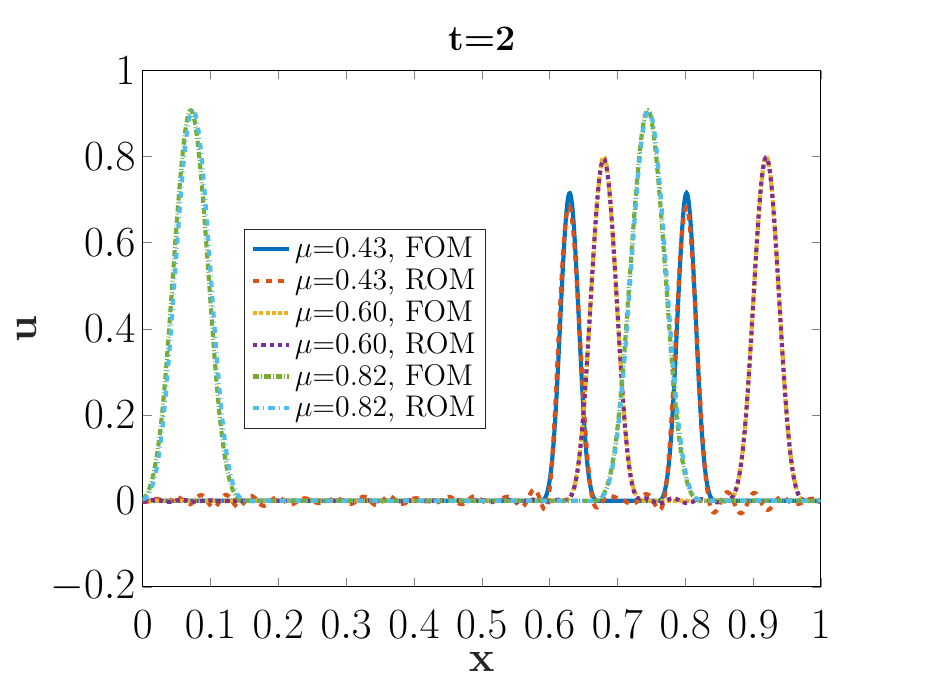}
\end{minipage}
\begin{minipage}[h]{0.32\linewidth}
	\includegraphics[width=1.1\textwidth]{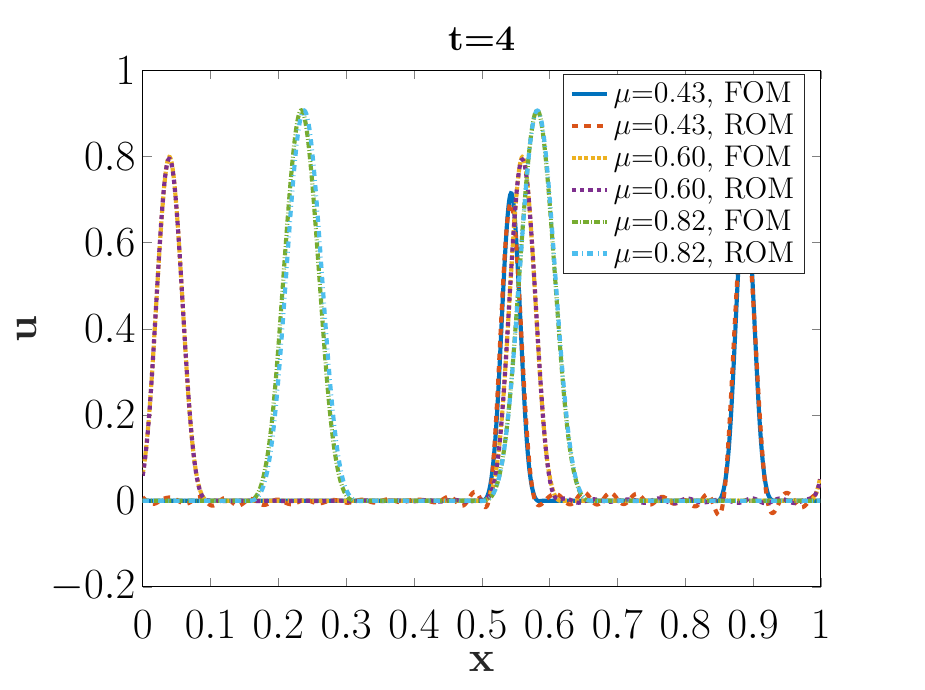}
\end{minipage}
\begin{minipage}[h]{0.32\linewidth}
	\includegraphics[width=1.1\textwidth]{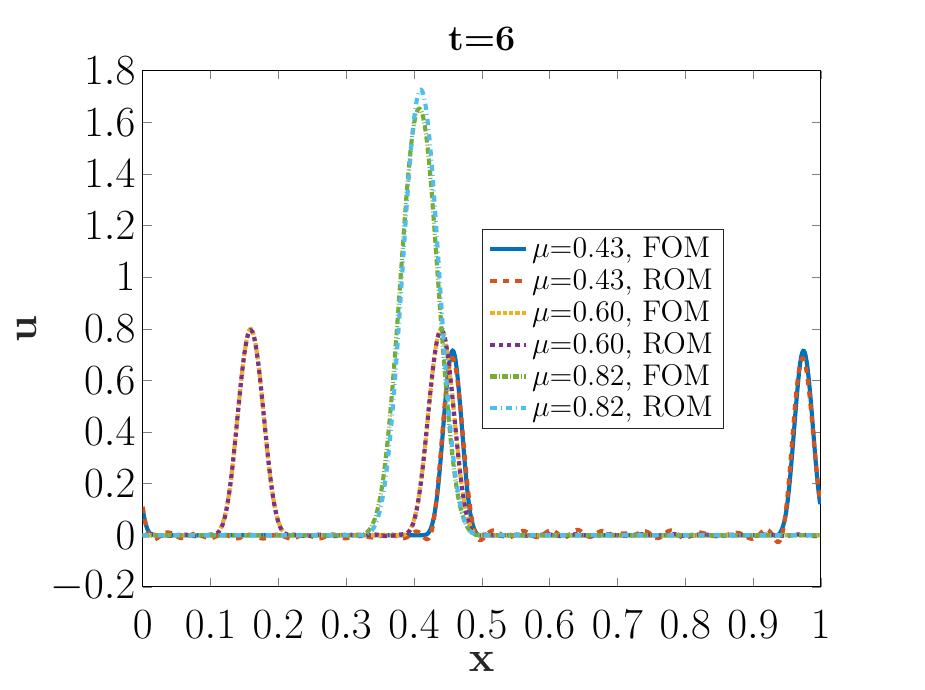}
\end{minipage}

\caption{(Wave Equation) Performance at three test parameters: (top left) comparison of the ROM approximation errors \eqref{eq:e_approx} for the {\em GP-OpInf} with the {\em SP-G} ROM of the same dimension $2r$; (top right and bottom) the {\em GP-OpInf} ROM solution and the FOM solution at $t=0, 2, 4, 6$. 
}\label{fig:lin_wave_mu2}
\end{figure}

\Ra{We next consider a more challenging parametric problem in which both the wave speed and the initial condition are parameterized: $c(\mu) = 0.1\mu$ and $u_0(x; \mu) = (1+\mu)h(s(x;\mu))$, where $s(x; \mu) = \frac{20}{\mu} |x - \frac{1}{2} - \frac{\mu}{10}|$ and $\mu\in \mathcal{P} = [5/12, 5/6]$. Compared to the preceding test, the support of the initial pulse is thinner and its location, magnitude, and travelling speed vary with the parameter $\mu$. We infer the {\em GP-OpInf} ROM based on the training set of 11 uniformly distributed parameter samples in $\mathcal{P}$. For three randomly chosen testing samples outside the training set, $\mu = 0.43$, $0.60$ and $0.82$, we plot the approximation errors of the {\em GP-OpInf} ROM in Figure~\ref{fig:lin_wave_mu2} (top left) together with those of the {\em SP-G} ROM.  Several time snapshots of the {\em GP-OpInf} with $r=55$ and the FOM solutions are presented in Figure~\ref{fig:lin_wave_mu2} (top right and bottom), demonstrating the good performance of the {\em GP-OpInf} ROM in this case as well.}

\subsubsection{Parameterized Korteweg-de Vries (KdV) equation\label{sec: kdv}}
\noindent \indent Next, we consider the one-dimensional KdV equation
\begin{equation*}
u_t = \alpha uu_x + \nu u_{xxx}, \quad (x, t)\in \Omega\times I, 
\end{equation*}
where $u=u(x, t; \mu)$ is the parameter-dependent solution and $\mu\in P$ a parameter appearing in the initial condition, defined below alongside the boundary and initial conditions. The KdV equation can be rewritten as a non-canonical Hamiltonian system (see, e.g., \cite{celledoni2012preserving}):
\begin{equation*}
u_t = \mathcal{S} \frac{\delta \mathcal{H}}{\delta u},
\end{equation*}
where $\mathcal{S}= \partial_x$ denotes the first-order derivative operator with respect to space, and the Hamiltonian function
$\mathcal{H}(u)=\int_\Omega \left(\frac{\alpha}{6}u^3-\frac{\nu}{2}u_x^2 \right) \, {\rm d}x$. 
After taking a spatial discretization with a uniform mesh size $\Delta x$ and defining a consistent discrete Hamiltonian $\overline{H}(t) \coloneqq \Delta x H(\bu)$ with $H(\bu)= \sum_j\left[\frac{\alpha}{6}u_j^3 - \frac{\nu}{2}(u_{j+1}-u_j)^2\right]$, we obtain the semi-discrete Hamiltonian system 
\begin{flalign}
\dt{\bu} = \bA\left(\frac{\alpha}{2}\bu^2+\nu\bB\bu\right),
\label{eq: kdv_fom1}
\end{flalign}
where $\bA$ and $\bB$ are the matrices associated to the discretization of the skew-adjoint operator $\mathcal{S}$ and the second-order derivative by central differences, respectively.  
Since $\bA$ is skew-symmetric, this dynamical system conserves the discrete Hamiltonian $\overline{H}(t)$.

\paragraph{Computational setting}
We set $\Omega = [-20, 20]$, 
$\alpha= -6$, $\nu= -1$, and consider periodic boundary conditions $u(-20, t; \mu)= u(20, t; \mu)$ for any $t\in I$ and we parameterize the initial condition $u_0(x; \mu) = \text{sech}^2\left(\frac{x}{\mu}\right)$, for any $\mu \in {P}= [1,5]$. 
In the FOM~\eqref{eq: kdv_fom1}, we have  
\begin{equation*}
\bA= \frac{1}{2\Delta x}\mathbf{S}_n,
\quad \text{and} \quad 
\bB= \frac{1}{\Delta x^2} \mathbf{L}_n. 
\end{equation*} 
The mesh size $\Delta x = 0.01$ (correspondingly, $n=4,000$) is used in all simulations. 
For time integration, we use the average vector field (AVF) method \Rb{\cite{celledoni2012preserving}} together with a Picard iteration to solve the nonlinear system.

\begin{figure}[htb!]
\centering
\begin{minipage}[ht]{0.48\linewidth}
\includegraphics[width=1\textwidth]{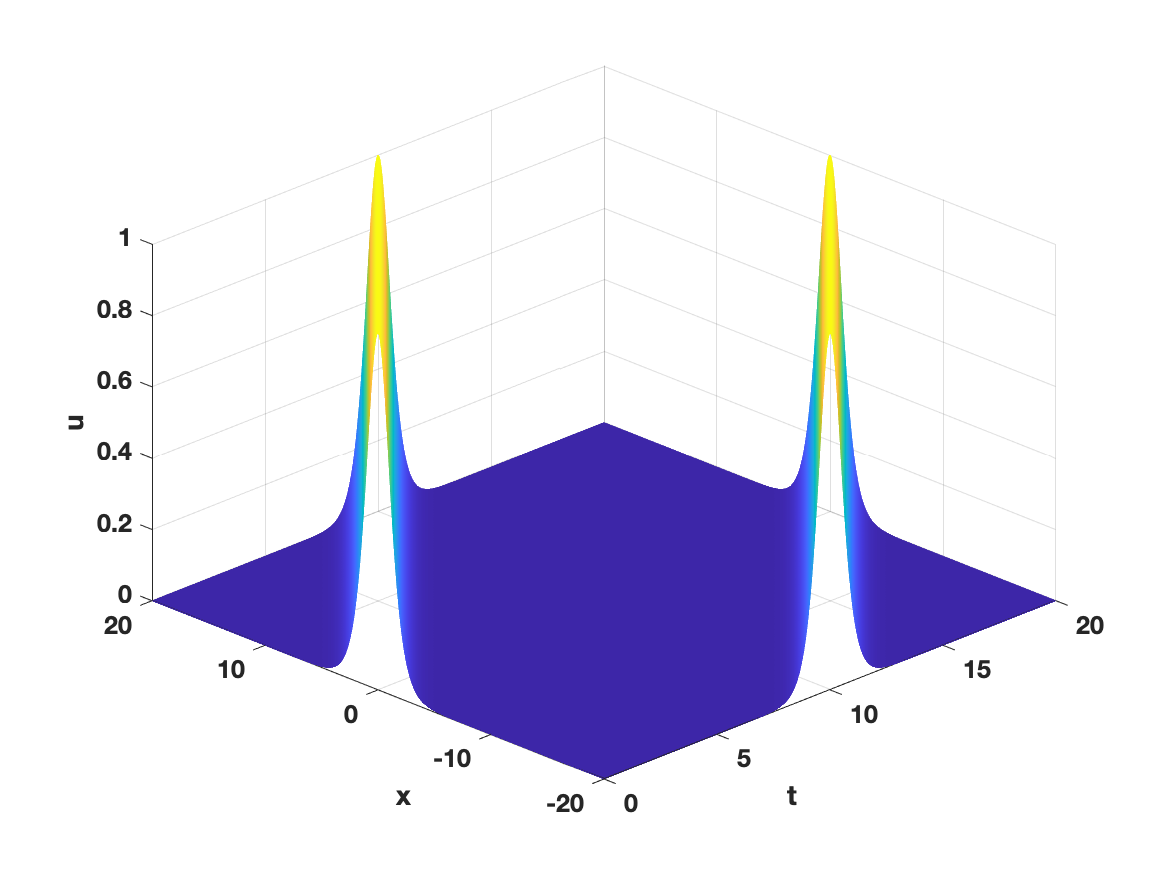}
\end{minipage}
\begin{minipage}[ht]{0.48\linewidth}
\includegraphics[width=1\textwidth]{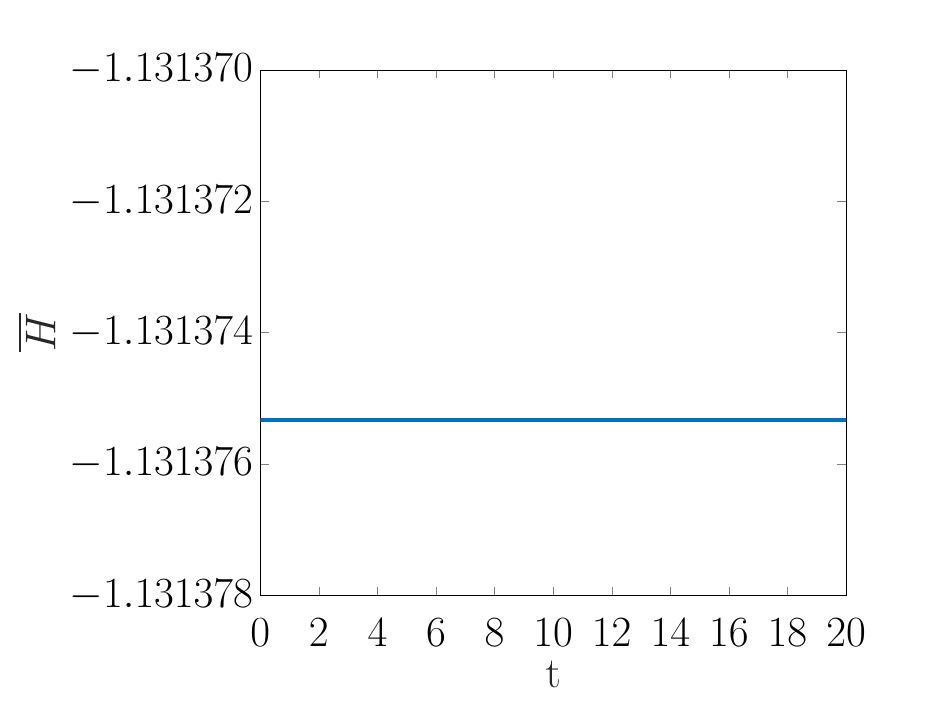}
\end{minipage}
\caption{(KdV Equation) Full-order model simulation up to final time $T_{\text{FOM}} = 20$: time evolution of $\bu$ (left) and approximate Hamiltonian $\overline{H}(t)$ (right).
}\label{fig:kdv_fom}
\end{figure}

When $T_{\text{FOM}} = 20$, $\Delta t= 0.01$ and $\mu = \sqrt{2}$, the time evolution of the full-order state and approximate Hamiltonian $\overline{H}(t)$ are shown in Figure \ref{fig:kdv_fom}. We observe that the system is conservative, with the Hamiltonian staying around $\overline{H}(t) \approx -1.13$ during the simulation. 

\begin{figure}[htb!]
\centering
\begin{minipage}[ht]{0.48\linewidth}
\includegraphics[width=1\textwidth]{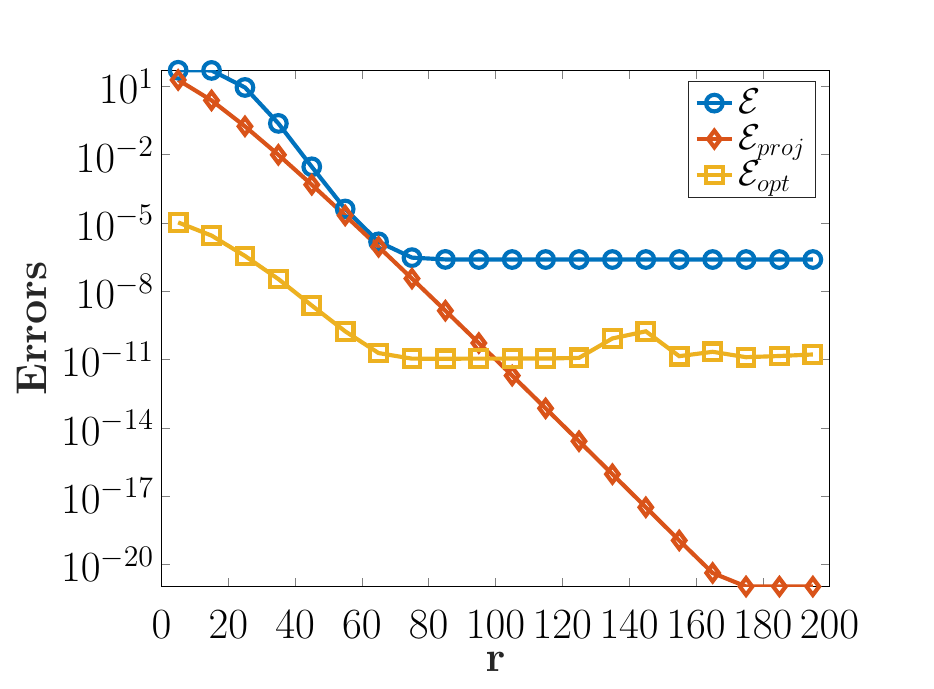}
\end{minipage}
\begin{minipage}[ht]{0.48\linewidth}
\includegraphics[width=1\textwidth]{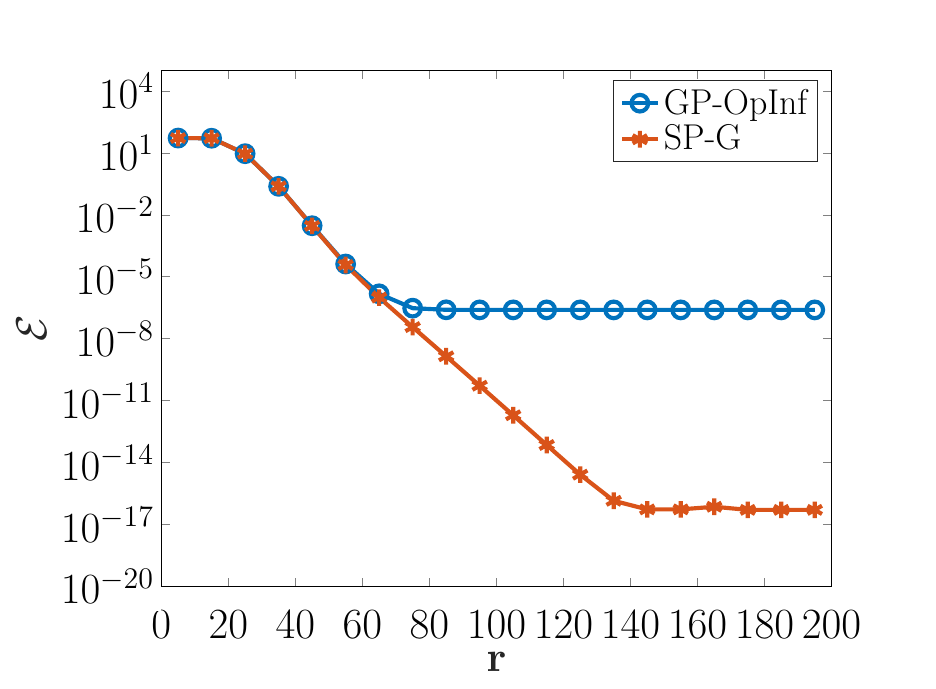}
\end{minipage}
\caption{(KdV Equation) Numerical errors of the $r$-dimensional {\em GP-OpInf} ROM when $T_{\text{FOM}}=T_{\text{ROM}}=20$ and $\Delta t = 2.5\times 10^{-3}$: (left) the ROM approximation error \eqref{eq:e_approx} together with POD projection error \eqref{eq:e_proj} and optimization error \eqref{eq:e_opt}; (right) comparison of the ROM approximation error of the {\em GP-OpInf} ROM  with the {\em SP-G} ROM. 
}\label{fig:kdv_T20Tr20err}
\end{figure}

\begin{figure}[htb!]
\centering
\begin{minipage}[ht]{0.48\linewidth}
\includegraphics[width=1\textwidth]{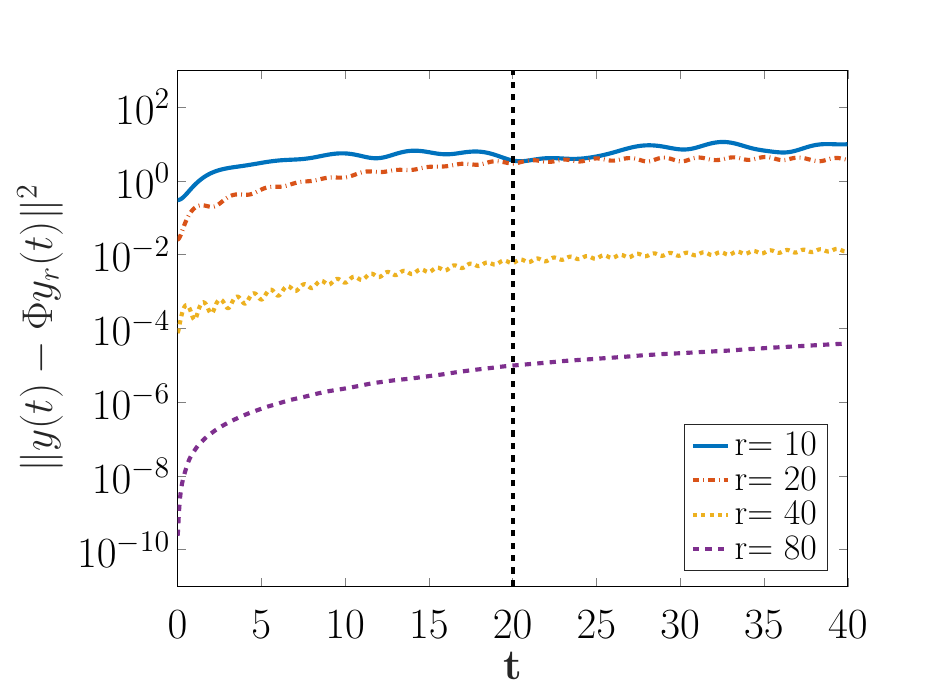}
\end{minipage}
\begin{minipage}[ht]{0.48\linewidth}
\includegraphics[width=1\textwidth]{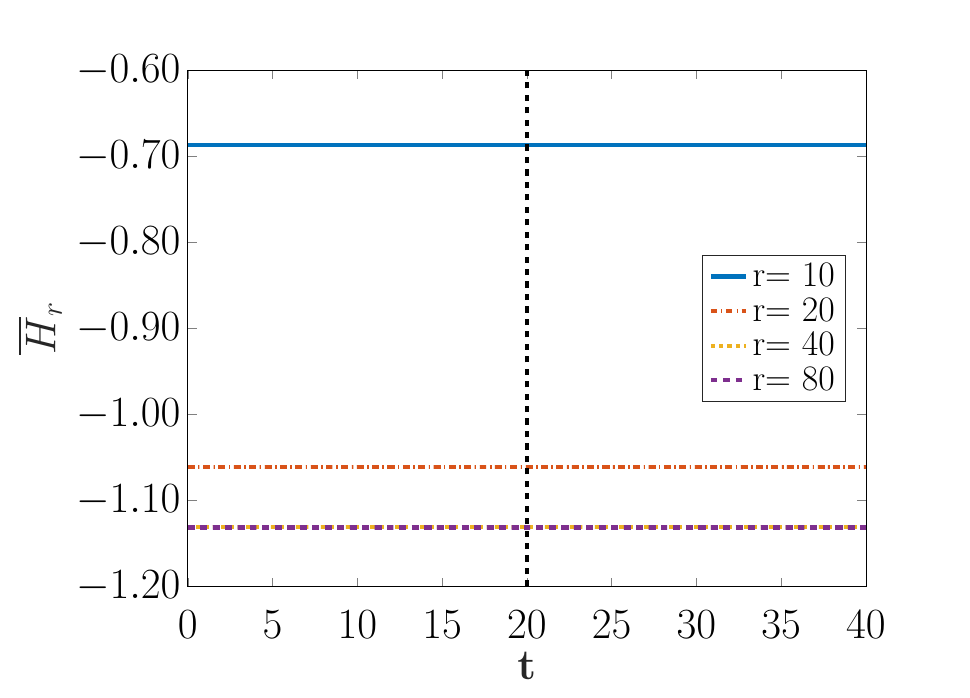}
\end{minipage}
\caption{(KdV Equation) Numerical accuracy of the $r$-dimensional {\em GP-OpInf} ROM: (left) time evolution time evolution of the summand of the ROM approximation error~\eqref{eq:e_approx}; (right) time evolution of the values of reduced-order approximate Hamiltonian functions. The dashed line indicates the end of the training interval. 
}\label{fig:kdv_T20Tr40comp}
\end{figure}

\paragraph{Test 1. Illustration of error estimation}
We demonstrate the estimation of the {\em GP-OpInf} ROM approximation error in a non-parametric case with $\mu= \sqrt{2}$, and compare the approximation error with the {\em SP-G} ROM of the same dimension. 
To verify the error estimation, we choose a small time step $\Delta t=2.5\times 10^{-3}$, let $T_{\text{FOM}}=T_{\text{ROM}}=20$ and vary the dimension of the ROM. 
Figure~\ref{fig:kdv_T20Tr20err} (left) shows the three error measures from equations~\eqref{eq:e_approx}--\eqref{eq:e_opt}. The ROM approximation error decays monotonically until $r>80$ and then levels off. The optimization error follows a similar trend, yet the projection error continues to decay until $r\approx 170$.
Figure~\ref{fig:kdv_T20Tr20err} (right) compares the {\em GP-OpInf} ROM approximation error with that of the {\em SP-G} ROM. We observe  that {\em SP-G} and {\em GP-OpInf} achieve the same accuracy for $r<70$, which covers the model dimension where most practical ROMs would be selected from, whereas for larger $r$ the {\em SP-G} ROM performs better. 
\paragraph{Test 2. Long-term predictive capabilities of the GP-OpInf ROM}
We set $\Delta t = 10^{-2}$ and fix $\mu=\sqrt{2}$ and generate snapshots from FOM with the final time $T_{\text{FOM}}= 20$. We simulate the {\em GP-OpInf} ROM 100\% past the training interval, so $T_{\text{ROM}}= 40$. 
Figure~\ref{fig:kdv_T20Tr40comp} (left) shows the time evolution of the ROM approximation errors at $r= 10, 20, 40$ and $80$, where the FOM solution is the benchmark. The associated approximate Hamiltonian values are plotted in Figure~\ref{fig:kdv_T20Tr40comp} (right), where the dashed line indicates the end of the training interval. The errors increase gradually as $t$ increases due to usual error accumulation. Since the appropriate (gradient) structure is captured in the inferred {\em GP-OpInf} ROMs, the approximate Hamiltonian functions are constant and  approach the benchmark value for $r$ that is sufficiently large.  

\paragraph{Test 3. Parametric predictions away from training data}
We parameterize the initial condition as $u_0(x; \mu) = \text{sech}^2\left(\frac{x}{\mu}\right)$ for $\mu\in {P} = [1, 5]$. To generate snapshots, 9 training samples are uniformly collected from the interval and the FOM is simulated at all samples with $\Delta t = 10^{-2}$ and $T_{\text{FOM}}=20$. We compute the POD basis from the collection of snapshot matrices. 

\begin{figure}[t!]
\centering
\begin{minipage}[ht]{0.32\linewidth}
\includegraphics[width=1.1\textwidth]{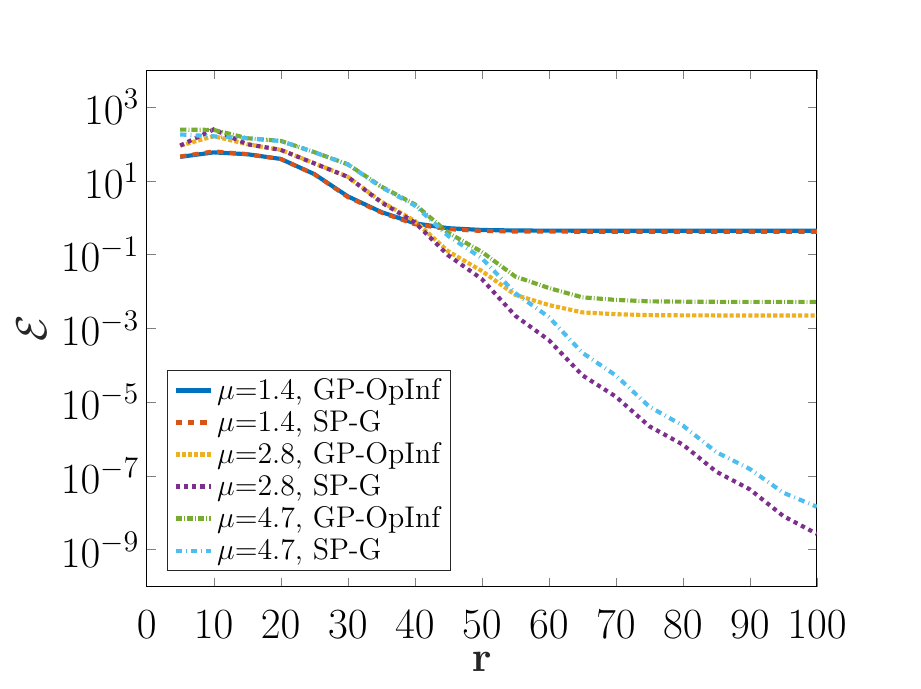}
\end{minipage}
\begin{minipage}[h]{0.32\linewidth}
	\includegraphics[width=1.1\textwidth]{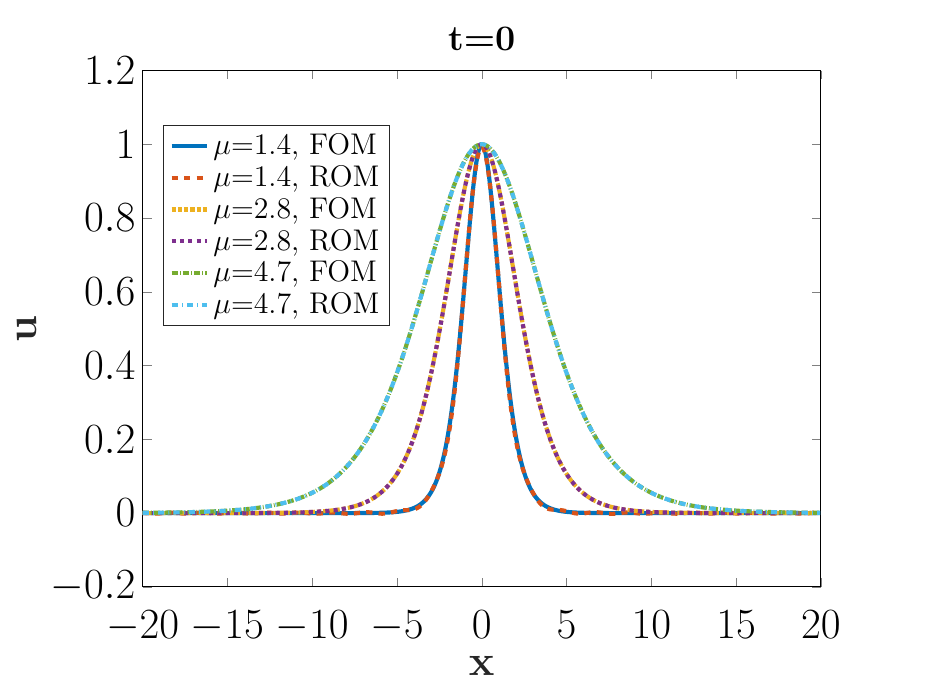}
\end{minipage}
\\
\begin{minipage}[h]{0.32\linewidth}
	\includegraphics[width=1.1\textwidth]{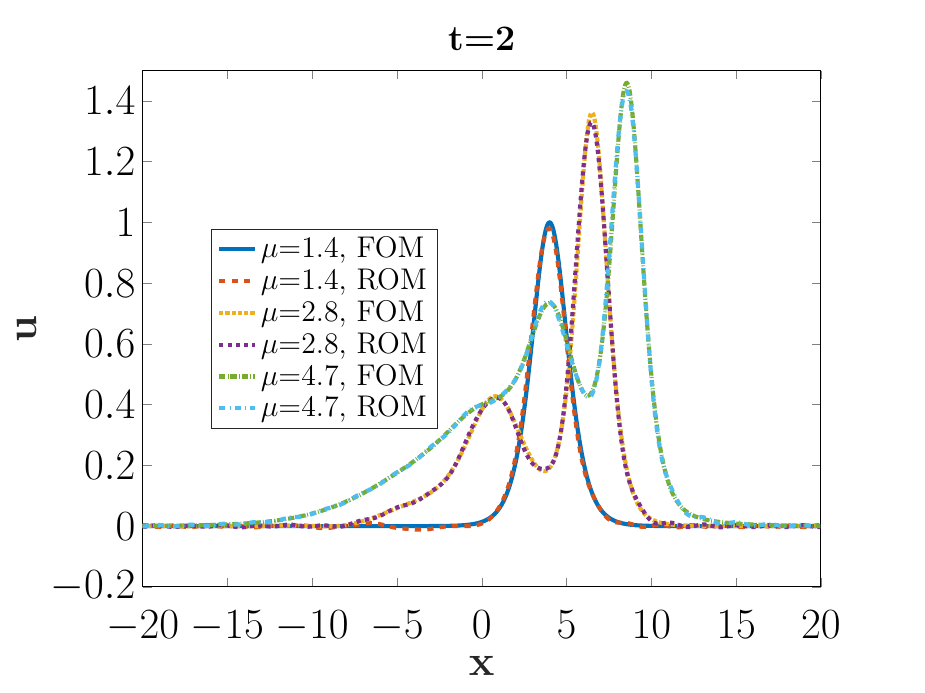}
\end{minipage}
\begin{minipage}[h]{0.32\linewidth}
	\includegraphics[width=1.1\textwidth]{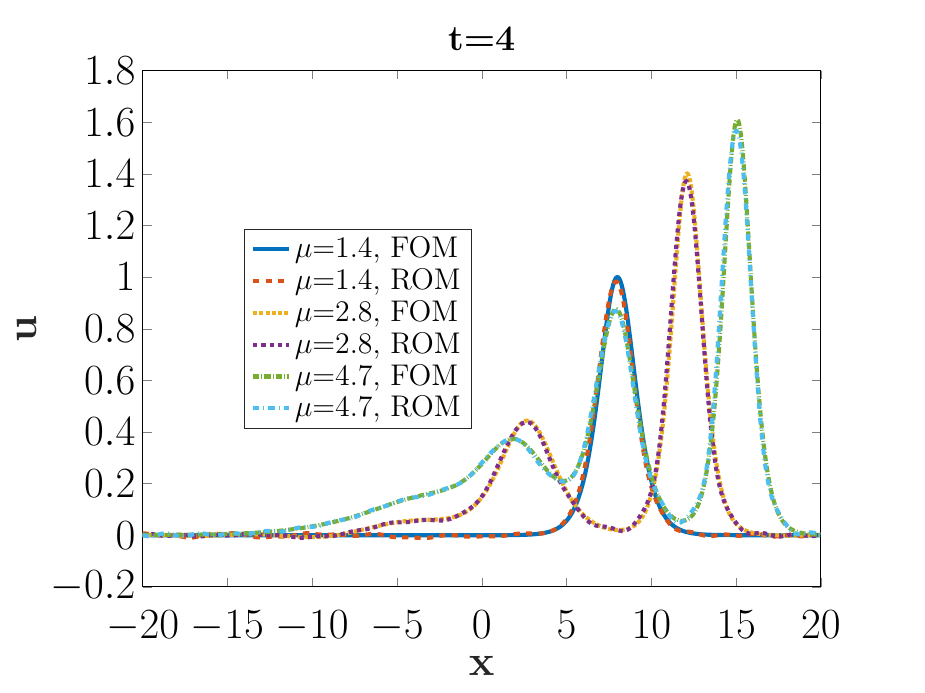}
\end{minipage}
\begin{minipage}[h]{0.32\linewidth}
	\includegraphics[width=1.1\textwidth]{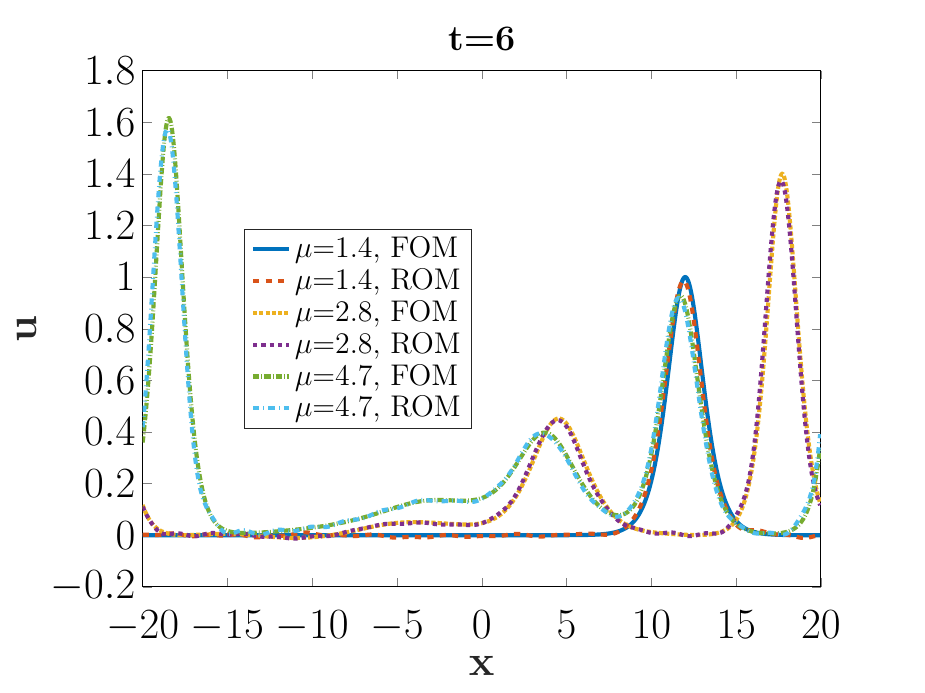}
\end{minipage}
\caption{(KdV Equation) Performance of the $r$-dimensional {\em GP-OpInf} ROM at several test parameters: (top left) comparison of the ROM approximation error \eqref{eq:e_approx} for the {\em SP-G} ROM of the same dimension $r$; (top right and bottom) snapshots of the {\em GP-OpInf} ROM solution and the FOM solution at several time instances. 
}\label{fig:kdv_mu1}
\end{figure}

Next, we infer the {\em GP-OpInf} ROM and use it to make predictions at the test parameters, where we keep the time interval fixed ($T_{\text{ROM}}=T_{\text{FOM}}$). We evaluate the ROMs accuracy at three randomly selected test parameters $\mu = 1.4, 2.8$ and $4.7$. 
Figure~\ref{fig:kdv_mu1} (top left) shows a comparison of the ROM approximation errors \eqref{eq:e_approx} of the {\em SP-G} and the {\em GP-OpInf} ROMs. At the test parameters, the error of the {\em GP-OpInf} ROM is close to the {\em SP-G} ROM when $r\leq 50$ (which is where most practical ROM model dimensions would be selected from), yet the error saturates when $r$ gets bigger and thus becomes larger than that of {\em SP-G}. 
We fix $r=40$ and plot the FOM solution and the {\em GP-OpInf} ROM solution at several selected time instances in Figure~\ref{fig:kdv_mu1} (top right and bottom). This shows that the ROM produces accurate approximations of the FOM solution at these test parameters. 

\Ra{
\subsubsection{Three-dimensional Linear Elastic Cantilever Plate}
We consider the 3D linear elastic cantilever plate problem from \cite{Gruber2023canonical}. 
The governing equation can be recast in the conservative Hamiltonian form
\[
\left[
\begin{array}{c}
\dot{\mathbf{u}}\\
\dot{\mathbf{v}}
\end{array}
\right] =
\left[
\begin{array}{cc}
0 & \mathbf{I_n} \\
-\mathbf{I_n} & 0
\end{array}
\right]
\left[
\begin{array}{cc}
\mathbf{K} & 0 \\
0 & \mathbf{M^{-1}}
\end{array}
\right]
\left[
\begin{array}{c}
{\mathbf{u}}\\
{\mathbf{v}}
\end{array}
\right],
\]
where $\bu$ is the displacement, $\bv= \mathbf{M}\dt{\bu}$ is the momentum, $\mathbf{M}$ and $\mathbf{K}$ are the mass and stiffness matrix obtained from a finite element discretization.
The Hamiltonian of the semi-discrete system is $H(\bu, \bv) = \frac{1}{2}\left(\mathbf{u^{\intercal}Ku} + \mathbf{v^{\intercal}M^{-1}v}\right)$.
Following \cite{Gruber2023canonical}, we consider a rectangular steel plate of size $0.2m \times 0.2m \times 0.03m$, whose ``left" side is set in place by a clamp, and the ``right" side is given an initial upward velocity $(0, 0, 100)$. 
For the spatial discretization, we use $\Delta x=\Delta y=\Delta z=0.01$ (correspondingly $n=5,292$), which is identical to the numerical setting in \cite{Gruber2023canonical} \footnote{We use the data from the github repository: \url{https://github.com/ikalash/HamiltonianOpInf}. We note that regularization is also used in their code to infer Hamiltonian reduced-order operators via optimization.}. We simulate the full-order model over the time interval $[0, 0.02]$ with the implicit midpoint method, with time step size $\Delta t = 1\times10^{-4}$ to collect snapshots.
Figure~\ref{figure:3Dplateplots} shows the 3D surface of the plate at time instances $t=0, 0.005, 0.0125$ and $0.02$.  

\begin{figure}[t!]
    \centering
    \includegraphics[width=0.24\textwidth]{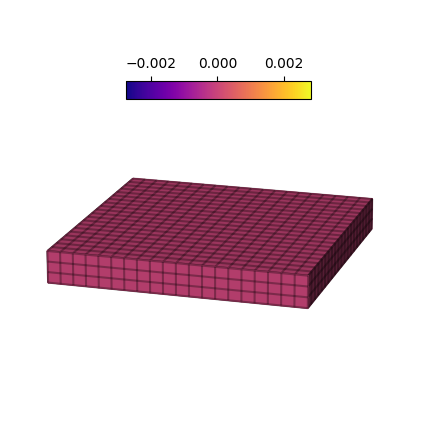} 
    \includegraphics[width=0.24\textwidth]{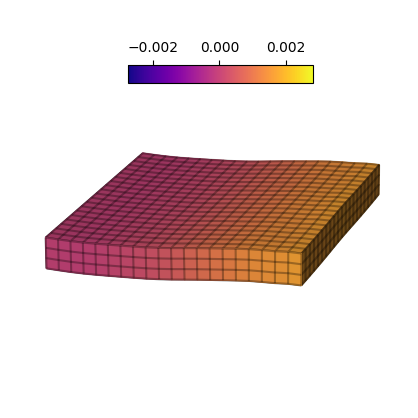} 
    \includegraphics[width=0.24\textwidth]{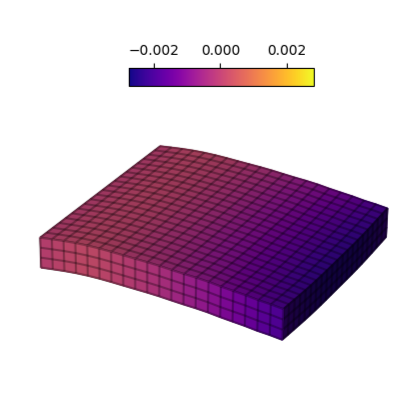}
    \includegraphics[width=0.24\textwidth]{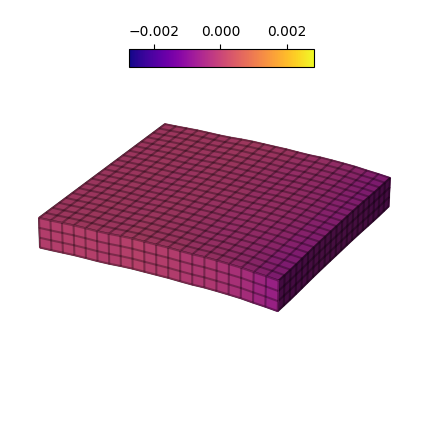}
    \caption{(3D Cantilever Plate) Plate's surface plots at $t=0, 0.005, 0.0125$, and $0.02$ (from left to right), in which the color bar indicates the displacement in the vertical direction. }
    \label{figure:3Dplateplots}
\end{figure}

\paragraph{Test 1. Effects of Regularization} Similar to the linear wave case, we obtain the reduced basis $\mathbf{\Phi} = \left[
\begin{array}{cc}
\mathbf{\Phi_u} & 0 \\
0 & \mathbf{\Phi_v}
\end{array}
\right]$ with $r_1=r_2=r$, and simulate the {\em GP-OpInf}, {\em GP-OpInf-V}, and {\em GP-OpInf-P} (with $\alpha=2$ and $c_0 = 10^{-13}$ in the definition of $\epsilon$ \eqref{eq:epsilon}) ROMs, respectively, within the training data regime. 
The corresponding approximation errors are shown in Figure~\ref{figure:3DPlateROMerrors} as $r$ varies. It is evident that \emph{GP-OpInf} surpasses the other two with its accuracy. Especially, as $r$ increases beyond $35$, the errors of \emph{GP-OpInf-V} and \emph{GP-OpInf-P} grow uncontrollably, while \emph{GP-OpInf} maintains a small error. 
\begin{figure}[htb]
\centering
\begin{minipage}[ht]{\linewidth}
\centering
\begin{minipage}[h]{.48\linewidth}
    \includegraphics[width=\linewidth]{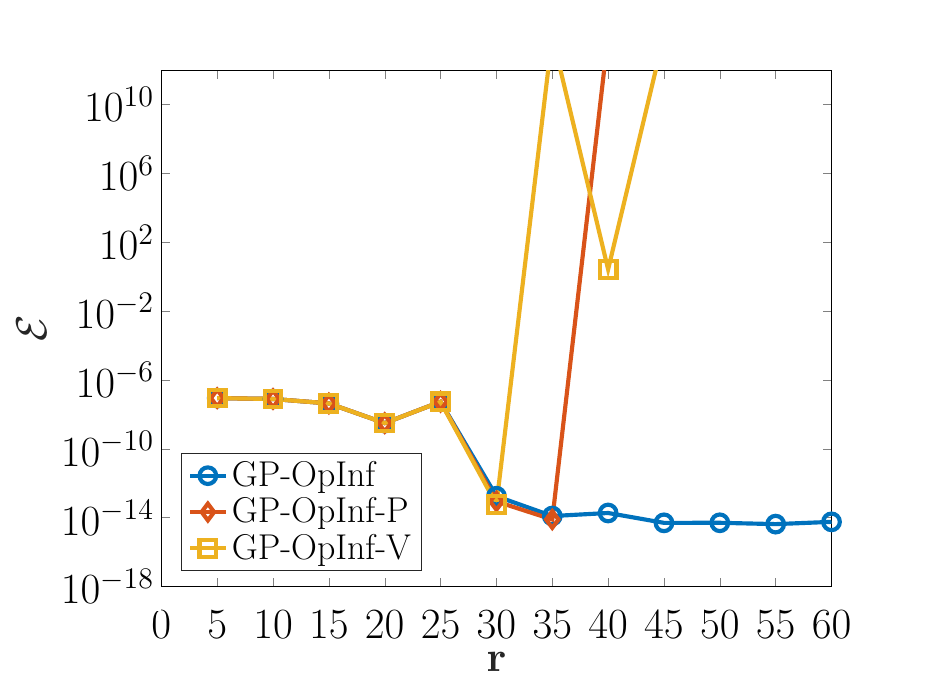}
\end{minipage}
\end{minipage}
\caption{(3D Cantilever Plate) Approximation errors for {\em GP-OpInf}, {\em GP-OpInf-V} and {\em GP-OpInf-P} as $r$ increases.}
\label{figure:3DPlateROMerrors}
\end{figure}

\paragraph{Test 2. Illustration of error estimation} We now investigate the \emph{GP-OpInf} approximation error, projection error, and optimization error in the training data regime. This comparison is illustrated in Figure~\ref{figure:3DPlateGPerrors} (left), while $r$ varies from $5$ to $60$, increasing incrementally by $5$. The optimization error stagnates for $r\geq 30$ at a level below $10^{-6}$ which for many applications is sufficiently low.

\begin{figure}[htb!]
\centering
\begin{minipage}[ht]{.32\linewidth}
\centering
    \includegraphics[width=1.1\linewidth]{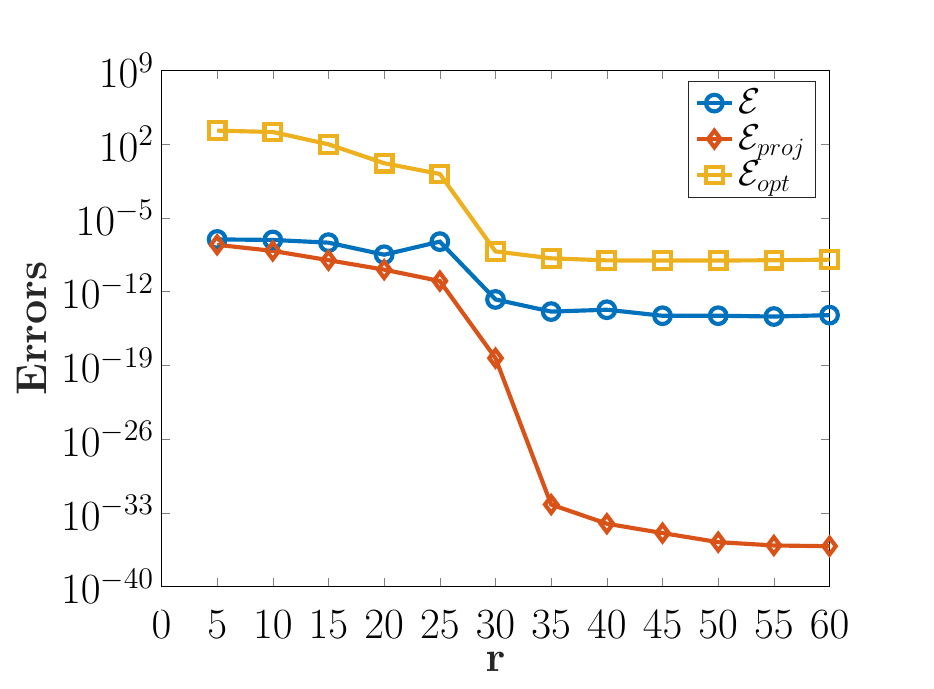}
\end{minipage}
\begin{minipage}[h]{.32\linewidth}
    \includegraphics[width=1.1\linewidth]{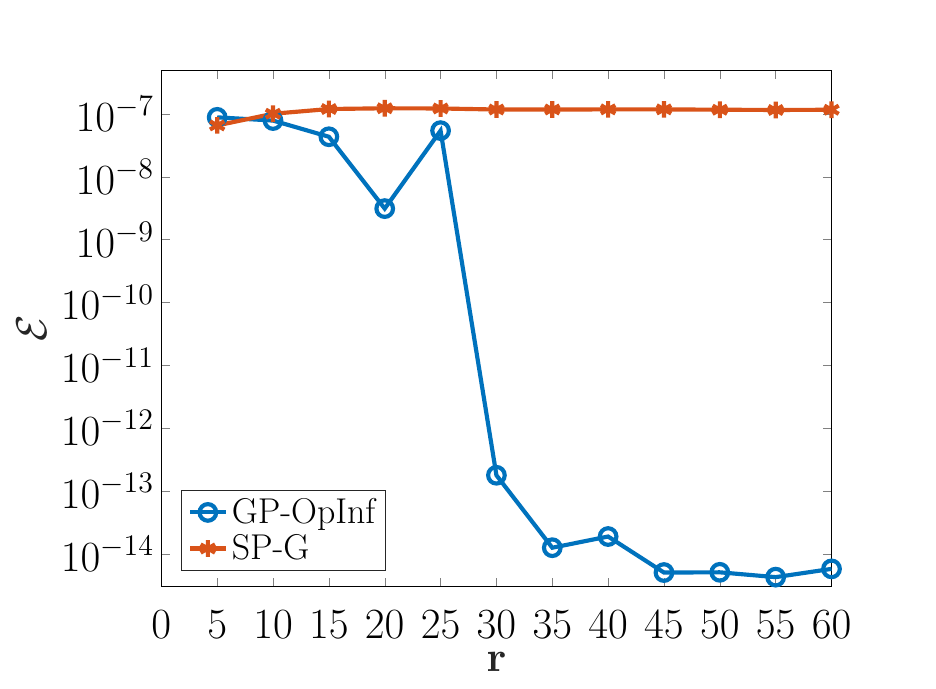}
\end{minipage}
\begin{minipage}[h]{.32\linewidth}
    \includegraphics[width=1.1\linewidth]{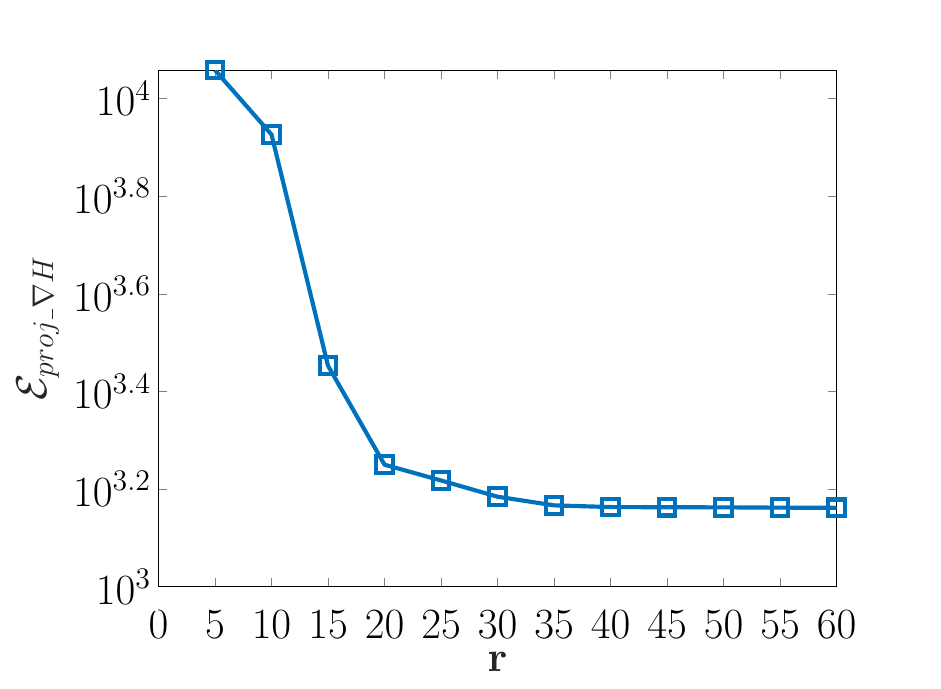}
\end{minipage}
\caption{(3D Cantilever Plate) Numerical errors of the {\em GP-OpInf}: (left) the ROM approximation error \eqref{eq:e_approx} together with POD projection error \eqref{eq:e_proj} and optimization error \eqref{eq:e_opt}; (Middle) comparison with approximation errors of the  {\em SP-G}; (Right) the projection error of $\nabla H$ onto the reduced basis space.}
\label{figure:3DPlateGPerrors}
\end{figure}

We further compare the accuracy of the nonintrusive \emph{GP-OpInf} with the intrusive \emph{SP-G} in Figure~\ref{figure:3DPlateGPerrors} (middle). 
Despite increasing $r$, the error of \emph{SP-G} does not improve, whereas the \emph{GP-OpInf} error decays to almost machine precision. A potential explanation for this phenomenon is discussed in Remark~\ref{re:2}. 
In light of that discussion, we calculate the projection errors of $\nabla H$-snapshots onto the reduced basis subspace,
\[
\mathcal{E}_{{\text{proj}}\_{\nabla H}}(r) = \frac{T \Delta A}{N}\sum_{j=1}^N\lvert\lvert \nabla H(\mathbf{y^j}) - \mathbf{\Phi\Phi^{\intercal}}\nabla H(\mathbf{y^j}) \rvert\rvert^2,
\]
and plot it in Figure~\ref{figure:3DPlateGPerrors}~(right). 
We observe that the magnitude of this projection error is on the order of $10^{3}$ even when $r=60$, which could limit the accuracy of {\em SP-G}. However, it does not affect the {\em GP-OpInf} as indicated by the error bound in Theorem~\ref{thm:1}.

We remark that we also tested a similar parametric problem, in which the initial velocity vector was parameterized as $(0, 0, \mu)$, with the parameter $\mu\in \mathcal{P} \coloneqq [80, 120]$. The {\em GP-OpInf} ROM, inferred based on a training set of $11$ uniformly distributed parameter samples, accurately predicted the solution for randomly selected parameters from $\mathcal{P}$. Since this phenomenon closely resembles our observations in the linear wave test case, we omit further discussion on it here.
}


\subsection{Dissipative PDEs}
\label{sec:num_dissipative}
We consider the Allen-Cahn equation to assess the {\em GP-OpInf} ROM's performance on a dissipative system.

\subsubsection{One-dimensional Allen-Cahn equation}
\label{sec:1DAC}
Consider the one-dimensional Allen-Cahn equation 
\begin{flalign*} 
u_t = \epsilon^2 u_{xx}+(u-u^3), \quad (x, t)\in \Omega\times I \text{ and } \mu\in P, 
\end{flalign*}
where $u=u(x, t; \mu)$ with $\mu$ a parameter in the initial condition. The equation can be recast to the form \eqref{eq:h-pde}, that is,
$ u_t = - \frac{\delta \mathcal{H}}{\delta u}$,  
with the Lyapunov function 
$\mathcal{H}(u) = \int_\Omega \left(\frac{\epsilon^2}{2} u_x^2 -\frac{1}{2}u^2  +\frac{1}{4}u^4 \right)\, {\rm d} x$ (see, e.g., \cite{celledoni2012preserving}). 
After applying a uniform spatial discretization with the mesh size $\Delta x$ and defining a consistent discrete Lyapunov function $\overline{H}(t)\coloneqq \Delta x H(\bu)$ with $H(\bu) = \sum_j [ \frac{\epsilon^2}{2\Delta x^2} (u_{j+1}-u_j)^2 -\frac{1}{2}u_j^2  +\frac{1}{4}u_j^4 ]$, we obtain the semi-discrete system
\begin{flalign}
    \dt{\mathbf{u}} = -\mathbf{I}_n(  -\mathbf{Au}-\mathbf{u}+\mathbf{u}^3 ),
    \label{eq:ac-r}
\end{flalign}
where $\bI_n$ is the $n\times n$ identity matrix, $\bA$ is associated with a discrete, one-dimensional, second-order differential operator, and $\mathbf{u}^3$ is the component-wise cubic function of the vector $\mathbf{u}$. Since $-\bI_n$ is negative definite, the Lyapunov function decreases in time.

\begin{figure}[htb!]
\centering
\begin{minipage}[ht]{0.48\linewidth}
\includegraphics[width=1\textwidth]{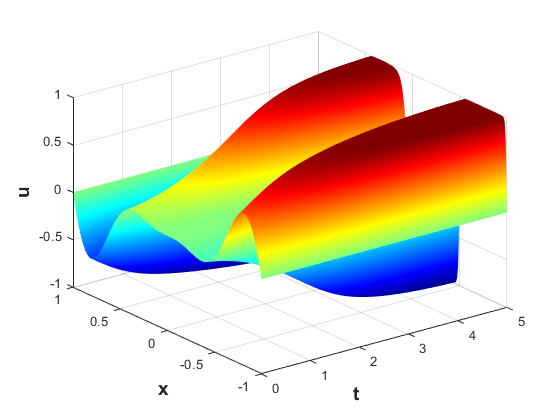}
\end{minipage}
\begin{minipage}[ht]{0.48\linewidth}
\includegraphics[width=1\textwidth]{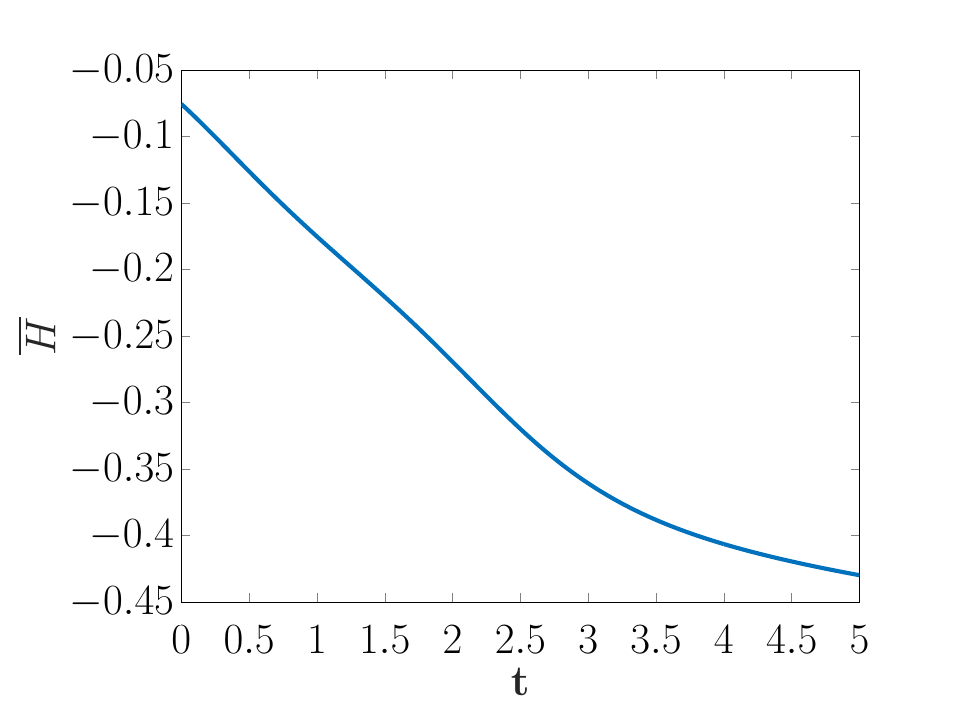}
\end{minipage}
\caption{(1D Allen-Cahn Equation) Full-order model simulation with final time $T_{\text{FOM}} = 5$: time evolution of the state $\bu$ (left) and the approximate Lyapunov function $\overline{H}(t)$ (right).
}\label{fig:1D_AC_FOM_3Dsolution_energy}
\end{figure}

\paragraph{Computational setting} 
Set $\Omega = [-1, 1]$, 
the interface parameter $\epsilon=0.01$, and consider periodic boundary conditions $u(-1,t; \mu)=u(1,t; \mu)$ and $\frac{\partial u}{\partial \mathbf{n}}(-1,t; \mu) = \frac{\partial u}{\partial \mathbf{n}}(1,t; \mu)$ for $t\in I$. The initial condition is $u_0(x; \mu)=\mu x^2 \sin(2\pi x)$, for $\mu\in {P}= [0.2,2]$ (see \cite{zhao2020solving} for numerical simulations of a non-parametric configuration). In the FOM~\eqref{eq:ac-r}, we have 
\begin{equation*}
\bA= \frac{\epsilon^2}{\Delta x^2}\mathbf{L_n}.
\end{equation*}
In all full-order simulations, we choose the mesh size $\Delta x= 10^{-3}$ and use the AVF method for time integration with a time step size $\Delta t$. To solve the nonlinear systems of equations, we apply Picard iteration. Figure~\ref{fig:1D_AC_FOM_3Dsolution_energy}  shows the time evolution of the full-order state and discrete Lyapunov function $\overline{H}(t)$ where we set $\mu =1$, $T_{\text{FOM}} = 5$ and $\Delta t=10^{-3}$.

\begin{figure}[htb!]
\centering
\begin{minipage}[ht]{0.48\linewidth}
\includegraphics[width=1\textwidth]{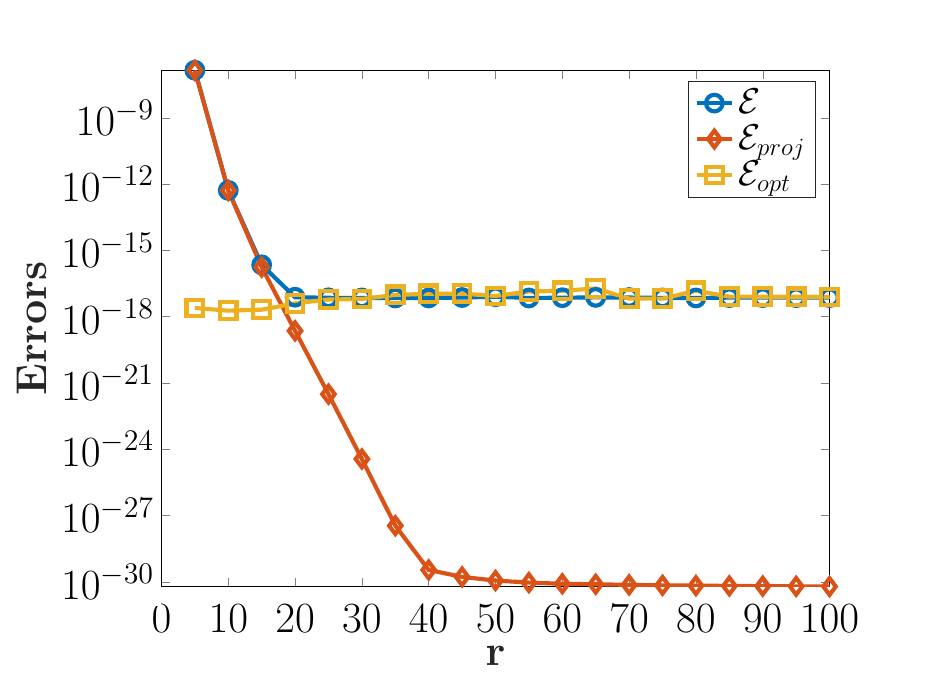}
\end{minipage}
\begin{minipage}[ht]{0.48\linewidth}
\includegraphics[width=1\textwidth]{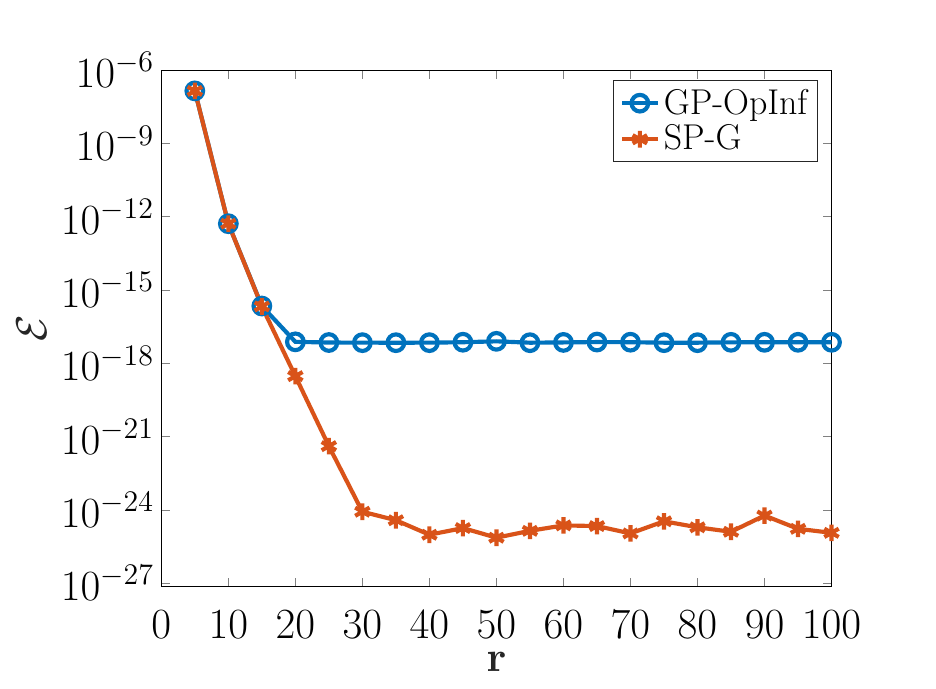}
\end{minipage}
\caption{(1D Allen-Cahn Equation) Numerical errors for the $r$-dimensional {\em GP-OpInf} when $T_{\text{FOM}}=T_{\text{ROM}} = 3$ and $\Delta t=2.5\times 10^{-4}$: (left) the ROM error together with the POD projection error \eqref{eq:e_proj} and the optimization error \eqref{eq:e_opt}; (right) comparison of the {\em GP-OpInf} ROM approximation error \eqref{eq:e_approx} with that of the {\em SP-G} ROM.}
\label{fig:1D_AC_error_estimation}
\end{figure}

\paragraph{Test 1. Illustration of error estimation}
We set $T_{\text{FOM}} = T_{\text{ROM}} = 3$, fix the parameter $\mu=1$, and choose a small time step $\Delta t = 2.5\times 10^{-4}$ so that the time discretization error is negligible in the error estimate. Varying the dimension $r$ of the ROM, we compute the {\em GP-OpInf} ROM approximation error \eqref{eq:e_approx}, projection error~\eqref{eq:e_proj} and optimization error~\eqref{eq:e_opt}, which are shown in Figure~\ref{fig:1D_AC_error_estimation} (left). 
In this case, for $r\leq 15$, the optimization error is smaller than the POD projection error; however, for a bigger $r$, the former becomes larger and levels off.  
Correspondingly, the {\em GP-OpInf} ROM approximation error first decreases monotonically and then reaches a plateau. 
Therefore, compared to the {\em SP-G} ROM, the {\em GP-OpInf} ROM achieves the same accuracy for small $r$, but yields larger numerical error when its dimension is large, as shown in Figure~\ref{fig:1D_AC_error_estimation} (right).

\begin{figure}[htb!]
\centering
\begin{minipage}[ht]{0.48\linewidth}
\includegraphics[width=1\textwidth]{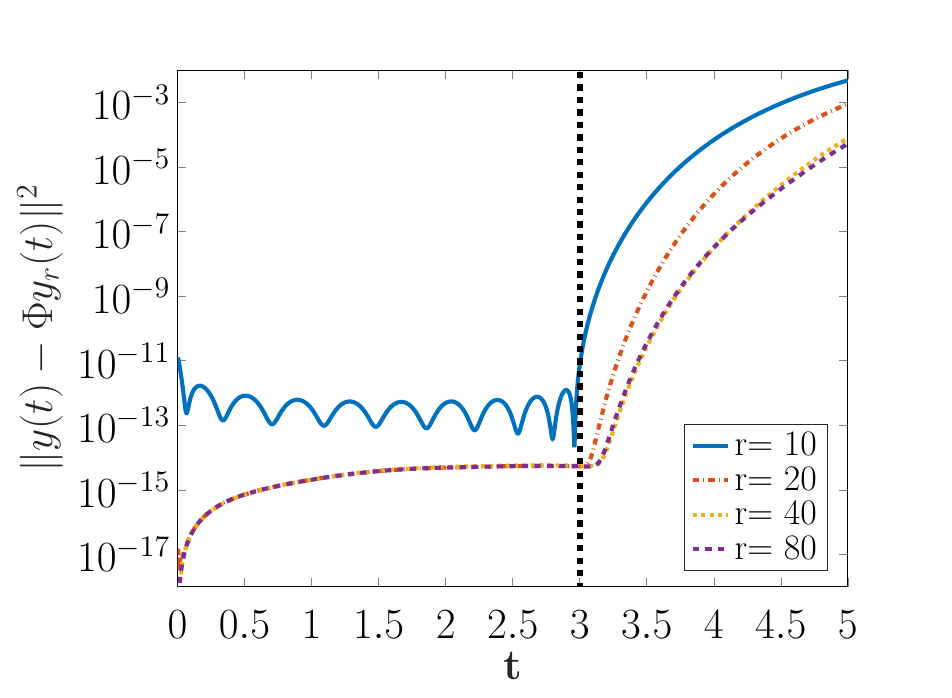}
\end{minipage}
\begin{minipage}[ht]{0.48\linewidth}
\includegraphics[width=1\textwidth]{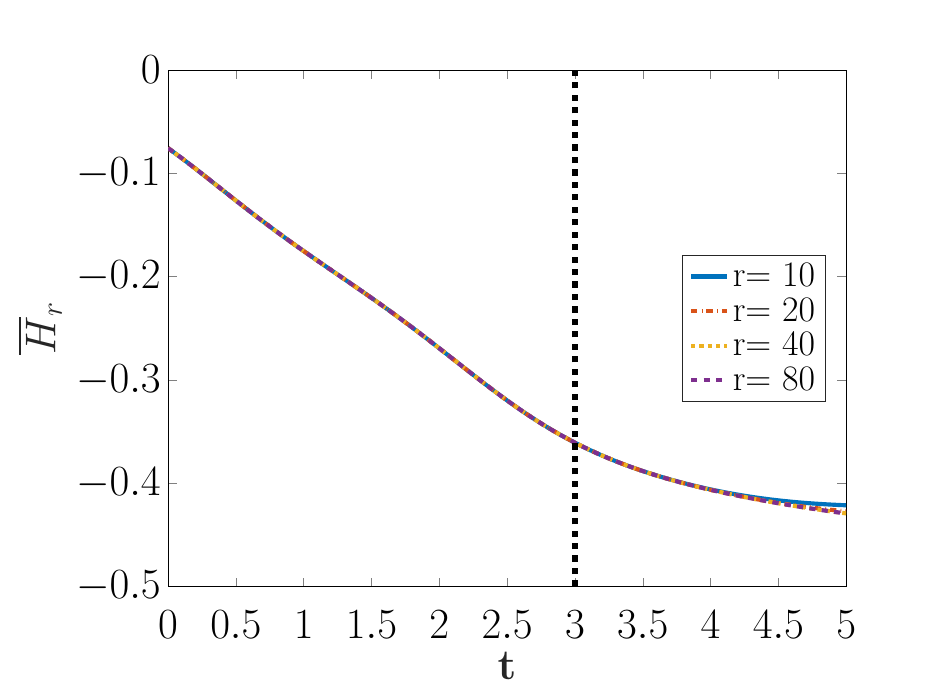}
\end{minipage}
\caption{(1D Allen-Cahn Equation) Numerical accuracy of the $r$-dimensional {\em GP-OpInf} ROM: (left) time evolution of the summand of the ROM approximation error~\eqref{eq:e_approx}; (right) time evolution of the reduced-order approximate Lyapunov function. The dashed line indicates the end of the training interval.}
\label{fig:1D_AC_long_simulation}
\end{figure}

\paragraph{Test 2. Prediction-in-time capabilities of the GP-OpInf ROM}
Setting $\Delta t = 10^{-3}$, we generate snapshots from the full-order simulation with the final time $T_{\text{FOM}}=3$ and simulate the ROM to the final time $T_{\text{ROM}}=5$, which is 66\% past the training data. The ROM approximation errors \eqref{eq:e_approx} at $r= 10, 20, 40$ and $80$ are plotted in Figure~\ref{fig:1D_AC_long_simulation} (left) alongside the FOM solutions. The associated approximate Lyapunov function is plotted in Figure~\ref{fig:1D_AC_long_simulation} (right). The dashed line indicates the end of the training time interval. It is evident that the ROM is accurate within the interval $[0, 3]$, but its accuracy degrades when the simulation time exceeds this range. Meanwhile, the approximate Lyapunov is decreasing in time as guaranteed by the gradient structure of the ROM.

\paragraph{Test 3. Parametric predictions away from training data} 

We parameterize the problem with the interface parameter $\mu\in {P}= [0.2,2]$, which also affects the initial condition. To generate snapshots, we simulated the FOM at 10 uniformly distributed training parameter samples in ${P}$ with the time step $\Delta t = 10^{-3}$ and final simulation time $T_{\text{FOM}}=3$. We construct the {\em GP-OpInf} ROM from the collected snapshot data and use it to obtain predictions at any given test parameter over the same time interval ($T_{\text{ROM}}=T_{\text{FOM}}$). Here, we select three random test parameters, $\mu = 0.34, 0.96$ and $1.87$.

\begin{figure}[t!]
\centering
\begin{minipage}[ht]{0.32\linewidth}
\includegraphics[width=1.1\textwidth]{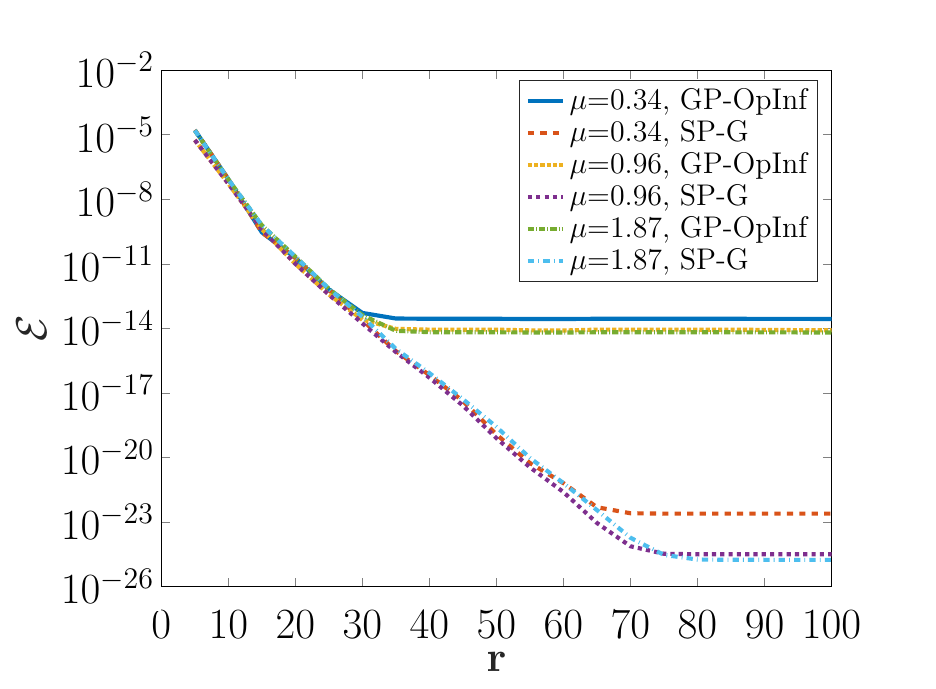}
\end{minipage}
\begin{minipage}[h]{0.32\linewidth}
	\includegraphics[width=1.1\textwidth]{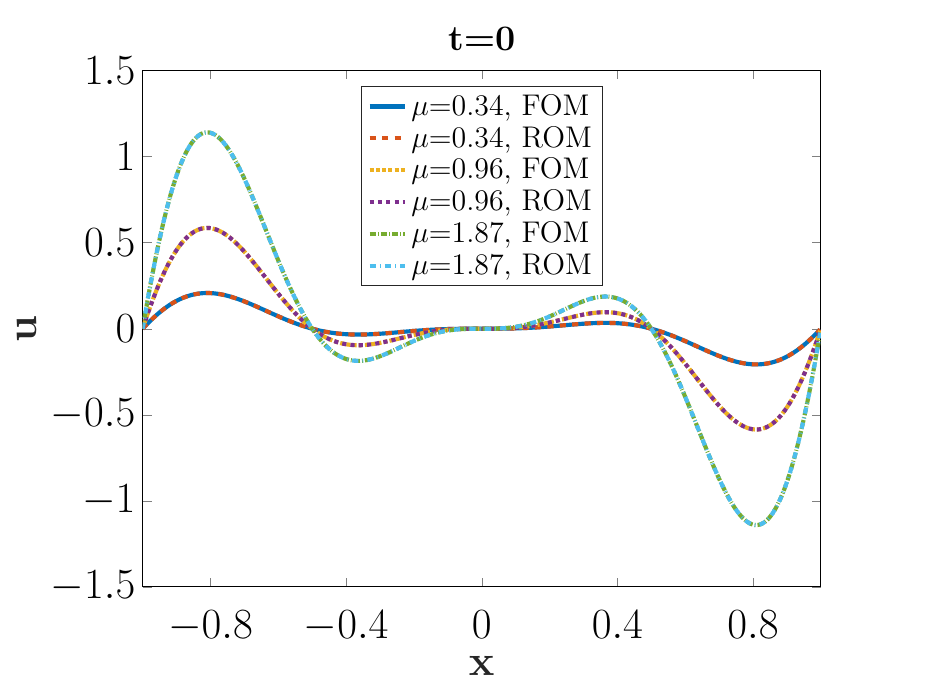}
\end{minipage}

\begin{minipage}[h]{0.32\linewidth}
	\includegraphics[width=1.1\textwidth]{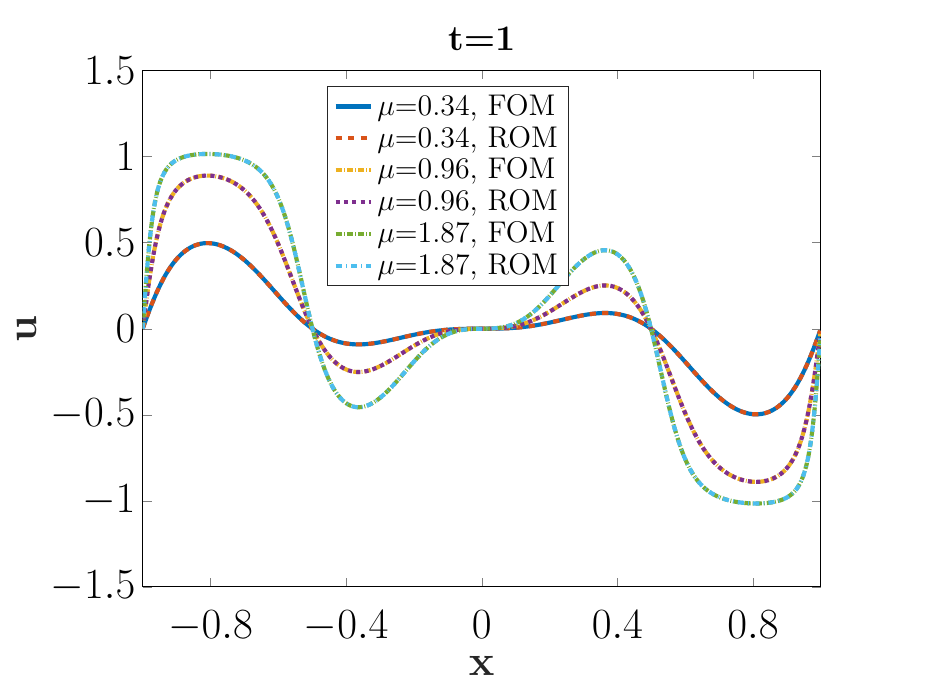}
\end{minipage}
	\begin{minipage}[h]{0.32\linewidth}
	\includegraphics[width=1.1\textwidth]{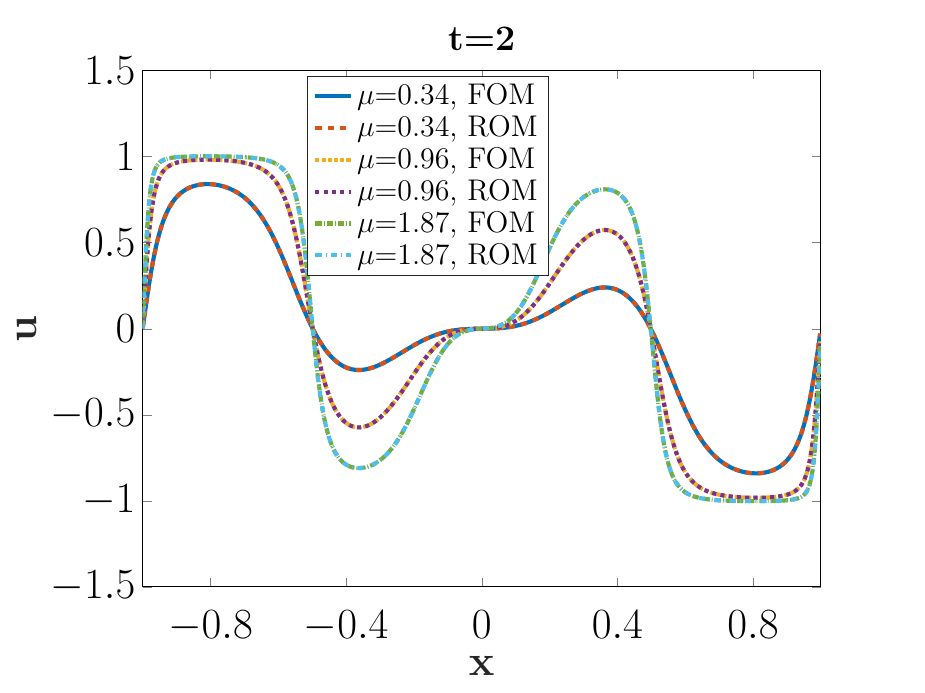}
\end{minipage}
\begin{minipage}[h]{0.32\linewidth}
	\includegraphics[width=1.1\textwidth]{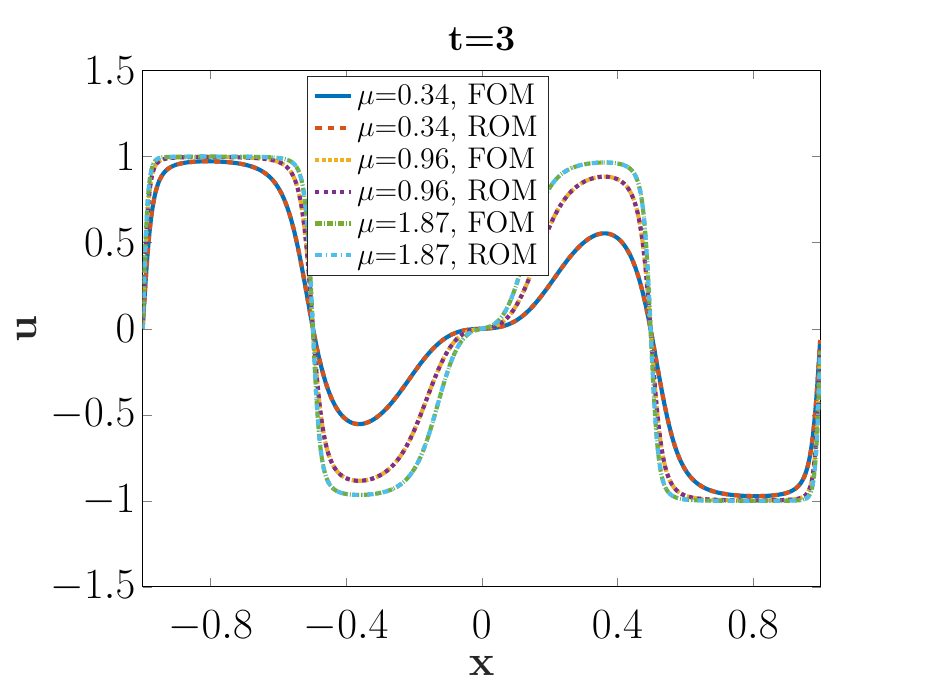}
\end{minipage}

\caption{(1D Allen-Cahn Equation) Performance of $r$-dimensional {\em GP-OpInf} ROM: (top left) comparison of the ROM approximation error~\eqref{eq:e_approx} for the {\em SP-G} ROM in problems with different $\mu$; (top right and bottom) snapshots of the ROM and FOM solutions at $t=0, 1, 2$ and $3$. 
}\label{fig:1D_AC_diff_mu}
\end{figure}

Figure~\ref{fig:1D_AC_diff_mu} (top left) shows the approximation errors of the {\em SP-G} and {\em GP-OpInf} ROMs. This illustrates that, at all test parameters, the error of {\em GP-OpInf} is close to that of {\em SP-G} when $r\leq 30$. The former levels off at values around $10^{-15}$ when $r$ becomes large, which is well beyond a typically required accuracy. Therefore, the {\em GP-OpInf} is very accurate.  
Figure~\ref{fig:1D_AC_diff_mu} (top right and bottom) compares the FOM solution with the {\em GP-OpInf} ROM solution of dimension $r=30$ at several time instances. This confirms that the {\em GP-OpInf} ROM provides an accurate approximation to the FOM solution in all cases. 

\subsubsection{Two-dimensional parameterized Allen-Cahn equation}
Consider the two-dimensional Allen-Cahn equation on a rectangular domain $\Omega$:

\begin{equation*}\label{eq1:2D_ac}
    u_t = \epsilon^2 \Delta u + u-u^3,\quad (x, y) \in \Omega, t \in I, 
\end{equation*}
where $u=u(x, t; \mu)$ is the parameter-dependent solution and $\mu\in P$ a parameter appearing in the initial condition, defined below alongside the boundary and initial conditions.
The Allen-Cahn equation can be recast into the form \eqref{eq:h-pde}, in a similar way as in the 1D case, 
with the Lyapunov function 
$\mathcal{H}(u) = \int_\Omega \left(\frac{\epsilon^2}{2} |\nabla u|^2 -\frac{1}{2}u^2  +\frac{1}{4}u^4 \right)\, {\rm d}x$. 
Let $\bu$ be a vectorization of the approximate state defined on the rectangular grid, obtained by applying a spatial discretization with a uniform mesh size $\Delta x$ and $\Delta y$ in the horizontal and vertical directions. The semi-discrete system is then
$$
\mathbf{\dt u}=-\mathbf{I}_{n^2}(-\epsilon^2 \mathbf{L}\mathbf{u} -\mathbf{u}+\mathbf{u}^3),
$$
where $\mathbf{L} = \mathbf{I}_n\otimes \bD_x^2+\bD_y^2 \otimes \mathbf{I}_n$ with $\bD_x^2= \frac{1}{(\Delta x)^2} \mathbf{L}_n$, 
$\bD_y^2= \frac{1}{(\Delta y)^2} \mathbf{L}_n$ and $\otimes$ denotes the Kronecker product \cite{cai2019structure}.

\paragraph{Computational setting} 

Let the problem be defined on $\Omega = [-0.5,0.5]^2$ and $I = [0,T_{\text{FOM}}]$. Consider a periodic boundary condition and $\epsilon = 0.02$. The initial condition (following \cite{geng2023deep}) is given by  
\begin{equation*}
    u_0(x,y; \mu)= \max \left(\tanh\left(\frac{r-\sqrt{(x-r(0.7-\mu))^2+(y+r\mu)^2}}{\epsilon}\right), 
    \tanh\left(\frac{r-\sqrt{(x+r(0.7-\mu))^2+(y-r\mu)^2}}{\epsilon}\right) \right)
    \label{eq:2D_AC_inital}
\end{equation*}
with $r=0.2$ and $\mu$ a parameter in ${P}=[0, 0.7]$. The initial condition has two disks centered at $((0.7-\mu)r,-\mu r)$ and $(-(0.7-\mu)r, \mu r)$ with radius $r$.
In all FOM simulations, the spatial domain is partitioned by a rectangular grid with the mesh size $\Delta x =\Delta y = 2^{-6}$ and the AVF method is used for time integration with a time step $\Delta t$. The nonlinear system is solved by the Picard iteration.

\begin{figure}[htb!]
    \centering
    \begin{minipage}[ht]{1\linewidth}
	\includegraphics[width=1\textwidth]{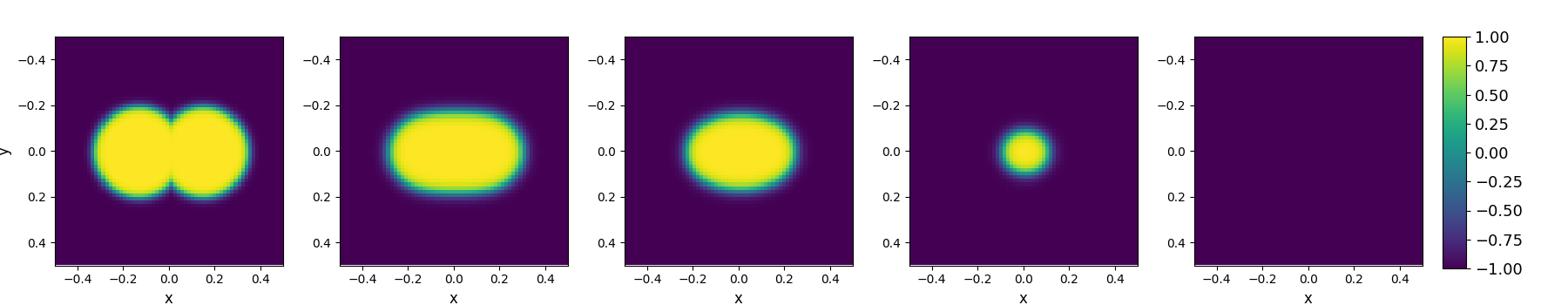}
    \end{minipage}
    \begin{minipage}[h]{0.48\linewidth}
    \includegraphics[width=1\textwidth]{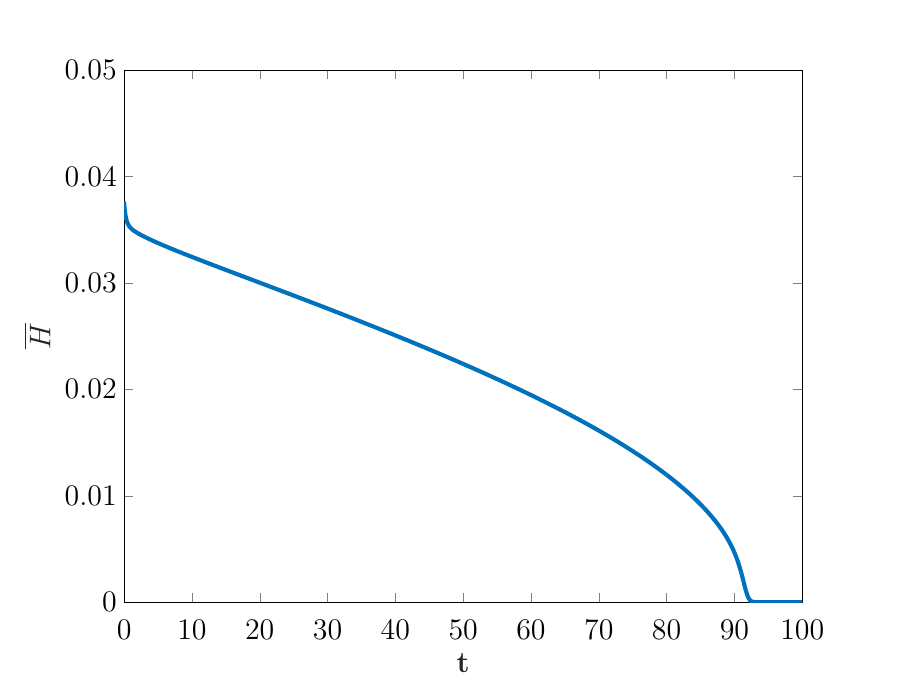}
    \end{minipage}
    \caption{(2D Allen-Cahn Equation) Full-order simulation when $T_{\text{FOM}}=100$ and $\Delta t = 10^{-2}$: (top) approximate states at $t= 0, 20, 40, 80, 100$; (bottom) time evolution of the discrete Lyapunov function.}
\label{fig:2D_AC_FOM_solution}
\end{figure}

Figure~\ref{fig:2D_AC_FOM_solution} shows the time evolution of the FOM state and the discrete Lyapunov function $\overline{H}(t)$ when $T_{\text{FOM}}=100$, $\Delta t=1\times10^{-2}$ and $\mu = 0$. 

Similar to the 1D Allen-Cahn equation discussed in Section~\ref{sec:1DAC}, we investigated the performance of the {\em GP-OpInf} ROM for 2D Allen-Cahn equation through three tests where we illustrate the error estimation, evaluate the {\em GP-OpInf}'s capability for time prediction, and examine parametric predictions beyond the training data. We observed that the numerical behaviors of the {\em GP-OpInf} ROM for the 2D Allen-Cahn equation closely align with those reported for the 1D case. Therefore, in the sequel, we only present the numerical test on the parametric problem in which $T_{\text{FOM}}=T_{\text{ROM}}=20$ and $\Delta t = 10^{-2}$. 

\begin{figure}[htb!]
\centering
\begin{minipage}[ht]{0.48\linewidth}
\includegraphics[width=1\textwidth]{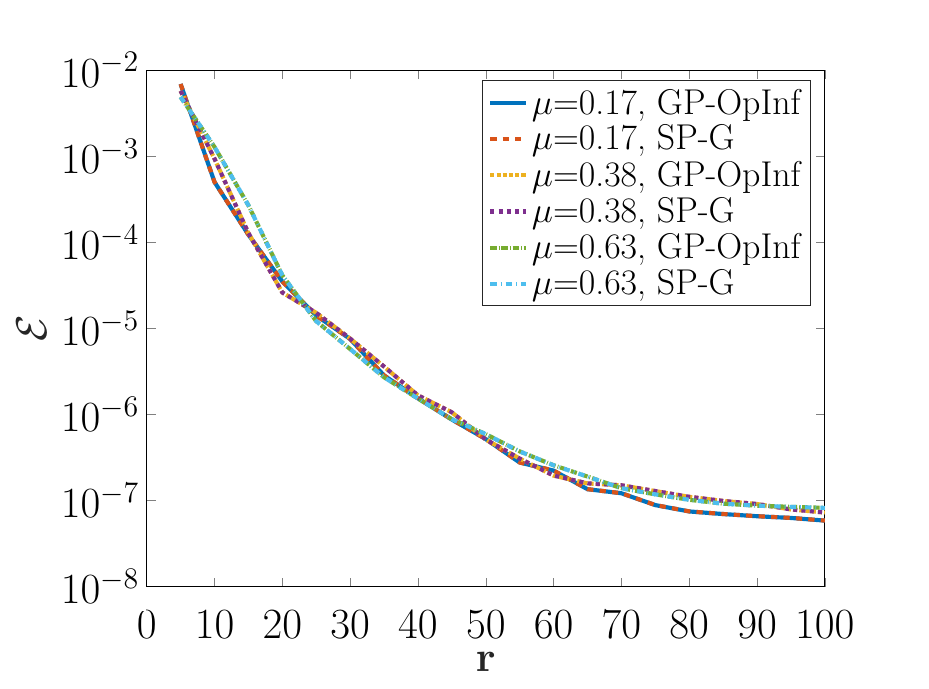}
\end{minipage}
\begin{minipage}[ht]{0.48\linewidth}
\includegraphics[width=1\textwidth]{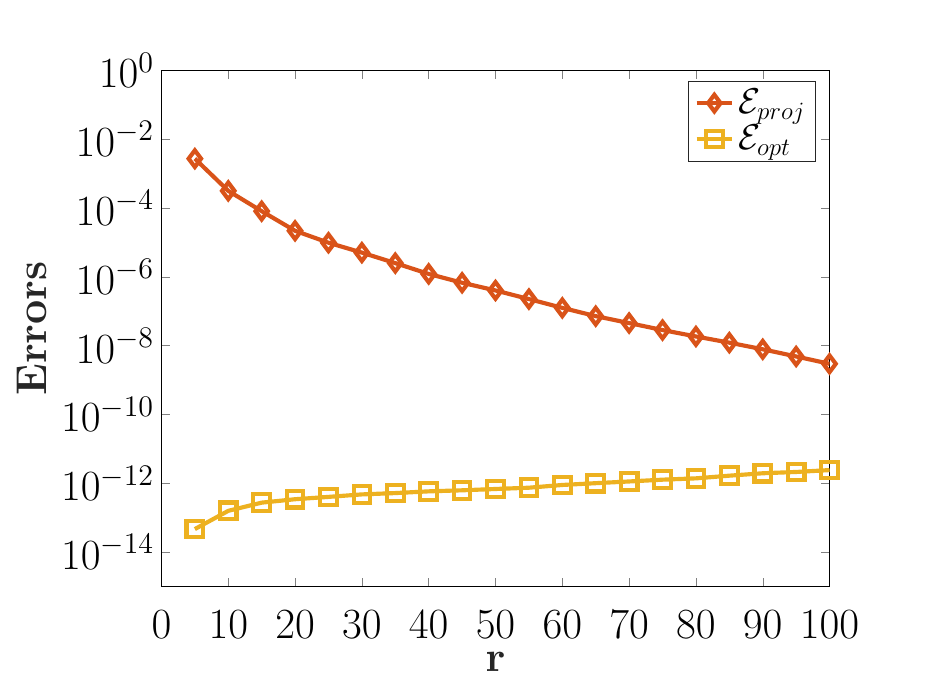}
\end{minipage}
\caption{
(2D Allen-Cahn Equation) Numerical errors for the $r$-dimensional {\em GP-OpInf} ROM: (left) comparison of the {\em GP-OpInf} ROM approximation error \eqref{eq:e_approx} with that of the {\em SP-G} ROM; (right) the projection error \eqref{eq:e_proj} and optimization error \eqref{eq:e_opt}.
}
\label{fig:2D_AC_diff_mu}
\end{figure}

\begin{figure}[htb!]
\centering
\includegraphics[width=\textwidth]{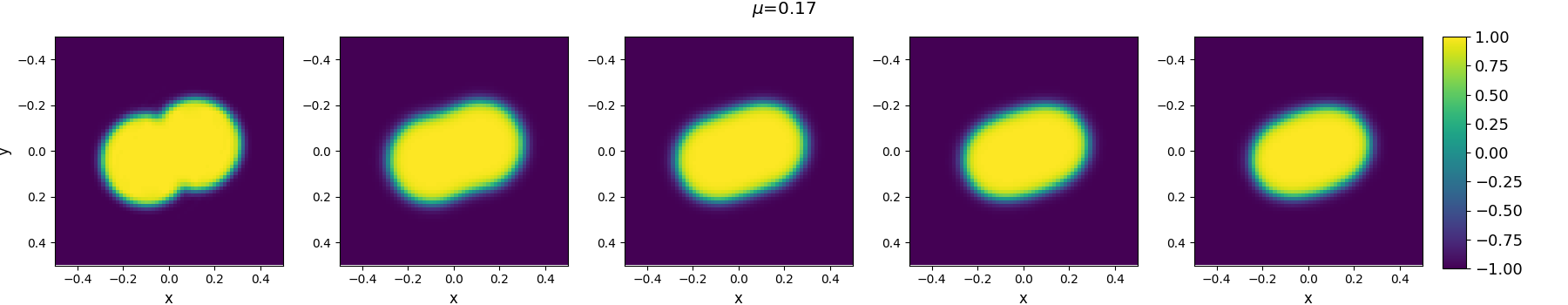}
\includegraphics[width=\textwidth]{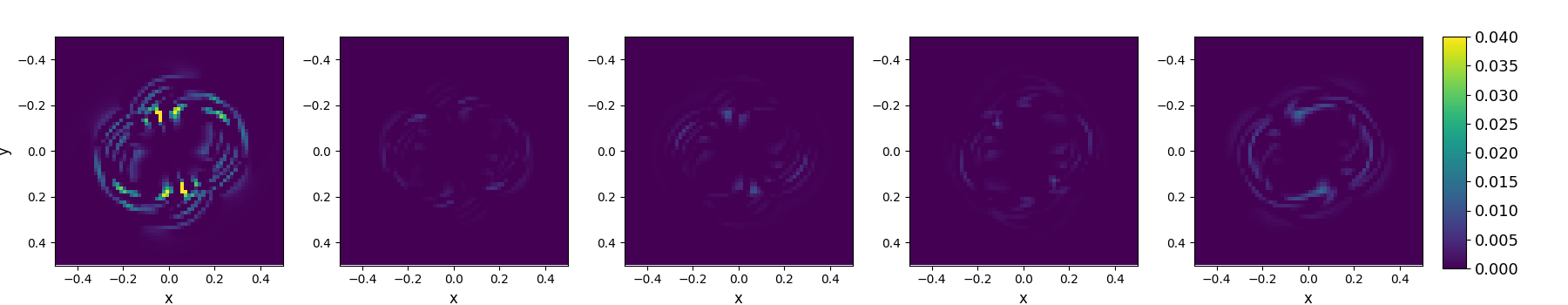}
\includegraphics[width=\textwidth]{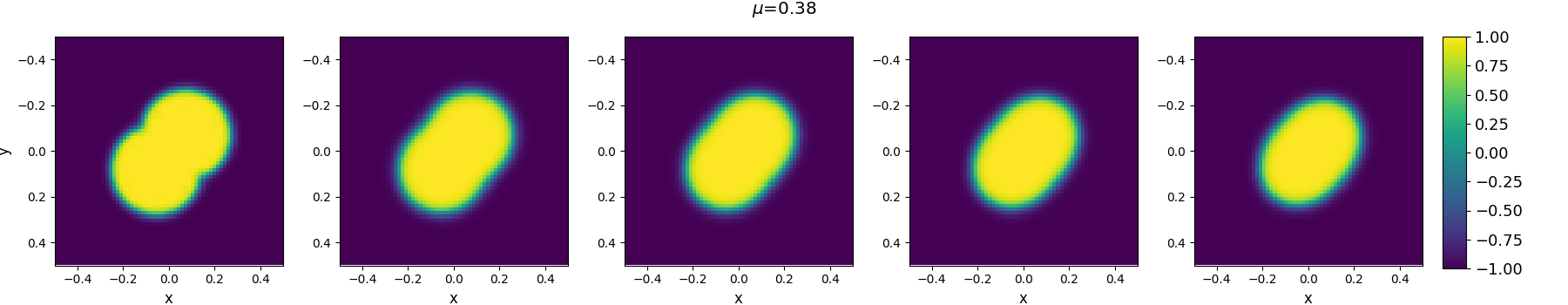}
\includegraphics[width=\textwidth]{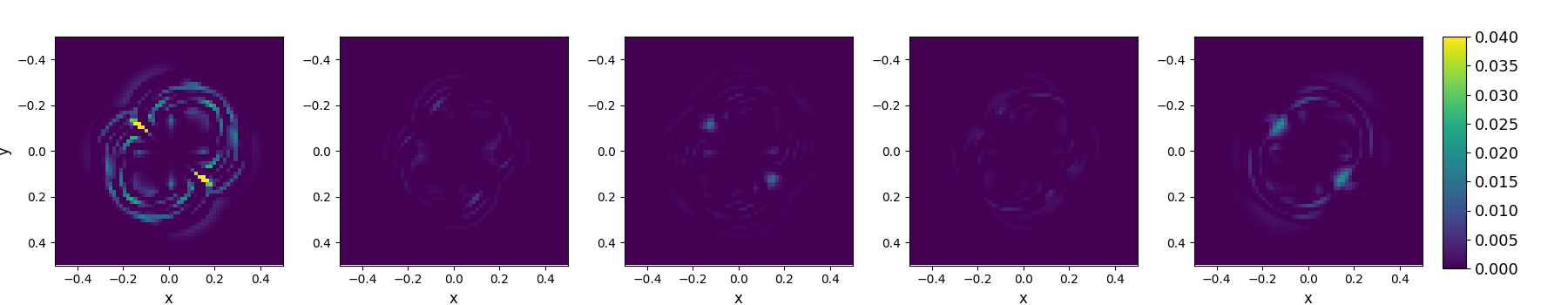}
\includegraphics[width=\textwidth]{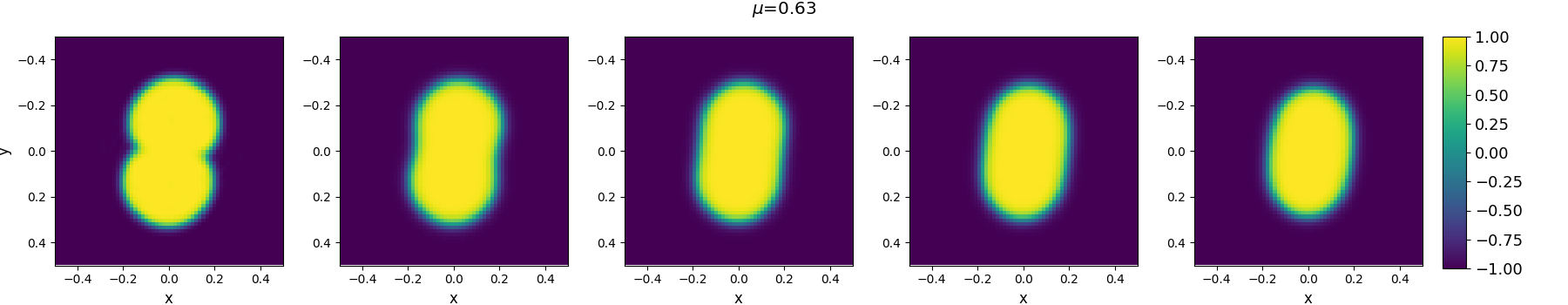}
\includegraphics[width=\textwidth]{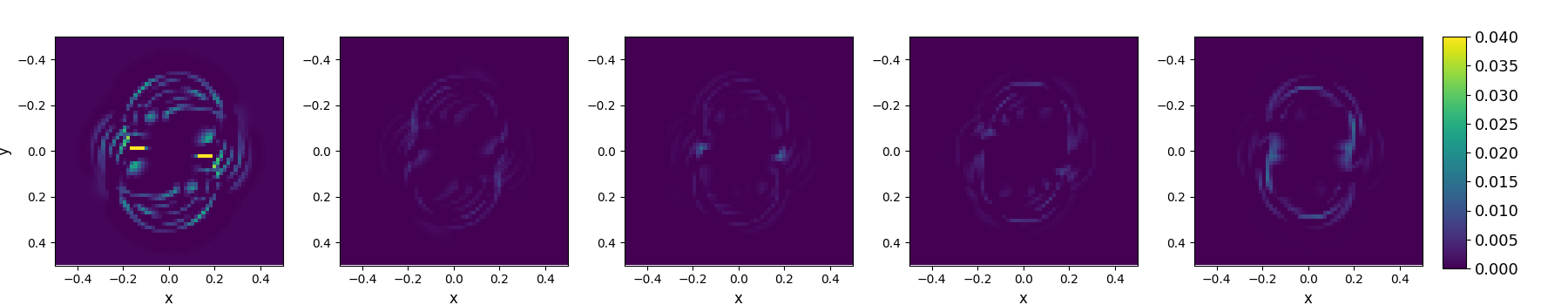}
\caption{(2D Allen-Cahn Equation) The {\em GP-OpInf} ROM of ($r= 40$) approximate solution (displayed in odd rows) and the associated approximation error $|\by(t,\mu) - \by_r(t,\mu)|$ (displayed in even rows) when $\mu = 0.17, 0.38$ and $0.63$, respectively, and $t=0, 5, 10, 15, 20$ (from left to right). 
}
\label{fig:2D_AC_diff_mu2}
\end{figure}

\paragraph{Test 1. Parametric predictions away from training data} 
To generate snapshots, 15 training samples are uniformly selected from ${P}$. The {\em GP-OpInf} ROM is then inferred from the data and used to predict solutions at any given testing samples from ${P}$. Here, the approximation errors of {\em SP-G} and {\em GP-OpInf} at three test parameters $\mu = 0.17, 0.38$ and $0.63$ are compared in Figure~\ref{fig:2D_AC_diff_mu} (left). The figure shows that the error of {\em GP-OpInf} is close to that of {\em SP-G} with the same dimension $r$ at all test parameters. 
As indicated in Figure~\ref{fig:2D_AC_diff_mu} (right), the optimization error is significantly lower than the POD projection error, thereby aligning the trends of errors for both {\em GP-OpInf} and {\em SP-G} with the POD projection error. 
Fixing $r= 40$, we plot the ROM solutions and associated absolute errors at time instances $t= 0, 5, 10, 15, 20$ (from top to bottom) in Figure ~\ref{fig:2D_AC_diff_mu2}. We observe that the {\em GP-OpInf} ROM produces accurate solutions for the parametric problem at the selected testing samples.

\newpage

\section{Conclusions \label{sec:col}}
In this work, we considered evolutionary PDEs with a gradient structure with the goal to infer reduced-order models with the same gradient structure from simulated data of the semi-discretized PDEs.
We first projected the high-dimensional snapshot data onto a POD basis, and from the projected data we infer a low-dimensional operator by solving a constrained optimization problem. The resulting operators then define a low-dimensional model, termed the gradient-preserving Operator Inference ({\em GP-OpInf}) ROM. 
By enforcing the proper constraints in the inference of the low-dimensional operators, we ensure that the {\em GP-OpInf} ROM has the appropriate gradient structure, thereby preserving the essential physics in the reduced-order dynamics. 
We further analyzed the associated approximation error of the {\em GP-OpInf} ROM. The analysis shows that an upper bound for the error consists of the combination of the POD projection error, the data error, and the OpInf optimization error. 
We test the accuracy, structure-preservation properties, and predictive capabilities of the {\em GP-OpInf} ROM on several PDE examples. For conservative test cases, we consider the parameterized wave equation and Korteweg-de-Vries equation, \Ra{as well as a three-dimensional linear elasticity problem}, and for dissipative test cases we consider the one-and two-dimensional Allen-Cahn equation. 
The numerical experiments show that {\em GP-OpInf} can achieve the same performance as the intrusive, projection-based ROM when the optimization error and data error are dominated by the POD projection error, which is typically true if the dimension of the ROM is low. This observation aligns with our error estimation. 
%
Moreover, for parametric problems, the low-dimensional {\em GP-OpInf} ROMs can attain accuracy close to that of the full-order simulations even for unseen parameters. 

\Ra{Going forward, it would be interesting to learn the Hamiltonian or its gradient function entirely from data, without relying on any prior knowledge. For this purpose, the work in \cite{ShNaToKr-Lagrangian-OPINF-SPML-2024} could be extended to the Hamiltonian case. Moreover, a joint learning of the dissipative and conservative matrices, similar to the port-Hamiltonian DMD framework from \cite{morandin2023port} would be an important next step.} 
\section*{CRediT authorship contribution statement}
\textbf{Yuwei Geng:} Software, Investigation, Visualization.
\textbf{Jasdeep Singh:} Software, Investigation.
\textbf{Lili Ju:} Methodology, Writing - Review \& Editing.
\textbf{Boris Kramer:} Validation, Writing - Review \& Editing, Funding Acquisition. 
\textbf{Zhu Wang:} Conceptualization, Methodology, Writing - Original Draft, Review \& Editing.

\section*{Declaration of competing interest}
Boris Kramer reports a relationship with ASML Holding US that includes: consulting or advisory.
The other authors declare that they have no known competing financial interests or personal relationships that could have appeared to influence the work reported in this paper.

\section*{Data availability}
Data will be made available upon request.

\section*{Acknowledgements}
\Ra{We would like to thank the anonymous reviewers for their valuable suggestions and feedback. We are also grateful to Dr.~Irena Tezaur (Sandia National Laboratories) for making their codes and data for the 3D plate simulations available online, and for patiently providing us with instructions on setting up the experiment.}

B.K. and Z.W. were in part supported by the U.S. Office of Naval Research under award number N00014-22-1-2624. Moreover, B.K. was in part  supported by the Ministry of Trade, Industry and Energy (MOTIE) and the Korea Institute for Advancement of Technology (KIAT) through the International Cooperative R\&D program (No. P0019804, Digital twin
based intelligent unmanned facility inspection solutions). 
Z. W. was in part supported by the U.S. National Science Foundataion under award number DMS-2038080 and an ASPIRE grant from the Office of the Vice President for Research at the University of South Carolina. 
L. J. was in part supported by the U.S. National Science Foundation  under award number DMS-2109633 and the U.S. Department of Energy  under award number  DE-SC0022254.

\clearpage


\appendix
\section{Iterative algorithm for solving \eqref{eq:sdp}}
\label{app1}

    To solve \eqref{eq:sdp} ,  we first find the gradient of the objective function. 
    The barrier function $g(\bX)$ can be written as follows. 
    \begin{align*}
        g(\bX) = \sum_{i=1}^r \log(-\lambda_i(\bX))
         = \log \left(\prod_{i=1}^{r}(-\lambda_i(\bX))\right)
         = \log((-1)^r det(\bX)). 
    \end{align*}
    Without loss of generality, we consider an even $r$ here. 
    Differentiating it with respect to the $ij$ entry yields 
    \begin{align*}
        \frac{\partial \log(det(\bX))}{\partial \bX_{ij}} = \frac{1}{det(\bX)}\frac{\partial det(\bX)}{\partial \bX_{ij}}.
    \end{align*}  
    Based on the cofactor formula for determinant, $det(\bX)=\sum_{k=1}^{s} \bX_{ik} \bC_{ik}$, we have 
    \begin{align*}
        \frac{\partial det(\bX)}{\partial \bX_{ij}} &= \frac{\partial (\sum_{k=1}^{r} \bX_{ik} \bC_{ik})}{\partial \bX_{ij}}
         = \sum_{k=1}^s \left(\frac{\partial \bX_{ik}}{\bX_{ij}}\bC_{ik} +\bX_{ik} \frac{\partial \bC_{ik}}{\partial \bX_{ij}}\right)
         = \bC_{ij}.
    \end{align*}
    Therefore, we have 
    $$\frac{\partial \log(det(\bX))}{\partial \bX_{ij} }=  \frac{1}{det(\bX)} \bC_{ij}.$$ 
    Due to the fact  that $\bX^{-1} = \frac{1}{det(\bX)} \bC^\intercal$, we have $\nabla g(\bX) = \bX^{-\intercal}$, and 
    \begin{equation}\label{eq:gradg}
        \nabla g(\bD_r) = \bD_r^{-\intercal}.
    \end{equation}
   Next, $f(\bD_r)$ is recast into the following form using the trace: 
    \begin{align*}
        f(\bD_r) &= \text{tr} \left((\dt\bY_r - \bD_r\bF_r)(\dt\bY_r - \bD_r\bF_r)^\intercal \right) \\
                & = \text{tr} \left(\dt\bY_r \dt\bY_r^\intercal - \dt\bY_r\bF_r^\intercal  \bD_r^\intercal - \bD_r\bF_r \dt\bY_r^\intercal  +\bD_r \bF_r \bF_r^\intercal  \bD_r^\intercal \right).
    \end{align*}
Then the gradient of $f(\bD_r)$ can be obtained as follows:  
    \begin{align}
    \label{eq:gradf}
        \nabla f(\bD_r)  
        = -\dt\bY_r\bF_r^\intercal - \dt\bY_r\bF_r^\intercal + (\bF_r \bF_r^\intercal \bD_r^\intercal)^\intercal +\bD_r \bF_r \bF_r^\intercal
        = -2(\dt\bY_r - \bD_r\bF_r)\bF_r^\intercal. 
    \end{align} 
By combining \eqref{eq:gradf} and \eqref{eq:gradg} , we have 
    \begin{align*}
        \nabla f_{\beta}(\bD_r) = \nabla f(\bD_r) - \beta \nabla g(\bD_r)
            = -2(\dt\bY_r - \bD_r\bF_r)\bF_r^\intercal - \beta \bD_r^{-\intercal}.
    \end{align*}  
Secondly, we  use the gradient descent method \cite{amari1993backpropagation} to approximate the solution $\bD_r^{k+1}$ by the given $\bD_r^{k}$, which follows
\begin{equation*}
     \bD_r^{k+1}  = \bD_r^{k} - \tau^k \nabla f_{\beta}(\bD_r^{k})
\end{equation*}
where $\tau^k \geq 0 $ is the step size. 
        Instead of a constant barrier parameter $\beta$, one could use a decaying parameter with $\beta^{k+1} = \sigma \beta^{k}$, where $\sigma \in (0,1)$ and its value depends on $r$ \cite{potra2000interior,gonzaga1989algorithm}. 
        
This iteration continues till the magnitude of gradient $\nabla f_{\beta}(\bD_r)$ is sufficiently small and the barrier parameter $\beta$ drops below a user-defined tolerance.

\bibliographystyle{abbrv}
\bibliography{reference.bib}

\begin{thebibliography}{10}

\bibitem{afkham2017structure}
B.~M. Afkham and J.~S. Hesthaven.
\newblock Structure preserving model reduction of parametric {H}amiltonian
  systems.
\newblock {\em SIAM Journal on Scientific Computing}, 39(6):A2616--A2644, 2017.

\bibitem{agrawal2018rewriting}
A.~Agrawal, R.~Verschueren, S.~Diamond, and S.~Boyd.
\newblock A rewriting system for convex optimization problems.
\newblock {\em Journal of Control and Decision}, 5(1):42--60, 2018.

\bibitem{amari1993backpropagation}
S.-i. Amari.
\newblock Backpropagation and stochastic gradient descent method.
\newblock {\em Neurocomputing}, 5(4-5):185--196, 1993.

\bibitem{antoulas2005approximation}
A.~C. Antoulas.
\newblock {\em Approximation of large-scale dynamical systems}.
\newblock SIAM, 2005.

\bibitem{antoulas2020interpolatory}
A.~C. Antoulas, C.~A. Beattie, and S.~G{\"u}{\u{g}}ercin.
\newblock {\em Interpolatory methods for model reduction}.
\newblock SIAM, 2020.

\bibitem{antsaklis1997linear}
P.~J. Antsaklis and A.~N. Michel.
\newblock {\em Linear systems}, volume~8.
\newblock Springer, 1997.

\bibitem{mosek}
M.~ApS.
\newblock {\em The MOSEK optimization toolbox for MATLAB manual. Version 9.0.},
  2019.

\bibitem{arakawa_1966_incompressible}
A.~Arakawa.
\newblock {Computational design for long-term numerical integration of the
  equations of fluid motion: Two-dimensional incompressible flow. Part I}.
\newblock {\em Journal of Computational Physics}, 1(1):119--143, 1966.

\bibitem{bartels1972algorithm}
R.~H. Bartels and G.~W. Stewart.
\newblock Algorithm 432 [c2]: solution of the matrix equation ax+ xb= c [f4].
\newblock {\em Communications of the ACM}, 15(9):820--826, 1972.

\bibitem{beattie2009interpolatory}
C.~Beattie and S.~Gugercin.
\newblock Interpolatory projection methods for structure-preserving model
  reduction.
\newblock {\em Systems \& Control Letters}, 58(3):225--232, 2009.

\bibitem{beattie2011structure}
C.~Beattie and S.~Gugercin.
\newblock Structure-preserving model reduction for nonlinear port-{H}amiltonian
  systems.
\newblock In {\em 2011 50th IEEE conference on decision and control and
  European control conference}, pages 6564--6569. IEEE, 2011.

\bibitem{benner2015survey}
P.~Benner, S.~Gugercin, and K.~Willcox.
\newblock A survey of projection-based model reduction methods for parametric
  dynamical systems.
\newblock {\em SIAM review}, 57(4):483--531, 2015.

\bibitem{berkooz1993proper}
G.~Berkooz, P.~Holmes, and J.~L. Lumley.
\newblock The proper orthogonal decomposition in the analysis of turbulent
  flows.
\newblock {\em Annual Review of Fluid Mechanics}, 25(1):539--575, 1993.

\bibitem{brunton2022data}
S.~L. Brunton and J.~N. Kutz.
\newblock {\em Data-driven science and engineering: Machine learning, dynamical
  systems, and control}.
\newblock Cambridge University Press, 2022.

\bibitem{buchfink2023symplectic}
P.~Buchfink, S.~Glas, and B.~Haasdonk.
\newblock Symplectic model reduction of {H}amiltonian systems on nonlinear
  manifolds and approximation with weakly symplectic autoencoder.
\newblock {\em SIAM Journal on Scientific Computing}, 45(2):A289--A311, 2023.

\bibitem{cai2019structure}
W.~Cai, C.~Jiang, Y.~Wang, and Y.~Song.
\newblock Structure-preserving algorithms for the two-dimensional {sine-Gordon}
  equation with {N}eumann boundary conditions.
\newblock {\em Journal of Computational Physics}, 395:166--185, 2019.

\bibitem{carlberg2015preserving}
K.~Carlberg, R.~Tuminaro, and P.~Boggs.
\newblock Preserving {L}agrangian structure in nonlinear model reduction with
  application to structural dynamics.
\newblock {\em SIAM Journal on Scientific Computing}, 37(2):B153--B184, 2015.

\bibitem{celledoni2012preserving}
E.~Celledoni, V.~Grimm, R.~I. McLachlan, D.~McLaren, D.~O’Neale, B.~Owren,
  and G.~Quispel.
\newblock Preserving energy resp. dissipation in numerical {PDE}s using the
  “average vector field” method.
\newblock {\em Journal of Computational Physics}, 231(20):6770--6789, 2012.

\bibitem{chaturantabut2016structure}
S.~Chaturantabut, C.~Beattie, and S.~Gugercin.
\newblock Structure-preserving model reduction for nonlinear port-{H}amiltonian
  systems.
\newblock {\em SIAM Journal on Scientific Computing}, 38(5):B837--B865, 2016.

\bibitem{chaturantabut2012state}
S.~Chaturantabut and D.~C. Sorensen.
\newblock A state space error estimate for {POD-DEIM} nonlinear model
  reduction.
\newblock {\em SIAM Journal on numerical analysis}, 50(1):46--63, 2012.

\bibitem{chen2021conservative}
Q.~Chen, L.~Ju, and R.~Temam.
\newblock Conservative numerical schemes with optimal dispersive wave
  relations: Part i. derivation and analysis.
\newblock {\em Numerische Mathematik}, 149(1):43--85, 2021.

\bibitem{christiansen2011topics}
S.~H. Christiansen, H.~Z. Munthe-Kaas, and B.~Owren.
\newblock Topics in structure-preserving discretization.
\newblock {\em Acta Numerica}, 20:1--119, 2011.

\bibitem{diamond2016cvxpy}
S.~Diamond and S.~Boyd.
\newblock {CVXPY}: {A} {P}ython-embedded modeling language for convex
  optimization.
\newblock {\em Journal of Machine Learning Research}, 17(83):1--5, 2016.

\bibitem{evans2022partial}
L.~C. Evans.
\newblock {\em Partial differential equations}, volume~19.
\newblock American Mathematical Society, 2022.

\bibitem{filanova2023operator}
Y.~Filanova, I.~P. Duff, P.~Goyal, and P.~Benner.
\newblock An operator inference oriented approach for linear mechanical
  systems.
\newblock {\em Mechanical Systems and Signal Processing}, 200:110620, 2023.

\bibitem{fukao2017structure}
T.~Fukao, S.~Yoshikawa, and S.~Wada.
\newblock Structure-preserving finite difference schemes for the
  {C}ahn-{H}illiard equation with dynamic boundary conditions in the
  one-dimensional case.
\newblock {\em Commun. Pure Appl. Anal}, 16(5):1915--1938, 2017.

\bibitem{furihata2010discrete}
D.~Furihata and T.~Matsuo.
\newblock {\em Discrete variational derivative method: a structure-preserving
  numerical method for partial differential equations}.
\newblock CRC Press, 2010.

\bibitem{geng2023deep}
Y.~Geng, Y.~Teng, Z.~Wang, and L.~Ju.
\newblock A deep learning method for the dynamics of classic and conservative
  {A}llen-{C}ahn equations based on fully-discrete operators.
\newblock {\em Journal of Computational Physics}, page 112589, 2023.

\bibitem{gong2017structure}
Y.~Gong, Q.~Wang, and Z.~Wang.
\newblock Structure-preserving {G}alerkin {POD} reduced-order modeling of
  {H}amiltonian systems.
\newblock {\em Computer Methods in Applied Mechanics and Engineering},
  315:780--798, 2017.

\bibitem{gonzaga1989algorithm}
C.~C. Gonzaga.
\newblock An algorithm for solving linear programming problems in o (n 3 l)
  operations.
\newblock In {\em Progress in Mathematical Programming: Interior-Point and
  Related Methods}, pages 1--28. Springer, 1989.

\bibitem{gross2005conservation}
M.~Gro{\ss}, P.~Betsch, and P.~Steinmann.
\newblock Conservation properties of a time {FE} method. {Part IV}: Higher
  order energy and momentum conserving schemes.
\newblock {\em International Journal for Numerical Methods in Engineering},
  63(13):1849--1897, 2005.

\bibitem{Gruber2023canonical}
A.~Gruber and I.~Tezaur.
\newblock Canonical and noncanonical {H}amiltonian operator inference.
\newblock {\em Computer Methods in Applied Mechanics and Engineering},
  416:116334, 2023.

\bibitem{gruber2024variationally}
A.~Gruber and I.~Tezaur.
\newblock Variationally consistent hamiltonian model reduction.
\newblock {\em arXiv preprint arXiv:2404.15315}, 2024.

\bibitem{gugercin2012structure}
S.~Gugercin, R.~V. Polyuga, C.~Beattie, and A.~Van Der~Schaft.
\newblock Structure-preserving tangential interpolation for model reduction of
  port-{H}amiltonian systems.
\newblock {\em Automatica}, 48(9):1963--1974, 2012.

\bibitem{gunzburger2007reduced}
M.~D. Gunzburger, J.~S. Peterson, and J.~N. Shadid.
\newblock Reduced-order modeling of time-dependent pdes with multiple
  parameters in the boundary data.
\newblock {\em Computer methods in applied mechanics and engineering},
  196(4-6):1030--1047, 2007.

\bibitem{hairer2006geometric}
E.~Hairer, M.~Hochbruck, A.~Iserles, and C.~Lubich.
\newblock Geometric numerical integration.
\newblock {\em Oberwolfach Reports}, 3(1):805--882, 2006.

\bibitem{helmberg1996interior}
C.~Helmberg, F.~Rendl, R.~J. Vanderbei, and H.~Wolkowicz.
\newblock An interior-point method for semidefinite programming.
\newblock {\em SIAM Journal on optimization}, 6(2):342--361, 1996.

\bibitem{hesthaven2021structure}
J.~Hesthaven and C.~Pagliantini.
\newblock Structure-preserving reduced basis methods for poisson systems.
\newblock {\em Mathematics of Computation}, 90(330):1701--1740, 2021.

\bibitem{hesthaven2022reduced}
J.~S. Hesthaven, C.~Pagliantini, and G.~Rozza.
\newblock Reduced basis methods for time-dependent problems.
\newblock {\em Acta Numerica}, 31:265--345, 2022.

\bibitem{kaagstrom1996lapack}
B.~K{\aa}gstr{\"o}m and P.~Poromaa.
\newblock Lapack-style algorithms and software for solving the generalized
  sylvester equation and estimating the separation between regular matrix
  pairs.
\newblock {\em ACM Transactions on Mathematical Software (TOMS)},
  22(1):78--103, 1996.

\bibitem{karasozen2018energy}
B.~Karas{\"o}zen and M.~Uzunca.
\newblock Energy preserving model order reduction of the nonlinear
  schr{\"o}dinger equation.
\newblock {\em Advances in Computational Mathematics}, 44:1769--1796, 2018.

\bibitem{KPW_OPINF_survey2024}
B.~Kramer, B.~Peherstorfer, and K.~E. Willcox.
\newblock Learning nonlinear reduced models from data with operator inference.
\newblock {\em Annual Review of Fluid Mechanics}, 56(1):521--548, 2024.

\bibitem{kunisch2001galerkin}
K.~Kunisch and S.~Volkwein.
\newblock Galerkin proper orthogonal decomposition methods for parabolic
  problems.
\newblock {\em Numerische mathematik}, 90:117--148, 2001.

\bibitem{lall2003structure}
S.~Lall, P.~Krysl, and J.~E. Marsden.
\newblock Structure-preserving model reduction for mechanical systems.
\newblock {\em Physica D: Nonlinear Phenomena}, 184(1-4):304--318, 2003.

\bibitem{marsden2013introduction}
J.~E. Marsden and T.~S. Ratiu.
\newblock {\em Introduction to Mechanics and Symmetry: A Basic Exposition of
  Classical Mechanical Systems}, volume~17.
\newblock Springer Science \& Business Media, 2013.

\bibitem{mclachlan1998six}
R.~I. McLachlan, G.~Quispel, et~al.
\newblock {\em Six lectures on the geometric integration of ODEs}.
\newblock Citeseer, 1998.

\bibitem{mclachlan2006geometric}
R.~I. McLachlan and G.~R.~W. Quispel.
\newblock Geometric integrators for {ODEs}.
\newblock {\em Journal of Physics A: Mathematical and General}, 39(19):5251,
  2006.

\bibitem{mclachlan1999geometric}
R.~I. McLachlan, G.~R.~W. Quispel, and N.~Robidoux.
\newblock Geometric integration using discrete gradients.
\newblock {\em Philosophical Transactions of the Royal Society of London.
  Series A: Mathematical, Physical and Engineering Sciences},
  357(1754):1021--1045, 1999.

\bibitem{mcquarrie2023nonintrusive}
S.~A. McQuarrie, P.~Khodabakhshi, and K.~E. Willcox.
\newblock Nonintrusive reduced-order models for parametric partial differential
  equations via data-driven operator inference.
\newblock {\em SIAM Journal on Scientific Computing}, 45(4):A1917--A1946, 2023.

\bibitem{miyatake2019structure}
Y.~Miyatake.
\newblock Structure-preserving model reduction for dynamical systems with a
  first integral.
\newblock {\em Japan Journal of Industrial and Applied Mathematics},
  36:1021--1037, 2019.

\bibitem{morandin2023port}
R.~Morandin, J.~Nicodemus, and B.~Unger.
\newblock Port-{H}amiltonian dynamic mode decomposition.
\newblock {\em SIAM Journal on Scientific Computing}, 45(4):A1690--A1710, 2023.

\bibitem{pagliantini2021dynamical}
C.~Pagliantini.
\newblock Dynamical reduced basis methods for {H}amiltonian systems.
\newblock {\em Numerische Mathematik}, 148(2):409--448, 2021.

\bibitem{pagliantini_2023_time_integration}
C.~Pagliantini, G.~Manzini, O.~Koshkarov, G.~L. Delzanno, and V.~Roytershteyn.
\newblock {Energy-conserving explicit and implicit time integration methods for
  the multi-dimensional {H}ermite-{DG} discretization of the {V}lasov-{M}axwell
  equations}.
\newblock {\em Computer Physics Communications}, 284:108604, 2023.

\bibitem{pagliantini2023gradient}
C.~Pagliantini and F.~Vismara.
\newblock Gradient-preserving hyper-reduction of nonlinear dynamical systems
  via discrete empirical interpolation.
\newblock {\em SIAM Journal on Scientific Computing}, 45(5):A2725--A2754, 2023.

\bibitem{peherstorfer2016data}
B.~Peherstorfer and K.~Willcox.
\newblock Data-driven operator inference for nonintrusive projection-based
  model reduction.
\newblock {\em Computer Methods in Applied Mechanics and Engineering},
  306:196--215, 2016.

\bibitem{peng2016symplectic}
L.~Peng and K.~Mohseni.
\newblock Symplectic model reduction of {H}amiltonian systems.
\newblock {\em SIAM Journal on Scientific Computing}, 38(1):A1--A27, 2016.

\bibitem{potra2000interior}
F.~A. Potra and S.~J. Wright.
\newblock Interior-point methods.
\newblock {\em Journal of computational and applied mathematics},
  124(1-2):281--302, 2000.

\bibitem{quarteroni2015reduced}
A.~Quarteroni, A.~Manzoni, and F.~Negri.
\newblock {\em Reduced basis methods for partial differential equations: an
  introduction}, volume~92.
\newblock Springer, 2015.

\bibitem{quispel2008new}
G.~Quispel and D.~I. McLaren.
\newblock A new class of energy-preserving numerical integration methods.
\newblock {\em Journal of Physics A: Mathematical and Theoretical},
  41(4):045206, 2008.

\bibitem{quispel1996discrete}
G.~Quispel and G.~S. Turner.
\newblock Discrete gradient methods for solving {ODE}s numerically while
  preserving a first integral.
\newblock {\em Journal of Physics A: Mathematical and General}, 29(13):L341,
  1996.

\bibitem{renegar2001mathematical}
J.~Renegar.
\newblock {\em A mathematical view of interior-point methods in convex
  optimization}.
\newblock SIAM, 2001.

\bibitem{salmon1982shape}
R.~Salmon.
\newblock The shape of the main thermocline.
\newblock {\em Journal of Physical Oceanography}, 12(12):1458--1479, 1982.

\bibitem{sanchez2021symplectic}
M.~A. S{\'a}nchez, B.~Cockburn, N.-C. Nguyen, and J.~Peraire.
\newblock Symplectic {H}amiltonian finite element methods for linear
  elastodynamics.
\newblock {\em Computer Methods in Applied Mechanics and Engineering},
  381:113843, 2021.

\bibitem{SK_LagrangianOPINF}
H.~Sharma and B.~Kramer.
\newblock Preserving {L}agrangian structure in data-driven reduced-order
  modeling of large-scale mechanical systems.
\newblock {\em Physica D: Nonlinear Phenomena}, 462:134128, 2024.

\bibitem{sharma2023symplectic}
H.~Sharma, H.~Mu, P.~Buchfink, R.~Geelen, S.~Glas, and B.~Kramer.
\newblock Symplectic model reduction of {H}amiltonian systems using data-driven
  quadratic manifolds.
\newblock {\em Computer Methods in Applied Mechanics and Engineering},
  417:116402, 2023.

\bibitem{sharma2020review}
H.~Sharma, M.~Patil, and C.~Woolsey.
\newblock A review of structure-preserving numerical methods for engineering
  applications.
\newblock {\em Computer Methods in Applied Mechanics and Engineering},
  366:113067, 2020.

\bibitem{sharma2022hamiltonian}
H.~Sharma, Z.~Wang, and B.~Kramer.
\newblock Hamiltonian operator inference: Physics-preserving learning of
  reduced-order models for canonical {H}amiltonian systems.
\newblock {\em Physica D: Nonlinear Phenomena}, 431:133122, 2022.

\bibitem{ShNaToKr-Lagrangian-OPINF-SPML-2024}
K.~Sharma, D.~A. Najera-Flores, M.~D. Todd, and B.~Kramer.
\newblock Lagrangian operator inference enhanced with structure-preserving
  machine learning for nonintrusive model reduction of mechanical systems.
\newblock {\em Computer Methods in Applied Mechanics and Engineering},
  423:116865, 2024.

\bibitem{singler2014new}
J.~R. Singler.
\newblock New {POD} error expressions, error bounds, and asymptotic results for
  reduced order models of parabolic {PDE}s.
\newblock {\em SIAM Journal on Numerical Analysis}, 52(2):852--876, 2014.

\bibitem{sirovich1987turbulence}
L.~Sirovich.
\newblock Turbulence and the dynamics of coherent structures. i. {C}oherent
  structures.
\newblock {\em Quarterly of applied mathematics}, 45(3):561--571, 1987.

\bibitem{stewart2010multilayer}
A.~L. Stewart and P.~J. Dellar.
\newblock Multilayer shallow water equations with complete {C}oriolis force.
  part 1. derivation on a non-traditional beta-plane.
\newblock {\em Journal of Fluid Mechanics}, 651:387--413, 2010.

\bibitem{stewart2016energy}
A.~L. Stewart and P.~J. Dellar.
\newblock An energy and potential enstrophy conserving numerical scheme for the
  multi-layer shallow water equations with complete {C}oriolis force.
\newblock {\em Journal of Computational Physics}, 313:99--120, 2016.

\bibitem{sun2020structure}
Z.~Sun and Y.~Xing.
\newblock On structure-preserving discontinuous {G}alerkin methods for
  {H}amiltonian partial differential equations: energy conservation and
  multi-symplecticity.
\newblock {\em Journal of Computational Physics}, 419:109662, 2020.

\bibitem{xu2008discontinuous}
Y.~Xu, J.~J. van~der Vegt, and O.~Bokhove.
\newblock Discontinuous {H}amiltonian finite element method for linear
  hyperbolic systems.
\newblock {\em Journal of scientific computing}, 35:241--265, 2008.

\bibitem{yildiz2023data}
S.~Yildiz, P.~Goyal, T.~Bendokat, and P.~Benner.
\newblock Data-driven identification of quadratic symplectic representations of
  nonlinear {H}amiltonian systems.
\newblock {\em arXiv preprint arXiv:2308.01084}, 2023.

\bibitem{zhao2020solving}
C.~L. Zhao.
\newblock Solving {A}llen-{C}ahn and {C}ahn-{H}illiard equations using the
  adaptive physics informed neural networks.
\newblock {\em Communications in Computational Physics}, 29(3), 2020.

\end{thebibliography}

\end{document}